\documentclass[twoside,12pt]{article}
\usepackage[all]{xy}
\usepackage{amssymb,amsmath,latexsym}
\usepackage{theorem}
\usepackage{amsfonts}
\usepackage{ulem}
\usepackage{amscd}
\usepackage{amsxtra}
\usepackage{mathrsfs}

\newcommand{\field}[1]{\mathbb{#1}}

\newcommand{\R}{\field{R}}
\newcommand{\Q }{\field{Q}}
\newcommand{\Z }{\field{Z}}
\newcommand{\N }{\field{N}} 

\newtheorem{thm}{Theorem}[section]
\newtheorem{cor}[thm]{Corollary}
\newtheorem{lem}[thm]{Lemma}
\newtheorem{prop}[thm]{Proposition}

\newtheorem{defn}[thm]{Definition}

\newtheorem{rem}[thm]{Remark}

\newtheorem{exm}[thm]{Example} 
\newtheorem{observ}[thm]{Observation} 

\newtheorem{conv}[thm]{Convention}

\numberwithin{equation}{section}

\newcommand{\cqfd}
{\hspace{1cm}
\rule{2mm}{2mm}%
\medbreak%
\par%
}

\def\pr{{\parindent0pt {\bf Proof.\ }}}
\def\cqfd
{\hspace{1cm}
\rule{2mm}{2mm}%
\medbreak%
\par%
}

\def\nil{{\rm Nil}}

\author{}

\begin{document}
\title{On $n$-trivial Extensions of Rings}

\date{}
\maketitle \vspace*{-1.5cm}

\thispagestyle{empty}


\begin{center}
\author{D.  D. Anderson$^1$, Driss Bennis$^{2,\textbf{a}}$, Brahim Fahid$^{2,\textbf{b}}$ and Abdulaziz Shaiea$^{2,\textbf{c}}$}\bigskip

\small{1: Department of Mathematics, The University of Iowa, 
Iowa City, IA 52242-1419, USA.}\end{center}\vspace{-0,3cm}
\hspace{4cm}  \small{dan-anderson@uiowa.edu}\smallskip
\begin{center}
\small{2:  Department of Mathematics, Faculty of Sciences,   Mohammed V University in Rabat, Morocco.}
\end{center}\vspace{-0,3cm}
\hspace{4cm} \small{\textbf{a}: d.bennis@fsr.ac.ma; driss$\_$bennis@hotmail.com}\\
$\mbox{}$\hspace{4cm} \small{\textbf{b}:  fahid.brahim@yahoo.fr}\\
$\mbox{}$\hspace{4cm} \small{\textbf{c}: abdulazizsha@yahoo.com; a.a.shaiea@gmail.com}

\bigskip\bigskip
\noindent{\large\bf Abstract.} The notion of  trivial extension of a ring by a module has been extensively studied and used in   ring theory  as well as in various other areas of research like   cohomology theory,   representation theory,   category theory and homological algebra.  In this paper we extend this classical ring construction by associating a ring to a ring $R$ and a family  $M=(M_i)_{i=1}^{n}$ of $n$ $R$-modules for a given integer $n\geq 1$.  We call this new ring construction an $n$-trivial extension of $R$ by $M$. In particular, the classical trivial extension will be just the $1$-trivial extension. Thus we generalize several known results on the classical trivial extension to the setting of $n$-trivial extensions and we give some new ones. Various ring-theoretic constructions and properties of $n$-trivial extensions are studied and a detailed investigation of the graded aspect of $n$-trivial extensions is also given. We end the paper with an investigation of various divisibily properties of $n$-trivial extensions. In this context several open questions arise.\bigskip

\small{\noindent{\bf 2010 Mathematics Subject Classification.}  primary 13A02, 13A05,  13A15, 13B99, 13E05, 13F05, 13F30; secondary  	16S99, 17A99.}

\small{\noindent{\bf Key Words.}  trivial extension; $n$-trivial extension; graded rings; homogeneous ideal.}

\section{Introduction}
 Except for a brief excursion in Section 2, all rings considered in this paper are assumed to be commutative with an identity; in particular, $R$
denotes such a ring, and all modules are assumed to be unitary left modules. Of course left-modules over a commutative ring $R$ are actually $R$-bimodules with $mr:=rm$.  Let $\Z$ (resp., $\N$) denotes the set of integers (resp., natural numbers). The set $\N \cup \{0\}$ will be denoted by $\N_0$. The ring $\Z/n\Z$ of the residues modulo an integer $n\in \N$ will be noted by $\Z_n$.\bigskip

Recall that the trivial extension of   $R$ by an $R$-module
$M$ is the ring denoted by $R\ltimes M$ whose underlying additive group is
$R\oplus M$ with multiplication given by $(r,m)(r',m')=(rr',rm'+mr')$.  Since its introduction by Nagata in \cite{MN},  the trivial extension of rings (also called idealization since it reduces questions about modules to ideals) has been used by many authors and in various contexts in order to produce examples of rings satisfying preassigned conditions (see, for instance, \cite{AW} and \cite{HH}).\medskip

 It is known that  the trivial extension $R\ltimes M$  is related to the following two ring constructions (see for instance \cite[Section 2]{AW}):\medskip

\textbf{Generalized triangular matrix ring.} Let $\mathscr{R}:=  (R_i)_{i=1}^{n}$ be a family of rings and $\mathscr{M}:= (M_{i,j})_{ 1\leq i < j \leq n }$ be a family of modules such that for each  $1\leq i < j \leq n$,  $M_{i,j}$ is an $ (R_i,R_j)$-bimodule. Assume for every $1\leq i < j <k \leq n$, there exists an $(R_i, R_k)$-bimodule homomorphism $$ M_{i,j} \otimes_{R_j }M_{j,k} \longrightarrow  M_{i,k}$$ denoted multiplicatively such that   $(m_{i,j} m_{j,k})m_{k,l}  =m_{i,j}( m_{j,k}m_{k,l})  $ for every  $ (m_{i,j}, m_{j,k},m_{k,l}   ) \in  M_{i,j} \times M_{j,k} \times M_{k,l}$. 
Then  the set 
$$  \begin{pmatrix}
                       R_{1 }\ & M_{1,2}\ &  \cdots \ &\cdots & M_{1,n-1} \;\quad M_{1,n}\\\
                       0& R_{ 2}\ &  \cdots\ &\cdots &M_{2,n-1} \;\quad M_{2,n}\\\
                       \vdots &\ddots &\ddots &\ddots &\vdots\;\quad \quad\vdots\\
                        \vdots &\ddots &\ddots&\ddots& \vdots\;\quad \quad\vdots\\
                       0&0& \cdots&0 & \quad R_{ n-1}\ \ M_{n-1,n}\ \\
                       0&0& \cdots&0 & \  \quad 0  \quad \quad R_{n }\  
\end{pmatrix}$$ 
consisting of matrices\\ 
 $  \begin{pmatrix}
                       m_{1,1}\ & m_{1,2}\ &  \cdots \ &\cdots & m_{1,n-1} \;\quad m_{1,n}\\\
                       0& m_{2,2}\ &  \cdots\ &\cdots &m_{2,n-1} \;\quad m_{2,n}\\\
                       \vdots &\ddots &\ddots &\ddots &\vdots\;\quad \ \quad\vdots\\
                        \vdots &\ddots &\ddots&\ddots& \vdots\;\quad \ \quad\vdots\\
                       0&0& \cdots&0 &   m_{n-1,n-1}\ \ m_{n-1,n}\ \\
                       0&0& \cdots&0 &   \quad 0\;  \quad \quad m_{n,n}\  
\end{pmatrix}$, 
   $m_{i,i} \in R_i$ and  $m_{i,j} \in M_{i,j}$ ($1\leq i < j \leq n$),\\
  with the usual matrix addition and multiplication is a ring called a \textit{generalized (or formal) triangular matrix} ring  and  denoted also by $T_n(\mathscr{R},\mathscr{M})$ (see \cite {GH}  and \cite{GJT}). 
 Then the trivial extension $R\ltimes M$ is naturally isomorphic to the subring of $ \begin{pmatrix}R&M\\
 0&R
\end{pmatrix} $ consisting of matrices $\begin{pmatrix} r&m\\
 0&r
\end{pmatrix} $ where $r\in R$ and $m\in M$ (note that, since $R$ is commutative $rm=mr$).\medskip

\textbf{Symmetric algebra.} Recall that the symmetric algebra associated to   $M$  is the graded ring quotient $S_R(M):=T_R(M)/H$ where $T_R(M)$ is the graded tensor $R$-algebra with $T^{n}_{R}(M)=M^{\otimes n}$ and $H$ is the homogeneous ideal of $T_R(M)$ generated by $ \{m\otimes n-n\otimes m|m,n\in M\}$. Note that  $S_R(M)=\overset{\infty}{\underset{n=0}{\oplus}} S^{n}_{R}(M)$ is a graded $R$-algebra with $S^{0}_{R}(M)=R$ and $S^{1}_{R}(M)=M$ and, in general, $S^{i}_{R}(M)$ is the image of $T^{i}_{R}(M)$ in $S_R(M)$. Then   $R\ltimes M$ and $S_R(M)/\underset{n\geq2}{\oplus} S^{n}_{R}(M)$ are  naturally isomorphic as graded $R$-algebras.\medskip

It is also worth recalling that when $M$ is a free $R$-module with a basis $B$, the   trivial extension $R\ltimes M$ is also naturally isomorphic to $ R[\{X_b\}_{b\in B}]/(\{X_b\}_{b\in B})^2$ where $\{X_b\}_{b\in B} $ is a set of indeterminates over $R$. In particular, $R\ltimes R \cong R[X]/(X^2)$.\medskip

Inspired by the facts above, we introduce an extension of the classical trivial extension of rings to   extensions associated to $n$ modules for  any integer  $n\geq 1$.\medskip

  In the literature,  particular cases of such extensions have been  used to solve some open questions.  In \cite{ANY} the authors introduced   an extension for $n=2$ and they  used it  to give a counterexample  of the so-called Faith conjecture. Also, in the case $n=2$, an extension is introduced in    \cite{CHN}     to give an example of a ring which has a non-self-injective injective hull with compatible multiplication. This gave a negative answer of a question posed by Osofsky. In \cite{Po} the author introduced and studied a particular extension  for the case $n=3$ to obtain a Galois coverings for the enveloping algebras of trivial extension algebras of triangular algebras.     Also, there is a master's thesis \cite{Mc} which introduced and studied factorization properties of an extension of the trivial extension of a ring by itself (i.e., self-idealization). In this paper, we introduce the following  extension ring construction for an arbitrary integer $n\geq 1$. \medskip

Let $M=(M_i)_{i=1}^{n}$ be a family of $R$-modules and  $ \varphi=\{\varphi_{i,j}\}_{ \underset{1\leq i,j\leq n-1}{i+j\leq n}}$ be a family of bilinear maps such that each $\varphi_{i,j}$ is  written multiplicatively:
 $$\begin{array}{rccl}
\varphi_{i,j}:& M_i\times M_j&\longrightarrow & M_{i+j}\\
               & (m_i,m_j) &\longmapsto & \varphi_{i,j}(m_i,m_j):= m_i m_j.
\end{array}$$ 
In particular, if all $M_i$ are submodules of the same $R$-algebra $L$, then the bilinear maps,   if they are not specified, are just the multiplication of $L$  (see examples in Section 2).  The $n$-$\varphi$-trivial extension of $R$ by  $M$    is the set denoted by $R\ltimes_{ \varphi} M_{1}\ltimes  \cdots \ltimes M_{n}$ or simply $R \ltimes_{ \varphi} M$  whose underlying additive group is
$R\oplus M_{1}\oplus  \cdots \oplus M_{n}$ with multiplication given by
$$(m_0,...,m_n) (m'_0,...,m'_n) =(\underset{j+k=i}{ \sum}  m_jm'_k)$$
for all $(m_i), \, (m'_i)\in R \ltimes_{\varphi} M$. We could also define the product $\varphi_{i,j}: M_i\times M_j \longrightarrow   M_{i+j}$  as an $R$-bimodule homomorphism   $\widetilde{\varphi}_{i,j}: M_i\otimes M_j \longrightarrow   M_{i+j}$; see Section 2 for details. For the sake of simplicity, it is convenient to  set $M_0=R$. In what follows, if no ambiguity arises, the $n$-$\varphi$-trivial extension of $R$ by  $M$ will be simply called an $n$-trivial extension of $R$ by  $M$ and denoted by  $R\ltimes_{n} M_{1}\ltimes  \cdots \ltimes M_{n}$ or simply $R \ltimes_{n} M$.\\
\indent  While in general $R \ltimes_{n} M$ need not to be a commutative ring, in Section 2, we give conditions on the maps $\varphi_{i,j}$ that force $R \ltimes_{n} M$ to be a ring. Unless otherwise stated, we assume the maps $\varphi_{i,j}$ have been defined so that $R \ltimes_{n} M$ is a commutative associative ring with identity. Thus   $R \ltimes_{n} M$ is a commutative ring with identity $(1,0,...,0)$. Moreover,  $R \ltimes_{n} M$  is  naturally isomorphic to  the subring of the generalized triangular matrix ring 
       $$ \begin{pmatrix}
                       R&M_1&M_2&\cdots\cdots&M_n\\
                       0&R&M_1&\cdots & M_{n-1}\\
                       \vdots&\ddots&\ddots&\ddots & \vdots\\
                        0 &0&0&\cdots&  M_{1}\\
                       0&0&0&\cdots  &    R 
\end{pmatrix} $$
 consisting of matrices\footnote{When $R$ is a field and $M_i=R$  for every $i\in\{1,...,n\}$, these matrices are   well-known as upper triangular Toeplitz matrices. In \cite{Mc}, the author used the same terminology for such matrices with entries in  a commutative ring.} 
 $$ \begin{pmatrix}
                      r&m_1&m_2&\cdots\cdots&m_n\\
                       0&r&m_1&\cdots & m_{n-1}\\
                       \vdots&\ddots&\ddots&\ddots & \vdots\\
                        0 &0&0&\cdots&  m_{1}\\
                       0&0&0&\cdots  &    r 
\end{pmatrix} $$
where $r\in R$ and $m_i\in M_i$ for every $i\in \{1,...,n\}$.\smallskip

When, for every $k\in \{1,..., n\}$, $M_k= S^{k}_{R}(M_1)$, the ring $R \ltimes_{n} M$  is  naturally isomorphic to  $ S_R(M_1)/\bigoplus_{k\geq n+1}S^{k}_{R}(M_1)$. 
 In particular, if $M_1=F$ is a free $R$-module with a basis $B$, then the   $n$-trivial extension   $R\ltimes F\ltimes  S^{2}_{R}(F)\ltimes \cdots \ltimes S^{n}_{R}(F)$ is also naturally isomorphic to $ R[\{X_b\}_{b\in B}]/(\{X_b\}_{b\in B})^{n+1}$  where $\{X_b\}_{b\in B} $ is a set of indeterminates over $R$. Namely, when $F\cong R$,  $$ R\ltimes_n R\ltimes   \cdots \ltimes R \cong R[X]/(X^{n+1}).$$

Also, in \cite{BAS},  the trivial extension of a ring $R$ by an ideal $I$ is   connected to   the  Rees algebra $\mathcal{R}_+ $ associated
to $ R$ and $I$ which is precisely the following graded subring of $R[t]$ (where $t$ is an indeterminate over $R$):
$$\mathcal{R}_+:=  \bigoplus_{n\geq 0}I^nt^n \, .$$
Using   \cite[Lemma 1.2 and Proposition 1.3]{BAS}, we get, similar to  \cite[Proposition 1.4]{BAS}, the following diagram of extensions and isomorphisms of rings:
$$\xymatrix{
   R \ar@{=}[d] \ar@{^{(}->}[r] &  \mathcal{R}_+/(I^{n+1}t^{n+1}) \ar@{^{(}->}[r] \ar[d]^{\cong}   & R[t]/(t^{n+1}) \ar[d]^{\cong}  \\
   R  \ar@{^{(}->}[r] &    R\ltimes_n I\ltimes  I^2\ltimes \cdots \ltimes I^{n} \ar@{^{(}->}[r]  & R\ltimes_n R\ltimes   \cdots \ltimes R
 }$$
 \medskip

In this paper,  we study some   properties of the ring $R\ltimes_n M$, extending well-known results on the classical trivial  extension of rings. The paper is organized as follows.\medskip

In Section 2, we carefully define the $n$-trivial extension  $R\ltimes_n M$ giving conditions on the maps $\varphi_{i,j}$ so that  $R\ltimes_n M$ is actually a commutative ring with identity. Actually, we investigate the situation in greater generality where $R$ is not assumed to be commutative and $M_i$ is an $R$-bimodule for $i=1,...,n$. We end the section with a number of examples.\medskip 

In Section 3, we investigate some ring-theoretic constructions of $n$-trivial extensions. We begin by showing that $R\ltimes_n M$ may be considered as a graded ring for three different grading monoids, in particular, $R\ltimes_n M$ may be considered as $\N_0$-graded ring or $\Z_{n+1}$-graded ring. We then show how $R\ltimes_n M$ behaves with respect to polynomials  (Corollary \ref{cor-poly})  and power series   (Theorem \ref{thm-power}) extensions and localization (Theorem \ref{thm-mult}). In Theorem \ref{thm-prod}, we show that the  $n$-trivial extension of a finite direct product of rings is a finite direct product of $n$-trivial extensions. We end with two results on inverse limits and direct limits of $n$-trivial extensions (Theorems \ref{thm-inverse-limit} and  \ref{thm-direct-limit}).\medskip

In Section 4, we  present some natural ring homomorphisms   related to  $n$-trivial extensions (see Proposition \ref{pro-2-n-exte}). Also,  we study some basic properties of $R\ltimes_n M$. Namely, we extend the characterization of   prime and   maximal ideals of the classical trivial extension to  $R\ltimes_n M$ (see Theorem \ref{thm-prime}). As a consequence,    the nilradical and the Jacobson radical are determined (see Corollary \ref{cor-prime}). Finally, as an extension of \cite[Theorems 3.5 and 3.7]{AW},  the  set of  zero divisors, the set of units  and the set of idempotents of    $R\ltimes_n M$ are also   characterized (see Proposition \ref{pro-zero-U-Id}).\medskip

In Section 5, we investigate the graded aspect of $n$-trivial extensions. The motivation behind this study is that, in  the classical case (where $n=1$), the study of trivial extensions  as $\Z_2$-graded rings has lead  to some interesting properties (see  \cite{AW}) and has shed more light on the structure of ideals of the trivial extensions. In Section 5, we extend some of results given in \cite{AW} and we give some new ones. Namely, among other results, we characterize the homogeneous ideals of $R\ltimes_n M$  (Theorem \ref{thm1-homoge}) and we investigate some of their properties  (Propositions \ref{pro-homogeneous} and \ref{pro-homogeneous-fg}). We devote the remainder of Section 5 to investigate the question ``When is every ideal of a given class $ \mathscr{I} $  of ideals of $R\ltimes_n M$    homogeneous?" (see the discussion after Proposition  \ref{pro-homogeneous-fg}). In this context various results and examples are established.\medskip

Section 6 is devoted to some classical  ring-theoretic properties.  Namely, we  characterize when $R\ltimes_n M$ is, respectively,   Noetherian, Artinian,   (Manis) valuation, Pr\"{u}fer, chained, arithmetical,  a $\pi$-ring, a  generalized ZPI-ring or a PIR. We end the section with a remark on a question posed in \cite{A86} concerning $m$-Boolean rings.\medskip

Finally, in Section 7 we study divisibility properties of $n$-trivial extensions. We are mainly interested in showing how one could extend results on the classical trivial extension  presented in \cite[Section 5]{AW} to the context of $n$-trivial extensions.
\medskip


\section{The general $n$-trivial extension construction and some examples}
The purpose of this section is to formally define the $n$-trivial extension ($n\geq 1$) $ R\ltimes_n M_1\ltimes \cdots \ltimes M_n$ where $R$ is a commutative ring with identity and each $M_i$ is an $R$-module, and to give some interesting examples of $n$-trivial extensions. However, to better understand the construction and the underlying multiplication maps $\varphi_{i,j}:  M_i\times M_j \longrightarrow   M_{i+j}$, we begin in the more general context of $R$ being an associative ring (not necessarily commutative) with identity and the $M_i$'s being $R$-bimodules. Also, as there is a significant difference in the cases $n=1$, $n=2$ and $n\geq 3$, we handle these three cases separately.\medskip

Let $R$ be an associative ring with identity and $M_1$,...,$M_n$ be unitary $R$-bimodules (in the case where $R$ is commutative we will always assume that $rm=mr$ unless stated otherwise).\medskip

\noindent\underline{Case $n=1$.}\\
 $ R\ltimes_1 M_1=  R\ltimes  M_1 =   R\oplus M_1$ is just the trivial extension with multiplication $ (r,m ) (r',m' )= (rr',rm'+mr')$. Here $ R\ltimes_1 M_1$ is an associative ring with identity where the associative and distributive laws follow from the ring and $R$-bimodule axioms. For $R$ commutative, we write $ (r,m ) (r',m' )= (rr',rm'+r'm)$ as $r'm=mr'$. Now 
$ R\ltimes_1 M_1$ is an $\N_0$-graded or a $\Z_2$-graded ring isomorphic to $T_R(M_1)/\underset{i\geq 2}{\oplus} T^{i}_{R}(M_1)$ or $S_R(M_1)/\underset{i\geq 2}{\oplus}  S^{i}_{R}(M_1)$ and to the matrix ring representation mentioned in the introduction. Note that we could drop the assumption that $R$ has an identity and $M_1$ is unitary. We then get that $ R\ltimes_1 M_1$ has an identity (namely $(1,0)$) if and only if $R$ has an identity and $M_1$ is unitary.\medskip 

\noindent\underline{Case $n=2$.}\\ Here $ R\ltimes_2 M_1\ltimes M_2=     R\oplus M_1 \oplus M_2$ with coordinate-wise addition and multiplication $$ (r,m_1,m_2 ) (r',m'_1,m'_2 )= (rr',rm'_1+m_1r' , rm'_2+m_1m'_1+m_2r')$$ where $m_1m'_1:=\varphi_{1,1}(m_1,m'_1)$ with the map $\varphi_{1,1}:  M_1\times M_1 \longrightarrow M_{2}$. We readily see that $ R\ltimes_2 M_1\ltimes M_2$ satisfying the distributive laws is equivalent to $\varphi_{1,1}$ being additive in each coordinate. Since $R$ is assumed to be associative and $M_1$ and $M_2$ to be $R$-bimodules,  $ R\ltimes_2 M_1\ltimes M_2$ is associative precisely when $(rm_1)m'_1=r(m_1m'_1)$, $(m_1r)m'_1=m_1(rm'_1)$, and $(m_1 m'_1)r=m_1( m'_1r)$ for $r\in R$ and $m_1, m'_1\in M_1$. This is equivalent to $\varphi_{1,1}(rm_1,m'_1)=r\varphi_{1,1}(m_1,m'_1)$, $\varphi_{1,1}(m_1r, m'_1)=\varphi_{1,1}(m_1,rm'_1)$, and $\varphi_{1,1}(m_1 ,m'_1)r=\varphi_{1,1}(m_1, m'_1r)$. For $R$-bimodules $M$, $N$ and $L$, we call  a function $f :  M\times N \longrightarrow L$ a \textit{pre-product map} if it is additive in each coordinate, is middle linear (i.e.,  $f(m r, m' )=f(m ,rm')$) and is left and right homogeneous (i.e.,   $f(rm,m')=rf(m,m')$ and $f(m,m'r)=f(m,m')r$). Note that a pre-product map $f:M\times N\longrightarrow L$ uniquely corresponds to an $R$-bimodule homomorphism $\tilde{f}:M\otimes_R N\longrightarrow L$ with $f(m,n)=\tilde{f}(m\otimes n)$. Thus a  pre-product map $\varphi_{1,1}:M_1\times M_1\longrightarrow M_2$ corresponds to an $R$-bimodule homomorphism $\tilde{\varphi}_{1,1}:M_1\otimes_R M_1\longrightarrow M_2$. So we could equivalently define $m_1m'_1:=\tilde{\varphi}_{1,1}(m_1\otimes m'_1)$.\medskip

 So  $ R\ltimes_2 M_1\ltimes M_2$ is an (associative) ring with identity precisely when $\varphi_{1,1}$ is a pre-product map or 
$\tilde{\varphi}_{1,1}:M_1\otimes_R M_1\longrightarrow M_2$ is an $R$-bimodule homomorphism. We can identify  $ R\ltimes_2 M_1\ltimes M_2$ with the matrix representation given in the introduction:  $ (r,m_1,m_2 )$ is identified with    $$
  \begin{pmatrix}
                       r&m_1&m_2\\
                       0&r&m_1  \\ 
                     0 &  0&  r \mbox{}
\end{pmatrix} $$ 
 But the relationship with a tensor algebra or symmetric algebra is more difficult. When $R\ltimes_2 M_1\ltimes M_2$ is an associative ring, we can define a ring epimorphism  
$$ 
  T_R(M_1\oplus M_2)/ \underset{ i\geq 3}{ \oplus } T^{i}_R(M_1\oplus M_2) \longrightarrow   R\ltimes_2 M_1\ltimes M_2$$
by
      $$(r,(m_1,m_2),       \overset{l}{\underset{i=1}{\sum}} (m_{1,i},m_{2,i})\otimes (m'_{1,i},m'_{2,i}))+ \underset{ i\geq 3}{ \oplus }  T^{i}_R(M_1\oplus M_2) \longmapsto (r,m_1,m_2+\overset{l}{\underset{i=1}{\sum}} m_{1,i}m'_{1,i}). $$
  For the commutative case,  we get a similar ring epimorphism $S_R(M_1\oplus M_2)/\underset{ i\geq 3}{ \oplus }  S^{i}_R(M_1\oplus M_2)\longrightarrow R\ltimes_2 M_1\ltimes M_2$.\medskip
 
For $ R\ltimes_2 M_1\ltimes M_2$ to be a commutative ring with identity we need $R$ to be commutative with identity and $m_1 m'_1=m'_1m_1 $ for $m_1, m'_1\in M_1$, or 
 $\varphi_{1,1}(m_1,m'_1)=\varphi_{1,1}(m'_1,m_1) $. Thus for $R$ commutative, $ R\ltimes_2 M_1\ltimes M_2$ is a commutative ring if and only if 
$\varphi_{1,1}$ is a symmetric $R$-bilinear map, or equivalently, $\tilde{\varphi}_{1,1}(m_1\otimes m'_1)=\tilde{\varphi}_{1,1}(m'_1\otimes m_1)$.\medskip

\noindent\underline{Case $n\geq 3$.}\\
 Here again $R$ is an associative ring with identity and  $   M_1,$...,$ M_n$  ($n\geq 3$) are   $R$-bimodules. So  $ R\ltimes_n M_1\ltimes \cdots \ltimes M_n =R\oplus  M_1 \oplus \cdots \oplus M_n$ with coordinate-wise addition. Assume we have pre-product maps $\varphi_{i,j}:  M_i\times M_j \longrightarrow    M_{i+j}$,  or equivalently, the corresponding  $R$-bimodule homomorphism $\tilde{\varphi}_{i,j}:M_i\otimes_R M_j\longrightarrow M_{i+j}$
 for $1\leq i,j\leq n-1$ with $i+j\leq n$. As usual set $$m_i  m_j:=\varphi_{i,j}(m_i, m_j)=\tilde{\varphi}_{i,j}(m_i\otimes m_j)$$ for $m_i\in M_i$ and   $m_j\in M_j$. Setting $R=M_0$, we can write the multiplication in  $R\ltimes_n M_1\ltimes \cdots \ltimes M_n$ as  $(m_0,...,m_n)(m'_0,...,m'_n)=(m''_0,...,m''_n) $ where $m''_i=\underset{j+k=i}{ \sum}   m_jm'_k$. Then 
$ R\ltimes_n M_1\ltimes \cdots \ltimes M_n$  satisfies the distributive laws because the maps $\varphi_{i,j}$ are additive in each coordinate. So $ R\ltimes_n M_1\ltimes \cdots \ltimes M_n$ is a not necessarily associative ring with identity $(1,0,...,0)$ (see Example \ref{exmp-2} for a case where $ R\ltimes_n M_1\ltimes \cdots \ltimes M_n$ is not associative). Note that $ R\ltimes_n M_1\ltimes \cdots \ltimes M_n$  is associative precisely when $(m_im_j)m_k = m_i(m_jm_k)$ for $m_i\in M_i$, $m_j\in M_j$  and $m_k\in M_k$ with $1\leq i,j,k\leq n-2$ and $i+j+k\leq n$. In terms of the pre-product maps, this says that $\varphi_{i+j,k}(\varphi_{i,j}(m_i,m_j),m_k) =\varphi_{i,j+k}(m_i,\varphi_{j,k}(m_j,m_k))$, or equivalently,  $$\widetilde{\varphi}_{i+j,k} \circ (\widetilde{\varphi}_{i,j} \otimes id_{M_k}) =\widetilde{\varphi}_{i,j+k}\circ( id_{M_i} \otimes \widetilde{\varphi}_{j,k})$$ where  $id_{M_l}$ is the identity map on $M_l$ for $l \in \{1,..., n\}$. In other words, the diagram below commutes:
$$\xymatrixcolsep{6pc}\xymatrix{
M_{i}\otimes M_j \otimes  M_k \ar[r(0.8)]^-{id_{M_i} \otimes \widetilde{\varphi}_{j,k}}  \ar[d]_{\widetilde{\varphi}_{i,j} \otimes id_{M_k} }  & M_{i} \otimes  M_{j+k}  \ar[d]^{\widetilde{\varphi}_{i,j+k} } \\ 
  M_{i+j} \otimes  M_k  \ar[r(0.8)]_-{ \widetilde{\varphi}_{i+j,k}  }  &  M_{i+j+k}
}$$
   Let us call a family $\{\varphi_{i,j}\}_{ \underset{1\leq i,j\leq n-1}{i+j\leq n}}$ (or $\{\widetilde{\varphi}_{i,j}\}_{ \underset{1\leq i,j\leq n-1}{i+j\leq n}}$) of pre-product maps satisfying the previously stated associativity condition a \textit{family of product maps}. So, when 
$\{\varphi_{i,j}\}_{ \underset{1\leq i,j\leq n-1}{i+j\leq n}}$ (or equivalently $\{\widetilde{\varphi}_{i,j}\}_{ \underset{1\leq i,j\leq n-1}{i+j\leq n}}$) is a family of product maps, $R\ltimes_n M_1\ltimes \cdots \ltimes M_n$  is an associative  ring with identity. Further, 
  for $ R\ltimes_n M_1\ltimes \cdots \ltimes M_n$  to be a commutative ring with identity we need $R$ to be commutative with identity and  $\varphi_{i,j}(m_i,m_j)=\varphi_{j,i}(m_j,m_i)$ for every $1\leq i,j\leq n-1$ with $i+j\leq n$, or equivalently, $\widetilde{\varphi}_{i,j}=\widetilde{\varphi}_{j,i} \circ \tau_{i,j}$ where $\tau_{i,j} : M_{i}\otimes M_j  \rightarrow M_{j}\otimes M_i $ is the `flip' map defined by $\tau_{i,j}(m_i \otimes m_j)=m_j \otimes m_i$ for every $m_i \otimes m_j\in M_{i}\otimes M_j$. In other words, the diagram below commutes:
$$ \xymatrix{
M_{i}\otimes M_j   \ar[r]^-{ \widetilde{\varphi}_{i,j}  }  \ar[d]_{ \tau_{i,j}  }  &   M_{i+j}  \\ 
  M_{j} \otimes  M_i \ar[ur]_-{ \widetilde{\varphi}_{j,i}  }  
}$$ 
 In this case, the family $\{\varphi_{i,j}\}_{ \underset{1\leq i,j\leq n-1}{i+j\leq n}}$ (or $\{\widetilde{\varphi}_{i,j}\}_{ \underset{1\leq i,j\leq n-1}{i+j\leq n}}$) will be called a \textit{family of commutative product maps}. So, when $R$ is commutative and $\{\varphi_{i,j}\}_{ \underset{1\leq i,j\leq n-1}{i+j\leq n}}$ (or equivalently $\{\widetilde{\varphi}_{i,j}\}_{ \underset{1\leq i,j\leq n-1}{i+j\leq n}}$) is   a family of commutative product maps, 
$ R\ltimes_n M_1\ltimes \cdots \ltimes M_n$ is a commutative ring with identity.\medskip

As in the case $n=2$, when  $R\ltimes_n M_1\ltimes \cdots\ltimes M_n$ is an   (associative) ring with identity, we can identify  $R\ltimes_n M_1\ltimes \cdots\ltimes M_n$ with the matrix representation given in the introduction:  $ (r,m_1,...,m_n )$ is identified with    
 $$ \begin{pmatrix}
                      r&m_1&m_2&\cdots\cdots&m_n\\
                       0&r&m_1&\cdots & m_{n-1}\\
                       \vdots&\ddots&\ddots&\ddots & \vdots\\
                        0 &0&0&\cdots&  m_{1}\\
                       0&0&0&\cdots  &    r 
\end{pmatrix} $$
\indent Also, as in the case $n=2$,   when $R\ltimes_n M_1\ltimes \cdots\ltimes M_n$  is an associative ring, we can define a ring epimorphism  
  $T_R(M_1\oplus \cdots \oplus M_n)/   \underset{i\geq n+1}{\oplus}T^{i}_R(M_1\oplus\cdots \oplus M_n)\longrightarrow R\ltimes_n M_1\ltimes \cdots \ltimes M_n$ and we have a similar result concerning the symmetric algebra when $R\ltimes_n M_1\ltimes \cdots \ltimes M_n$ is commutative.\medskip

\begin{rem}
\begin{enumerate}
    \item Let $R_1$  and $R_2$ be two rings and $H$ an  $( R_1, R_2)$-bimodule.  It is well-known that  every  generalized  triangular matrix ring is naturally isomorphic to the trivial extension of $ R_1\times R_2$ by $H$ where the actions  of   $ R_1\times R_2$ on $H$ are defined as follows:
  $(r_1,r_2)h=r_1 h$ and   $h(r_1,r_2)= h r_2 $ for every $(r_1,r_2)\in   R_1\times R_2$ and $h\in H$. Below we see that an observation on the product of two matrices of the generalized triangular matrix ring shows that this fact can be extended to   $n$-trivial extensions.\\
\indent  Consider the generalized triangular matrix ring  $$ T_n(\mathscr{R},\mathscr{M})=  \begin{pmatrix}
                       R_{1 }\ & M_{1,2}\ &  \cdots \ &\cdots & M_{1,n-1} \;\quad M_{1,n}\\\
                       0& R_{ 2}\ &  \cdots\ &\cdots &M_{2,n-1} \;\quad M_{2,n}\\\
                       \vdots &\ddots &\ddots &\ddots &\vdots\;\quad \quad\vdots\\
                        \vdots &\ddots &\ddots&\ddots& \vdots\;\quad \quad\vdots\\
                       0&0& \cdots&0 & \quad R_{ n-1}\ \ M_{n-1,n}\ \\
                       0&0& \cdots&0 & \  \quad 0  \quad \quad R_{n }\  
\end{pmatrix}$$
where $(R_i)_{i=1}^{n}$ is a family of rings and $(M_{i,j})_{ 1\leq i < j \leq n }$ is a family of modules such that  for each  $1\leq i < j \leq n$,  $M_{i,j}$ is an $ (R_i,R_j)$-bimodule. Assume, for every $1\leq i < j <k \leq n$, there exists an $(R_i, R_k)$-bimodule homomorphism $$ M_{i,j} \otimes_{R_j }M_{j,k} \longrightarrow  M_{i,k}$$ denoted multiplicatively such that   $$(m_{i,j} m_{j,k})m_{k,l}  =m_{i,j}( m_{j,k}m_{k,l}) $$ for every  $ (m_{i,j}, m_{j,k},m_{k,l}) \in  M_{i,j} \times M_{j,k} \times M_{k,l}$.\\
\indent Consider  the finite direct product of rings $R=R_1\times \cdots \times R_n$ and set, for $2\leq i \leq n$, $M_i=M_{1,i}\times  M_{2,i+1}\times \cdots \times M_{n-(i-1),n}$ (for $i=n$,  $M_n=M_{1,n}$). We need to define an action of $R$ on each $M_i$ and a family of product maps so  that $R\ltimes_{n-1} M_2\ltimes \cdots \ltimes M_n$  is an $n-1$-trivial extension isomorphic to $T_n(\mathscr{R},\mathscr{M})$.\\
\indent  First, note that, for every matrix $A=(a_{i,j})$ of $T_n(\mathscr{R},\mathscr{M})  $ and for every $2\leq i \leq n$, the $i$-th diagonal above the main diagonal of $A$ naturally corresponds to the following $(n-i+1)$-tuple  $(a_{1,i}, a_{2,i+1},. ..,a_{n-(i-1),n})$ of $M_i$.  On the other hand, consider two matrices $A=(a_{i,j})$ and $B=(b_{i,j})$ of $T_n(\mathscr{R},\mathscr{M})$, and denote  the product $AB$ by $C=(c_{i,j})$. Then using the above correspondence for   $2\leq i \leq n$, the $i$-th diagonal above the main diagonal of $C$ can be seen as the  $(n-i+1)$-tuple $c_i=(c_{j,i+j-1})_j \in M_i$ such that, for every $1 \leq j \leq n-i+1$,
 $$\begin{array}{ccl}
     c_{j,i+j-1} &= &   \overset{i+j-1 }{\underset{k=j}{\sum}} a_{ j,k} b_{k,i+j-1 }\\
                 &=&   \overset{ i }{\underset{k=1}{\sum}}  a_{ j,k+j-1} b_{k+j-1,i+j-1 }\, .
\end{array}
$$
Then 
 $$\begin{array}{ccl}
     c_i &= & (\overset{ i }{\underset{k=1}{\sum}}  a_{ j,k+j-1} b_{k+j-1,i+j-1 })_j  \\
                 &=&  \overset{ i }{\underset{k=1}{\sum}}  ( a_{ j,k+j-1} b_{k+j-1,i+j-1 })_j\, .
\end{array}
$$
Thus the cases $k=1$ and $k=i$ allow us to define the left and right actions of $R$ on $M_i$ as follows: 
For every $(r_l)_l\in R$ and $(m_{j,i+j-1 })_j\in M_i$,
$$( r_l)_l (m_{j,i+j-1 })_j:= (r_j m_{ j,i+j-1 })_j$$
and 
$$( m_{ j,i+j-1})_j ( r_l)_l := ( m_{ j,i+j-1} r_{i+j-1})_j.$$
The other cases of $k$ can be used to define the  product maps $M_k\times M_{i-k} \longrightarrow  M_{i}$ as follows: Fix $k$, $ 1<k<i$, and consider $e_k=( e_{ j,k+j-1})_{1\leq j  \leq n-k+1}  \in M_k$ and $f_{i-k}= (f_{j,i-k+j-1})_{1\leq j  \leq n-i+k+1} \in M_{i-k}$. Then $$e_k f_{i-k}:=  ( e_{ j,k+j-1} f_{k+j-1,i+j-1 })_{1\leq j  \leq n-i+1}.$$
Therefore, endowed with these products, $R\ltimes_{n-1} M_2\ltimes \cdots \ltimes M_n$  is   an $n-1$-trivial extension naturally isomorphic to the generalized triangular matrix ring  $T_n(\mathscr{R},\mathscr{M})$.
 \item It is known that the generalized triangular matrix ring  $T_n(\mathscr{R},\mathscr{M})$ can be seen as a  generalized triangular $2 \times 2$  matrix ring\footnote{We are indebted to J. R. Garc\'{\i}a Rozas (Universidad de Almer\'{i}a, Spain) who pointed out this remark.}. Namely, there is a natural ring isomorphism between  $T_n(\mathscr{R},\mathscr{M})$ and $ T_2(S,N)$ where $$S= (T_{n-1}((R_i)_{i=1}^{n-1},(M_{i,j})_{ 1\leq i < j \leq n-1 }),R_n)$$ and 
$$N= \begin{pmatrix}
                         M_{1,n}\\ 
                       M_{2,n}\\ 
                       \vdots\\ 
               M_{n-1,n}  
\end{pmatrix}$$
However, an $n$-trivial extension  is not necessarily a $1$-trival extension. For that consider, for instance, the $2$-trivial extension $S=\Z/2\Z\ltimes_2 \Z/2\Z \ltimes \Z/2\Z$. One can check easily that $S$ cannot be isomorphic to any $1$-trival extension. 
\end{enumerate} 
\end{rem}

We end this section with a number of examples.\medskip

\begin{exm}\label{exmp-2}  
Suppose that $R$ is  a commutative ring and consider $ R\ltimes_n R\ltimes   \cdots \ltimes R $ ($n\geq 1$) with  a family of product maps $\varphi_{i,j}: Re_i\times Re_j \longrightarrow Re_{i+j}$ where, for $k\in \{ 1,...,n\}$,  $ e_k=(0,...,0,1,0,...,0)$   with $1$ in the $k+1$'th place.\smallskip

\indent For $n=1$, $ R\ltimes_1 R\cong  R[X]/(X^{2})$. \smallskip

\indent Suppose, $n=2$ and $e_ 1^2=r_{1,1}e_2$. Then  $ R\ltimes_2 R \ltimes R  \cong  R[X,Y]/(X^2-r_{1,1}Y,XY,Y^2)$ where $X$ and $Y$ are commuting indeterminates. So in the case where $r_{1,1}=1$, we get 
 $ R\ltimes_2 R \ltimes R  \cong  R[X,Y]/(X^2-Y,XY,Y^2)\cong   R[X]/(X^{3})  $.\smallskip

\indent The case $n=3$ is more interesting. Now, for $1\leq i,j\leq 2$ with $i+j\leq 3$,  $\varphi_{i,j}: R \times R \longrightarrow R $ with $\varphi_{i,j}(r,s)=r\varphi_{i,j}(1,1)s$. Put $\varphi_{i,j}(1,1)=r_{i,j}$; so $(re_i)(se_j)=rr_{i,j}se_{i+j}$. Now,   $ R\ltimes_3 R\ltimes R \ltimes R $ is commutative if and only if $e_1e_2=e_2e_1$ or $r_{1,2}=r_{2,1}$. And $ R\ltimes_3 R\ltimes R \ltimes R $ is associative if and only if   $(e_1e_1)e_1=e_1(e_1e_1)$ or $ r_{1,1}r_{2,1}=r_{1,2}r_{1,1}$. Thus if $ R\ltimes_3 R\ltimes R \ltimes R $ is commutative, it is also associative. However, if $R$ is a commutative integral domain and $r_{1,1}\not = 0$, $ R\ltimes_3 R\ltimes R \ltimes R $ is associative if and only if it is commutative. Thus if we take $R=\Z$, $r_{1,1}=1$, $r_{1,2}=1$ and $r_{2,1}=2$, $ R\ltimes_3 R\ltimes R \ltimes R $  is a non-commutative and non-associative ring.\smallskip

\indent For $n=4$, the reader can easily check that $ R\ltimes_4 R\ltimes R \ltimes R \ltimes R $ is commutative if and only if  $r_{i,j} =r_{j,i}$ for   $1\leq i,j\leq 3$ with $i+j\leq 4$, and that  $ R\ltimes_4 R\ltimes R \ltimes R \ltimes R $ is associative if and only if $ r_{1,1}r_{2,1}=r_{1,2}r_{1,1}$, $r_{2,1}r_{3,1}=r_{2,2}r_{1,1}$, $ r_{1,2}r_{3,1}=r_{1,3}r_{2,1}$, and $ r_{1,1}r_{2,2}=r_{1,3}r_{1,2}$. Thus if $R$ is a commutative integral domain with  $ r_{1,1}\not = 0$, then $ r_{1,1}r_{2,1}=r_{1,2}r_{1,1}$ if and only if $r_{2,1}=r_{1,2}$. So if  $ r_{1,1}\not = 0$ and  $ r_{1,2}\not = 0$, then 
 $ r_{1,2}r_{3,1}=r_{1,3}r_{2,1}$ if and only if $r_{3,1}=r_{1,3}$. Thus if  $ r_{1,1}\not = 0$ and  $ r_{1,2}\not = 0$, then $ R\ltimes_4 R\ltimes R \ltimes R \ltimes R $ is associative forces  $ R\ltimes_4 R\ltimes R \ltimes R \ltimes R $  to be commutative and in this case  $ R\ltimes_4 R\ltimes R \ltimes R \ltimes R $ is associative if and only if $r_{1,1}r_{2,2}=r_{1,3}r_{1,2}$. Thus if three of the numbers $ r_{1,1}$, $r_{2,2}$, $r_{1,3}$ and $r_{1,2}$ are given and nonzero, then there is only one possible choice for the remaining $r_{i,j}$ for   $ R\ltimes_4 R\ltimes R \ltimes R \ltimes R $  to be associative. If we take $R=\Z$ and $ r_{1,1}=1$, $r_{2,1} =r_{1,2}=2$, $r_{2,2}=3$ and   $r_{1,3}=r_{3,1}=4$, then the resulting ring is commutative but not associative. \smallskip

\indent For $n\geq 5$, the reader can easily write conditions on the $r_{i,j}=\varphi_{i,j}(1,1)$ for $ R\ltimes_n R\ltimes   \cdots \ltimes R $ to be commutative or associative.
 \end{exm}

\begin{exm}\label{exmp-3} Let $R$ be  a commutative ring and 
  $N_1$,...,$N_n$ be $R$-submodules of an $R$-algebra $T$ with $N_iN_j \subseteq N_{i+j}$ for    $1\leq i,j\leq n-1$ with $i+j\leq n$. Then, using the multiplication from $T$,  $ R\ltimes_n N_1\ltimes    \cdots \ltimes N_n$ is a ring which is commutative if $T$ is commutative. The following are some interesting special cases:
\begin{itemize}
    \item[(a)] Let $R$ be a commutative ring and $I$   an ideal of $R$. Then $R\ltimes_n I\ltimes  I^2\ltimes \cdots \ltimes I^{n}$ is  the  quotient of the  Rees ring $R[It]/(I^{n+1}t^{n+1})$ mentioned  in the introduction.
    \item[(b)] Let $R$ be a commutative ring, $T$ an $R$-algebra and   $J_1\subseteq\cdots\subseteq J_n$   ideals of $T$. Then $R\ltimes_n J_1\ltimes    \cdots \ltimes J_n$ is an example of $n$-trivial extension since $J_iJ_j\subseteq J_i \subseteq J_{i+j}$ for $i+j\leq n$. For example, we could take   $R\ltimes_2 XR[X] \ltimes R[X]$.
    \item[(c)]  Suppose that $R_1\subseteq \cdots\subseteq R_n$ are     $R$-algebras where $R$ is a commutative ring. Let  $N$ be an $R_{n-1}$-submodule  of $R_n$ (in particular, we could take $N=R_n$). Then $R\ltimes_n R_1\ltimes    \cdots\ltimes R_{n-1} \ltimes N$ with the multiplication induced by $R_n$ is a ring. For example, we could take   $\Z \ltimes_3 \Q \ltimes \R \ltimes  N$  where $N$ is the $\R$-submodule of $\R[X]$ of polynomials of degree $\leq 5$.
\end{itemize} 
\end{exm}

\begin{exm}\label{exmp-4} Let $R$ be  a commutative ring and  $M$ an $R$-module. Let $S:= R\ltimes_n R\ltimes    \cdots\ltimes R\ltimes M$ with 
 $\varphi_{i,j}: R \times R \longrightarrow R $   the usual ring product in $R$ for $i+j\leq n-1$, but, for $i+j=n$ and $i,j\geq 1$, $\varphi_{i,j}$ is the zero map. So 
$$(r_0,...,r_{n-1},m_n)(r'_0,...,r'_{n-1},m'_n)=(r_0r'_0,r_0r'_1+r_1r'_0,...,r_0r'_{n-1}+\cdots+r_{n-1}r'_0,r_0m'_n+r'_0m_n). $$
Then $S\cong \R[X]/(X^n)\ltimes M$ where $M$ is considered as an $ \R[X]/(X^n)$-module with $\overline{f(X)}m=f(0)m$.
\end{exm}

\begin{exm}\label{exmp-5} Let $R$ be  a commutative ring and  $T$ an $R$-algebra. Let $J_1\subseteq \cdots\subseteq J_{n}$  be ideals of $T$. Then take $R\ltimes_n T/J_1\ltimes    \cdots \ltimes T/J_n$    where the product  $T/J_i\times T/J_j \longrightarrow T/J_{i+j} $ is given by  $(t_i +J_i) (t_j +J_j) = t_it_j +J_{i+j}$ for  $i+j\leq n$.
\end{exm}

\begin{exm}\label{exmp-6} Let $R$ be a commutative ring, $N_1$,..., $N_{n-1}$ ideals of $R$ and $N_n=Ra$ a cyclic $R$-module. Then consider 
$R\ltimes_n N_1\ltimes    \cdots \ltimes N_n$  where the products  $N_i\times N_j \longrightarrow N_{i+j} $  are the usual products for  $R$ when  $i+j\leq n-1$,  and for $i+j=n$ define     $ n_in_j = n_in_j a$.
\end{exm}

In what follows we adopt the following notation.\medskip

\noindent\textbf{Notation.}  Unless specified otherwise, $R$ denotes a non-trivial ring and, for an integer $n\geq 1$, $M=(M_i)_{i=1}^{n}$ is a family of $R$-modules with  bilinear maps as indicated in the definition of the $n$-trivial extension defined so that $R\ltimes_n M$ is a commutative associative ring with identity. So $R\ltimes_n M$ is indeed a commutative ring with identity. Let $S$ be a nonempty subset of  $R$ and   $N=(N_i)_{i=1}^{n}$ be a family of sets such that, for every $i $,  $N_i\subseteq M_i$. 
Then as a subset of $R\ltimes_n M$, $S\times N_1\times \cdots\times N_n$ will be denoted by $S\ltimes_n N_1\ltimes \cdots \ltimes N_n$ or simply $S\ltimes_n N$.

\section{Some ring-theoretic constructions of $n$-trivial extensions} 
In this section we investigate some ring-theoretic constructions of $n$-trivial extensions. First we investigate the graded aspect of $n$-trivial extensions.\smallskip

 For the convenience of the reader we recall the definition of graded rings. Let $\Gamma$ be a commutative additive monoid. Recall that  a ring $S$ is said to be a  \textit{$\Gamma$-graded} ring, if there is a family of subgroups of $S$,  $(S_{\alpha})_{\alpha \in \Gamma}$, such that
 $S=  \underset{\alpha \in \Gamma}{\oplus} S_{\alpha}$ as an abelian group, with $S_{\alpha}S_{\beta}\subseteq S_{\alpha+\beta}$ for all $\alpha,\beta\in \Gamma$. And an $S$-module $N$ is said to be \textit{ $\Gamma$-graded} if $N=\underset{\alpha \in \Gamma}{\oplus} N_{\alpha}$ (as an abelian group) and $S_{\alpha}N_{\beta}\subseteq N_{\alpha+\beta}$ for all $\alpha,\beta\in \Gamma$. Note   that  $S_0$ is a subring of $S$ and  each $N_{\alpha}$ is an $S_0$-module. When $\Gamma=\N_0$, a   $\Gamma$-graded ring (resp., a  $\Gamma$-graded module) will simply be called a graded ring (resp., a graded module). See, for instance, \cite{NO} and \cite{Nor}  for more details about graded rings although \cite{NO} deals with group graded rings.\medskip

Now,     $ R\ltimes_n M_1\ltimes \cdots \ltimes M_n =R\oplus  M_1 \oplus \cdots \oplus M_n$ may be considered as a graded ring for the following three different grading monoids:\medskip

\indent \textbf{\underline{As an $\N_0$-graded ring}}.  In this case we set $M_k=0$ for all $k\geq n+1$ and we extend the definition of $\varphi_{i,j}$ to all $i,j\geq 0$ as follows:   For     $i$ or $j= 0$,
\begin{center}
    $\begin{array}{rccl}
\varphi_{0,j}:& R\times M_j&\longrightarrow & M_{j}\\
               & (r,m_j) &\longmapsto & \varphi_{0,j}(r,m_j):= r m_j
\end{array}$   
and  
 $\begin{array}{rccl}
\varphi_{i,0}:& M_i\times R&\longrightarrow & M_{i}\\
               & (m_i,r) &\longmapsto & \varphi_{i,0}(m_i,r):= m_i r
\end{array}$ 
\end{center}    
are just  the multiplication of $R$  when $i=j=0$ or the $R$-actions on $M_j$ and $M_i$ respectively  when $j>0$ and $i>0$ rerspectively. For     $i,j\geq 0$ such that $i+j\geq n+1$,   we define  $\varphi_{i,j}: M_i\times M_j \longrightarrow M_{i+j}$ by $\varphi_{i,j}(m_i,m_j) =0$ for all $(m_i,m_j)\in M_i\times M_j$. Thus $ R\ltimes_n M_1\ltimes \cdots \ltimes M_n $ is an $\N_0$-graded ring $\overset{\infty}{\underset{i=0}{\oplus}}R_i$  where $R_0=R$ and $R_i=M_i$ for $i\in\N$.\medskip

\indent \textbf{\underline{As a  $\Z_{n+1}$-graded ring}}.  In this case we consider, for $a\in \Z$, the least nonnegative integer $\widehat{a}$ 
 with $\widehat{a}\equiv a\, mod(n+1)$, and we set $M_{\overline{a}}:=M_{\widehat{a}}$. Then for $a,b\in  \Z$,  we define maps $\overline{\varphi}_{\overline{a},\overline{b}}: M_{\overline{a}}\times M_{\overline{b}} \longrightarrow M_{\overline{a+b}}$ by $\overline{\varphi}_{\overline{a},\overline{b}}=\varphi_{\widehat{a},\widehat{b}}$ when $\widehat{a}+\widehat{b}\leq n$ and $\overline{\varphi}_{\overline{a},\overline{b}}$ to be the zero map when $\widehat{a}+\widehat{b}> n$. 
Then $ R\ltimes_n M_1\ltimes \cdots \ltimes M_n $ is   a $\Z_{n+1}$-graded ring $R_{\overline{0}}\oplus  R_{\overline{1}} \oplus \cdots \oplus R_{\overline{n}}$ 
  where $R_{\overline{0}}=R$ and $R_{\overline{a}}=M_a$ for $a=1,...,n$.\medskip

\indent \textbf{\underline{As a  $\Gamma_{n+1}$-graded ring}}.  Here $\Gamma_{n+1}=\{0,1,...,n\}$ is    a commutative monoid  with addition $i\widehat{+}j:=i+j$ if $i+j\leq n$ and $i\widehat{+}j:=0$ if  $i+j> n$ (so $\Z_2$ and $\Gamma_2$ are isomorphic). In this case, we define maps $\widehat{\varphi}_{i,j}$, for $i,j\in \Gamma_{n+1}$, by $\widehat{\varphi}_{i,j}=\varphi_{i,j}$ when $i=j=0$ or   $i\widehat{+}j\not =0$ and $\widehat{\varphi}_{i,j}: M_i\times M_j \longrightarrow M_{0}=R$ to be the zero map when $i\widehat{+}j  =0$. Then  $ R\ltimes_n M_1\ltimes \cdots \ltimes M_n $ is a $\Gamma_{n+1}$-graded ring $R_0\oplus  R_1 \oplus \cdots \oplus R_n$ where $R_0=R$ and $R_i=M_i$ for $1\leq i\leq n$.\medskip

 Note that each of these  gradings have the same set of homogeneous elements.\medskip

We have observed that   $ R\ltimes_n M_1\ltimes \cdots \ltimes M_n $ is an  $\N_0$-graded ring $\overset{\infty}{\underset{i=0}{\oplus}} R_i$ where $R_0=R$, $R_i=M_i$ for $i=1,...,n$ and $R_i=0$ for $i>n$. So $ R\ltimes_n M_1\ltimes \cdots \ltimes M_n $ is a graded ring isomorphic to $\overset{\infty}{\underset{i=0}{\oplus}}R_i/   \underset{i\geq n+1}{\oplus} R_i$. The following result presents the converse implication. Namely, it shows that the $n$-trivial extensions can be realised as  quotients of graded rings. 

\begin{prop}\label{prop-1} Let $\overset{\infty}{\underset{i=0}{\oplus}}S_i$ be   an $\N_0$-graded ring and $m\in  \N $.
 Then $ S_0\ltimes_m S_1\ltimes \cdots \ltimes S_m $ with the product induced by $\overset{\infty}{\underset{i=0}{\oplus}}S_i$  is naturally an $\N_0$-graded ring isomorphic to  $\overset{\infty}{\underset{i=0}{\oplus}}S_i/   \underset{i\geq m+1}{\oplus}S_i $. \end{prop}
\pr Obvious.\cqfd\smallskip

The following result presents a particular case of Proposition \ref{prop-1}.

\begin{prop}\label{prop-2}
For an $R$-module $N$, we have the following two natural ring isomorphisms:$$\begin{array}{ccl}
T_R(N)/ \underset{i\geq n+1}{\oplus}T^{i}_{R}(N) & \cong&  R\ltimes_n N\ltimes T^{2}_{R}(N)\ltimes \cdots \ltimes T^{n}_{R}(N),  \ and \\ 
S_R(N)/\underset{i\geq n+1}{\oplus}S^{i}_{R}(N) & \cong & R\ltimes_n N\ltimes S^{2}_{R}(N)\ltimes \cdots \ltimes S^{n}_{R}(N) .
\end{array}$$
\indent Moreover, suppose that $N$ is a free $R$-module with a basis $B$, then $R\ltimes_n N\ltimes S^{2}_{R}(N)\ltimes \cdots \ltimes S^{n}_{R}(N) $ is (graded) isomorphic to  $ R[\{X_b\}_{b\in B}]/(\{X_b\}_{b\in B})^{n+1}$.\\
\indent In particular,  $ R\ltimes_n R\ltimes   \cdots \ltimes R $ with the natural maps is isomorphic to $ R[X]/(X^{n+1})$.
\end{prop}
\pr Obvious.\cqfd

Our next result shows that the $n$-trivial extension of a graded ring by graded modules  has a natural grading. It is an extension  of \cite[Theorem 4.5]{AW}. 

\begin{thm} \label{thm-grad}   Let $\Gamma$ be a commutative additive monoid. Assume that   $R=\underset{\alpha \in \Gamma}{\oplus} R_i$ is     $\Gamma$-graded  and $M_{i}=\underset{\alpha \in \Gamma}{\oplus}M^{i}_{\alpha }$ is     $\Gamma$-graded as an $R$-module  for every $i\in \{1,..., n\}$, such that  $  \varphi_{i,j} (M^{i}_{\alpha} ,M^{j}_{\beta}) \subseteq M^{i+j}_{\alpha+\beta}$. Then $R\ltimes_n M_{1}\ltimes \cdots \ltimes M_{n}$   is  a   $\Gamma$-graded ring with $(R\ltimes_n M_{1}\ltimes \cdots \ltimes M_{n})_{\alpha}=R_{\alpha}\oplus M^{1}_{\alpha}\oplus \cdots \oplus M^{n}_{\alpha}$.
\end{thm}
\pr Similar  to the proof of  \cite[Theorem 4.5]{AW}.\cqfd\medskip


 In the case where $R$ is either a polynomial ring or  a Laurent polynomial ring we get  the following result in which the first assertion is an extension of \cite[Corollary 4.6 (1)]{AW}.

\begin{cor}\label{cor-poly}
The following statements are true.
\begin{enumerate}
  \item $(R\ltimes_n M_{1}\ltimes \cdots \ltimes M_{n})[\{X_\alpha\}]\cong R[\{X_\alpha\}]\ltimes_n M_{1}[\{X_\alpha\}]\ltimes\cdots\ltimes M_{n}[\{X_\alpha\}]$ for any set of indeterminates $\{X_\alpha\}$  over $R$.
  \item $(R\ltimes_n M_{1}\ltimes \cdots \ltimes M_{n})[\{X_\alpha^{\pm 1}\}]\cong R[\{X_\alpha^{\pm 1}\}]\ltimes_n M_{1}[\{X_\alpha^{\pm 1}\}]\ltimes\cdots\ltimes M_{n}[\{X_\alpha^{\pm 1}\}]$ for any set of indeterminates $\{X_\alpha\}$  over $R$.
\end{enumerate}
\end{cor}

Also, as in the classical case, we get the related (but not graded) power series case. It is a generalization of \cite[Corollary 4.6 (2)]{AW}. First recall that, for a given set of analytic indeterminates  $\{X_{\alpha} \}_{\alpha \in \Lambda}$ over $R$, we can consider three types of power series rings (see \cite{R}   for further details about generalized power series rings):
  $$ 
    R[[\{X_{\alpha} \}_{\alpha \in \Lambda}]]_{1} \subseteq  R[[\{X_\alpha \}_{\alpha \in \Lambda}]]_{2} \subseteq  R[[\{X_\alpha \}_{\alpha \in \Lambda}]]_{3}.$$

Here $$\begin{array}{lcl}
R[[\{X_\alpha \}_{\alpha \in \Lambda}]]_{1}\!&\!=\!&\!\cup \{R[[\{X_{\alpha_1},...,X_{\alpha_n} \}]] | \{\alpha_1, ... ,\alpha_{n}\}\subseteq \Lambda\},\\
R[[\{X_\alpha \}_{\alpha \in \Lambda}]]_{2}\!&\!=\!&\!\{\sum^{\infty}_{i=0} f_{i} | f_{i}\in R[\{X_\alpha \}_{\alpha \in \Lambda}] \; is \;  homogeneous \; of \; degree \; i\} \quad and\\ 
R[[\{X_\alpha \}_{\alpha \in \Lambda}]]_{3}\!&\!=\!&\!\{\sum^{\infty}_{i=0} f_{i} | f_{i} \; is \; a\;  possibly \; infinite \; sum \; of \; monomials \; of \; degree \; i\\ 
                                     &  & \ with \; at \; most \; one \; monomial \; of \; the \; form \; r_{\alpha_1, ... ,\alpha_n}X^{i_1}_{\alpha_1}\cdots X^{i_n}_{\alpha_n} for \\
                                           &  & \ each\; set\; \{\alpha_1, ... ,\alpha_n\}\; with\; i_1+ \cdots +i_n=i   \}.
\end{array}
$$ 

More generally, given a partially ordered additive monoid $(S,+,\leq)$, the \textit{generalized power series ring} $R[[X,S^{\leq}]]$ consists of all formal sums $f=\underset{s\in S}{\sum} a_sX^{s}$ where $supp(f)=\{s\in S|a_s\neq 0\}$  is Artinian and narrow (i.e., has no infinite family of incomparable elements) where addition and multiplication are carried out  in the usual way. If $\Lambda $ is a well-ordered set,  $S= \underset{\lambda \in \Lambda }{\oplus} \N_0$  and $\leq$ is the reverse lexicographic order on $S$, then $R[[X,S^{\leq}]]\cong R[[\{X_\alpha\}]]_{3}$.\smallskip

Note, that in a similar manner we can define three types of power series over a module. The routine proof of the following theorem is left to the reader.  
  
\begin{thm}\label{thm-power}
\begin{enumerate}
  \item Let $\{X_\alpha\}_{\alpha\in \Lambda}$ be a set of analytic indeterminates over $R$. Then,   for $i=1, 2, 3$, $$ (R\ltimes_n \ltimes M_1 \ltimes \cdots \ltimes M_n)[[\{X_\alpha\}_{\alpha\in \Lambda}]]_i\cong R[[\{X_\alpha\}_{\alpha\in \Lambda}]]_i \ltimes_n M_{1}[[\{X_\alpha\}_{\alpha\in \Lambda}]]_i\ltimes \cdots \ltimes M_{n}[[\{X_\alpha\}_{\alpha\in \Lambda}]]_i\, .$$
  \item Let $(S,+,\leq)$ be a partially ordered additive monoid. Then  $$(R\ltimes_n \ltimes M_1 \ltimes \cdots\ltimes M_n)[[X,S^{\leq}]]\cong R[[X,S^{\leq}]] \ltimes_n M_{1}[[X,S^{\leq}]] \ltimes \cdots \ltimes M_{n}[[X,S^{\leq}]].$$
\end{enumerate}
\end{thm}

Now, we give, as an extension of \cite[Theorem 4.1]{AW}, the following result which investigates the localization of an $n$-trivial extension. For this we need the following technical lemma.

\begin {lem}\label{lem-localization}
For every $(m_i) \in   R\ltimes_n M$ and every $k\in \{1,...,n\} $, 
$$(m_0,0,...,0,m_k,m_{k+1},...,m_n)(m_0,0,...,0,-m_k,0,...,0)=(m^{2}_0,0,...,0,e_{k+1},...,e_n)$$
where  $e_{l}=m_0m_{l}-m_k m_{l-k}$ for every $l\in \{k+1,...,n\} $. Consequently, there is an element $(f_i) $ of $R\ltimes_n M$,  such that  $$(m_i) (f_i)=(m^{2^{n}}_0,0,...,0).$$ 
\end{lem}

We will denote the   element $(f_i)$ in Lemma \ref{lem-localization} by $(\widetilde{m_i})$ so $(m_i) (\widetilde{m_i})=(m^{2^{n}}_0,0,...,0).$

\begin{thm}\label{thm-mult}  Let $S$ be a  multiplicatively closed subset of $R$ and $N=(N_i)$ be a family of $R$-modules where $N_i$ is a submodule of $M_i$ for each   $i\in\{1,...,n\}$ and $N_iN_j\subseteq N_{i+j}$ for every $1\leq i,j\leq n-1$ and $i +j\leq n$. Then the set $S\ltimes_n N$ is a multiplicatively closed subset of  $R\ltimes_n M$ and we have a ring isomorphism 
    $$(R\ltimes_n M)_{S\ltimes_n N} \cong  R_{S}\ltimes_n M_{S}$$
where $M_{S}=(M_{iS})$.
\end{thm}
\pr It is trivial to show that  $S\ltimes_n N$ is a multiplicatively closed subset of  $R\ltimes_n M$. Now in order to show the desired isomorphism, we need to make, as done in the proof of  \cite[Theorem  4.1 (1)]{AW},  the following observation:  Let $(m_i)\in R\ltimes_n M$ and  $(s_i)\in S\ltimes_n N$. Then using the notation of Lemma \ref{lem-localization}, 
$$\frac{(m_i)}{(s_i)}=\frac{(m_i)(\widetilde{s_i})}{(S_0,0,...,0)}=\frac{( m'_i )}{(S_0,0,...,0)}$$ 
where $(m'_i ) =(m_i) (\widetilde{s_i}) $  and $S_0=s_0^{2^n}$. Then the map
 $$\begin{array}{rccl}
 f:& (R\ltimes_n M)_{S\ltimes_n N}  &\longrightarrow & R_{S}\ltimes_{n} M_{S}   \\
               & \frac{(m_i) }{(s_i) } &\longmapsto & (\frac{m'_0}{S_0},\frac{m'_1}{S_0},...,\frac{m'_n}{S_0}) 
\end{array}$$ is the desired isomorphism. \cqfd

As a simple but important particular case of Theorem  \ref{thm-mult}, we get the following result which extends  \cite[Theorem 4.1 and Corollary 4.7]{AW}.
In Theorem \ref{thm-prime}, we will show that if $P$ is a prime ideal of $R$, then   $P\ltimes_n M$   is a  prime  ideal of $R\ltimes_n M$. This fact is used in the next result to show that the localization of an $n$-trivial extension at a prime ideal is isomorphic to an $n$-trivial extension. 
In what follows, we use $T(A)$ to denote the  total quotient ring of a ring $A$. In Proposition \ref{pro-zero-U-Id}, we will prove that   
$S\ltimes_n M$, where $S=R-(Z(R)\cup Z(M_1)\cup \cdots \cup Z(M_n))$, is the set of all regular elements  of $R\ltimes_n M$. Thus $T(R\ltimes_n M)=(R\ltimes_n M)_{S\ltimes_n M}$.

\begin{cor}\label{cor-mult}  The following assertions are true.
\begin{enumerate}
\item Let $P$ be a prime ideal of $R$. Then we have a ring isomorphism
  $$(R\ltimes_n M)_{P\ltimes_n M}\cong R_{P}\ltimes_n M_{P}$$
where $ M_{P}=(M_{iP}$).
  \item We have a ring isomorphism 
  $$T(R\ltimes_n M) \cong R_{S}\ltimes_n M_{S}$$
where  $S=R-(Z(R)\cup Z(M_1)\cup \cdots \cup Z(M_n))$.
  \item For  an indeterminate  $X$ over $R$, we have a ring isomorphism 
      $$(R\ltimes_n M_{1}\ltimes \cdots \ltimes M_{n})(X)\cong R(X)\ltimes_n M_{1}(X)\ltimes\cdots\ltimes M_{n}(X).$$ 
\end{enumerate}
\end{cor} 
\pr  All the proofs are similar to the corresponding ones for the classical case. \cqfd

Our next result   generalizes   \cite[Theorem 4.4]{AW}. It shows that the  $n$-trivial extension of a finite direct product of rings is a finite direct product of $n$-trivial extensions.  
For the reader's convenience we recall here some known facts on the structure of modules over a finite direct product of rings. Let $R=\displaystyle\prod_{i=1}^s R_i$ be a finite direct product of rings where $s\in \N$. For $j\in \{1,...,s\}$, we set $\bar{R_j}:=0\times \cdots \times 0 \times R_j\times 0 \times \cdots \times 0$ and, for an $R$-module $N$, $N_{j}:=\bar{R_j}N$. Then $N_j$ is a  submodule of $N$ and we have $N=N_1 \oplus \cdots \oplus N_s$. Namely,   every element  $x$ in $N$ can be written in the form   $x=x_1+ \cdots +x_s$ where  $x_j=e_jx \in N_j$ for every  $j\in \{1,...,s\}$  (here $ e_j=(0,...,0,1,0,...,0)$   with $1$ in the $j$'th place).  Note that each $N_j$ is also an $R_j$-module and    $N_{1}\times \cdots  \times N_{s}$ is an $R$-module ismorphic to $N$ via the following $R$-isomorphism: 
 \begin{center}
    $\begin{array}{ccl}
   N &\longrightarrow & N_{1}\times \cdots  \times N_{s}\\
                 \sum e_jx =x &\longmapsto &  (e_1x_1,...,e_s x_s)
\end{array} \quad \mathrm{and}\quad  
 \begin{array}{ ccl}
  N_{1}\times \cdots  \times N_{s}&\longrightarrow & N \\
                (y_{1 },...,y_{s})&\longmapsto &  \sum  y_{j}
\end{array}$
\end{center}

Now, consider   the family of commutative product maps  $ \varphi=\{\varphi_{i,j}\}_{ \underset{1\leq i,j\leq n-1}{i+j\leq n}}$ and define the following maps:
$$\begin{array}{rccl}
\varphi_{j,i,k}:& M_{j,i}\times M_{j,k}&\longrightarrow & M_{j,i+k}\\
               & (m_{j,i},m_{j,k})&\longmapsto &  \varphi_{j,i,k}(m_{j,i},m_{j,k})=e_j\varphi_{i,k}(m_{j,i},m_{j,k})
\end{array}$$
where $M_{j,i}:=\bar{R_j} M_i$ for $j\in \{1,...,s\}$ and $i\in \{1,...,n\}$.   It is easily checked that, for every $j\in \{1,...,s\}$,  $\varphi_j=\{\varphi_{j,i,k}\}_{ \underset{1\leq i,k\leq n}{i+k\leq n}}$ is a family of commutative product maps and $R_j\ltimes_{\varphi_j} M_{j,1}\ltimes  \cdots\ltimes M_{j,n}$ is a $n$-$\varphi_j$-trivial extension.
Furthermore,  $$\begin{array}{rccl}
\varphi_{i,k}:& M_{i}\times M_{k}&\longrightarrow & M_{i+k}\\
               & (m_{i},m_{k})&\longmapsto &  \varphi_{i,k}(m_{i},m_{k})=  \overset{s}{\underset{j=1}{ \sum}}  \varphi_{j,i,k}(m_{j,i},m_{j,k}).
\end{array}$$
 
With this notation  in mind, we are ready to give  the desired result.

\begin{thm}\label{thm-prod} Let $R=\displaystyle\prod_{i=1}^s R_i$ be a finite direct product of rings where $s\in \N $. Then  
$$R\ltimes_{\varphi} M_{1}\ltimes  \cdots \ltimes M_{n}\cong (R_1\ltimes_{\varphi_1} M_{1,1}\ltimes \cdots \ltimes M_{1,n})\times \cdots \times (R_s\ltimes_{\varphi_s} M_{s,1}\ltimes \cdots \ltimes M_{s,n}).$$
\end{thm}
\pr  It is easily checked that the map $ (r,m_1,...,m_n)\longmapsto((r_j,m_{j,1},...,m_{j,n}))_{1\leq j\leq s}$ is an isomorphism.\cqfd \medskip

We end this section with  two results which investigate the inverse limit and direct limit of a  system of $n$-trivial extensions. Namely, we show that, under some conditions, the inverse limit or direct limit of a    system of $n$-trivial extensions is isomorphic to an $n$-trivial extension. The inverse limit case is a generalization of \cite[Theorem 4.11]{AW}.\smallskip

Let $\Gamma$ be a directed set and  $ \{M_{\alpha } ; f_ {\alpha \beta  } \}$ be an inverse   system of   abelian groups  over $\Gamma$ (so for $\alpha\leq \beta$, $ f_ {\alpha \beta  } : M_{ \beta }\rightarrow M_{\alpha }$).  We know that the inverse limit $\underset{\longleftarrow}{lim} M_{\alpha }    $     is isomorphic to   the following subset of the direct product  $\underset{\alpha}{\prod}  M_{\alpha  }   $:
  $$M_{\infty  } :=\{(x _{\alpha  } ) _{\alpha \in \Gamma } | \lambda \leq \mu \Rightarrow     x _{\lambda } =f_{\lambda \mu }    (x _{\mu }) \}.$$
In the next result, by  $\underset{\longleftarrow}{lim} M_{\alpha }    $  we mean exactly the set $M_{\infty  }$.

\begin{thm}\label{thm-inverse-limit}
 Let $\Gamma$ be a directed set and $n\geq 1$ be an integer. Consider  a family of  inverse systems $ \{M_{i,\alpha } ; f_{i,\alpha, \beta  } \}$ over $\Gamma$    (for   $i\in \{0,...,n\}$)     which  satisfy the following conditions:
\begin{enumerate}  
    \item For every $  \alpha \in \Gamma$,   $M_{0,\alpha }=R_{\alpha }$ is a ring,
   \item For every $  \alpha \in \Gamma$ and every  $i\in \{1,...,n\}$,  $M_{i,\alpha }$  is an   $R_{\alpha }$-module, and
   \item For every   $  \alpha \in \Gamma$,  $ R_{\alpha }\ltimes_n M_{1,\alpha }\ltimes \cdots \ltimes M_{n,\alpha } $   is an  $n$-trivial extension  with a family of commutative product maps:
 $$\varphi_{i,j, \alpha}:  M_{i,\alpha }\times M_{j,\alpha }\longrightarrow   M_{i+j,\alpha }$$
which satisfy, for every $ \alpha\leq \beta$,  
$$ \varphi_{i,j, \alpha}  ( f_{i,\alpha, \beta  }(m_{i,\beta }), f_{j,\alpha, \beta  }(m_{j,\beta }) )=  f_{i+j,\alpha, \beta  }(  \varphi_{i,j, \beta}    (m_{i,\beta },m_{j,\beta })).$$
\end{enumerate}
Then  $\underset{\longleftarrow}{lim} R_{\alpha }\ltimes_n \underset{\longleftarrow}{lim} M_{1,\alpha }\ltimes \cdots \ltimes \underset{\longleftarrow}{lim} M_{n,\alpha } $   is an  $n$-trivial extension  with the following family of well-defined commutative product maps:
$$\begin{array}{cclc}
    \varphi_{i,j, \alpha}:& \underset{\longleftarrow}{lim}   M_{i,\alpha }\times \underset{\longleftarrow}{lim}  M_{j,\alpha }&\rightarrow &  \underset{\longleftarrow}{lim}  M_{i+j,\alpha }\\
                          &   ( (m_{i,\alpha })_{\alpha },  (m_{j,\alpha })_{\alpha } )&\mapsto &   (  \varphi_{i,j, \alpha}   (m_{i,\alpha }, m_{j,\alpha }))_{\alpha }.
\end{array}
 $$
Moreover, there is a natural ring ismorphism:
 $$\underset{\longleftarrow}{lim}( R_{\alpha }\ltimes_n M_{1,\alpha }\ltimes \cdots \ltimes M_{n,\alpha } )\cong   \underset{\longleftarrow}{lim} R_{\alpha }\ltimes_n \underset{\longleftarrow}{lim} M_{1,\alpha }\ltimes \cdots \ltimes \underset{\longleftarrow}{lim} M_{n,\alpha }. $$
\end{thm}
\pr The result follows using a standard argument.\cqfd


Let $\Gamma$ be a directed set and  $\{M_{\gamma };f_{\gamma\lambda}\}$ a direct system of abelian groups over $\Gamma$  (so for $\gamma \leq \lambda$,  $f_{\gamma\lambda}:M_\gamma \rightarrow M_\lambda$). We know  that the  direct limit $\varinjlim M_\gamma $ is isomorphic to
 $\underset{\gamma }{\oplus} M_\gamma  / S$ where $S$ is generated by all elements $\lambda_\beta (f_{\alpha\beta}(a_\alpha))-\lambda_\alpha(a_\alpha)$ where $\alpha \leq\beta$ and $\lambda_\lambda:M_\lambda \longrightarrow \underset{\gamma }{\oplus} M_\gamma $   is the natural inclusion map for $\lambda \in \Gamma$. Since $\Gamma$ is directed, every element of  $\underset{\gamma }{\oplus} M_\gamma/ S$ has the form $\lambda_\alpha(a_\alpha)+S$ for some $\alpha \in \Gamma$ and $a_\alpha \in M_\alpha$.

\begin{thm}\label{thm-direct-limit}
 Let $\Gamma$ be a directed set and $n\geq 1$ be an integer. Consider  a family of   direct systems $ \{M_{i,\alpha } ; f_{i,\alpha, \beta  } \}$ over $\Gamma$    (for   $i\in \{0,...,n\}$)     which  satisfy the following conditions:
\begin{enumerate}  
    \item For every $  \alpha \in \Gamma$,   $M_{0,\alpha }=R_{\alpha }$ is a ring,
   \item For every $  \alpha \in \Gamma$ and every  $i\in \{1,...,n\}$,  $M_{i,\alpha }$  is an   $R_{\alpha }$-module, and
   \item For every   $  \alpha \in \Gamma$,  $ R_{\alpha }\ltimes_n M_{1,\alpha }\ltimes \cdots \ltimes M_{n,\alpha } $   is an  $n$-trivial extension  with a family of commutative product maps:
 $$\varphi_{i,j, \alpha}:  M_{i,\alpha }\times M_{j,\alpha }\longrightarrow   M_{i+j,\alpha }$$
which satisfy, for every $ \beta\leq \alpha$,  
$$ \varphi_{i,j, \alpha}  ( f_{i,\beta,\alpha   }(m_{i,\beta }), f_{j,\beta,\alpha  }(m_{j,\beta }) )=  f_{i+j,\beta,\alpha   }(  \varphi_{i,j, \beta}    (m_{i,\beta },m_{j,\beta })).$$
\end{enumerate}
Then  $\underset{\longrightarrow}{lim} R_{\alpha }\ltimes_n \underset{\longrightarrow}{lim} M_{1,\alpha }\ltimes \cdots \ltimes \underset{\longrightarrow}{lim} M_{n,\alpha } $   is an  $n$-trivial extension  with the following family of well-defined commutative product maps:
$$\begin{array}{cclc}
    \varphi_{i,j, \alpha}:& \underset{\longrightarrow}{lim}   M_{i,\alpha }\times \underset{\longrightarrow}{lim}  M_{j,\alpha }&\rightarrow &  \underset{\longrightarrow}{lim}  M_{i+j,\alpha }\\
                          &   ( (m_{i,\alpha })_{\alpha },  (m_{j,\alpha })_{\alpha } )&\mapsto &   (  \varphi_{i,j, \alpha}   (m_{i,\alpha }, m_{j,\alpha }))_{\alpha }.
\end{array}
 $$
Moreover, there is a natural ring ismorphism:
 $$\underset{\longrightarrow}{lim}( R_{\alpha }\ltimes_n M_{1,\alpha }\ltimes \cdots \ltimes M_{n,\alpha } )\cong   \underset{\longrightarrow}{lim} R_{\alpha }\ltimes_n \underset{\longrightarrow}{lim} M_{1,\alpha }\ltimes \cdots \ltimes \underset{\longrightarrow}{lim} M_{n,\alpha }. $$
\end{thm}
\pr The result follows using a standard argument.\cqfd

\section{Some basic algebraic properties of $R \ltimes_{n} M$}

In this section we give some basic properties of $n$-trivial extensions.   Before giving the first result, we make the following  observations on situations where a subfamily of $M$ is trivial. 
\begin{observ}\label{observ}
\begin{enumerate}
\item If there is an integer $i\in \{1,..., n-1\}$ such that $M_j=0$ for every $j\in \{i+1,..., n\}$, then there is a natural ring isomorphism $$R\ltimes_n M_1\ltimes \cdots \ltimes M_{i} \ltimes 0 \ltimes \cdots \ltimes 0 \cong R\ltimes_{i} M_1\ltimes \cdots \ltimes M_{i}.$$
If $M_1=\cdots = M_{n-1}=0$, then $R\ltimes_n M$ can be represented as $R\ltimes_1 M_n$. 
However, if   $n\geq 3$ and there is an integer $i\in \{1,..., n-2\}$ such that, for $j\in \{1,..., n\}$, $M_j=0$ if and only if $j\in \{1,..., i\}$, then in general  $R\ltimes_n 0\ltimes  \cdots \ltimes 0 \ltimes M_{i+1} \ltimes  \cdots \ltimes M_n$ cannot be represented as an $n-i$-trivial extension as above. Indeed, if     for example,   $i$ satisfies $2i+2\leq n$,   then    $R\ltimes_{n-i}   M_{i+1} \ltimes  \cdots \ltimes M_n$  makes no sense (since  $\varphi_{i+1,i+1}(M_{i+1},M_{i+1})$ is a subset of $M_{2i+2}$ not of $M_{i+2}$).
  \item If $M_{2k}=0$ for every $k\in \N$ with $1\leq 2k\leq n$, then  $R\ltimes_n M$ can be   represented   as  the trivial  extension of $R$ by the $R$-module $M_1\times M_3 \times \cdots \times   M_{2n'+1}$  where $2n'+1$ is the biggest odd  integer in  $\{1,..., n\}$. Namely,   there is a natural ring isomorphism
 $$R\ltimes_n M  \cong R\ltimes_1 (M_1\times M_3 \times \cdots \times   M_{2n'+1}).$$  
   \item  If $M_{2k+1}=0$ for every $k\in \N$ with $1\leq 2k+1\leq n$, then there is a natural ring isomorphism
 $$R\ltimes_n M  \cong R\ltimes_{n''}  M_2\ltimes M_4 \ltimes \cdots \ltimes   M_{2n''}$$ where $2n''$ is 
the biggest even  integer in  $\{1,..., n\}$. In general, for every cyclic submonoid  $ G $    of $\Gamma_{n+1}$ generated by an element $g  \in \{1,..., n\}$, if $M_i=0$ if and only if $i\not \in G$, then there is a natural ring isomorphism
 $$R\ltimes_n M  \cong R\ltimes_{s}  M_g\ltimes M_{2g} \ltimes \cdots \ltimes   M_{sg}$$ where $sg$ is 
the biggest    integer in  $G\cap \{1,..., n\}$.  \medskip
\end{enumerate}
\end{observ}

As observed above, if one would discuss according to whether a subfamily of $M$ is trivial or not, then  various situations may occur. Thus, for the sake of simplicity,  we make the following convention.

 \begin{conv} \label{conv}   Unless explicitly stated otherwise, when we consider an $n$-trivial extension for a given $n$, then we implicitly suppose that $M_i \not =0$ for every $i\in \{1,..., n\}$.  This  will be used in the sequel without explicit mention.
\end{conv}

Note also that    the nature of the maps  $\varphi_{i,j}$ can affect the structure of the $n$-trivial extension.  For example, in case where $n=2$, if $\varphi_{1,1}=0$, then $R\ltimes_{2}M_{1}\ltimes M_2 \cong R\ltimes (M_1\times M_2) .$ For example, if $I\subseteq J$ is an extension of ideals of $R$, then $R\ltimes_2 I\ltimes R/J \cong R\ltimes (I\times R/J).$\\

Let us start with the following result which presents some  relations (easily established) between $n$-trivial extensions.

 \begin{prop} \label{pro-2-n-exte}
The following assertions are true.
\begin{enumerate}
    \item Let  $G $ be a submonoid of $\Gamma_{n+1}$ and consider the family of $R$-modules $M'=(M'_{i})_{i=1}^{n}$ such that $M'_i=M_i$ if $i\in G$ and $M'_i=0$ if $i\not\in G$.  Then we have the following  (natural) ring extensions:
$$R \hookrightarrow   R\ltimes_{n}   {M'}  \hookrightarrow  R\ltimes_n M.$$
In particular, for every $m\in \{1,...,n\}$,  we have the following  (natural) ring extensions:
$$R \hookrightarrow   R\ltimes_{n} 0\ltimes\cdots \ltimes 0\ltimes M_{m}\ltimes  \cdots \ltimes M_{n}   \hookrightarrow  R\ltimes_n M_{1}\ltimes  \cdots \ltimes M_{n}.$$
The extension $  R \hookrightarrow      R\ltimes_n M_{1}\ltimes  \cdots \ltimes M_{n}$ will be denoted by $i_n$.
    \item  For every   $m\in  \{1,...,n\} $,   $0 \ltimes_n 0 \ltimes \cdots\ltimes 0 \ltimes M_{m}\ltimes\cdots\ltimes M_n$ is an ideal of $R\ltimes_n M$ and an 
$R\ltimes_{j} M_{1}\ltimes  \cdots \ltimes M_{j}$-module   for every $j\in \{n-m,...,n\}$ via the action 
$$\begin{array}{ll}
(x_0,x_1,...,x_j) (0,...,0,y_m,...,y_n)&:= (x_0,x_1,...,x_j,0,...,0) (0,...,0,y_m,...,y_n)\\
&\;=(x_0,x_1,...,x_{n-m},0,...,0) (0,...,0,y_m,...,y_n).
\end{array}$$
 Moreover, the structure of  $0 \ltimes_n 0 \ltimes   \cdots \ltimes 0 \ltimes M_{m}\ltimes\cdots\ltimes M_n$ as an  ideal of  $R\ltimes_{n} M$ is the same as the $R\ltimes_{j} M_{1}\ltimes  \cdots \ltimes M_{j}$-module structure for every $j\in \{n-m,...,n\}$. 
  In particular,  the structure of the ideal $0 \ltimes_{n} 0 \ltimes    \cdots\ltimes 0\ltimes M_n$   is the same as  the one of the $R$-module $  M_n $.
\item For every $m\in  \{1,...,n\} $, we have the following natural ring isomorphisms: 
   $$  R\ltimes_n M_{1}\ltimes  \cdots \ltimes M_{n}/  0 \ltimes_n 0 \ltimes  \cdots\ltimes 0 \ltimes M_{m}\ltimes\cdots\ltimes M_n \cong R\ltimes_{m-1} M_{1}\ltimes  \cdots \ltimes M_{m-1} $$
obtained from the natural ring homomorphism:
$$\begin{array}{rccl}
\pi_{m-1}:& R\ltimes_n M_{1}\ltimes  \cdots \ltimes M_{n} &\longrightarrow & R\ltimes_{m-1} M_{1}\ltimes  \cdots \ltimes M_{m-1}  \\
               & (r,x_1,...,x_n) &\longmapsto & (r,x_1,...,x_{m-1})
\end{array}$$
where for $m=1$, $R\ltimes_{m-1} M_{1}\ltimes  \cdots \ltimes M_{m-1}=R$.
 \end{enumerate}
 \end{prop}

To give another example for the assertion $(1)$, one can show that, for $n=3$, $\{ 0,2\}$ is a submonoid  of $\Gamma_{4}$. Then   we have the following  (natural) ring extensions:
$$R \hookrightarrow   R\ltimes M_{2} \hookrightarrow R\ltimes_3 M_{1}\ltimes   M_{2} \ltimes M_{3}.$$

\begin{rem}
 We have seen  that, in the case of $n=1$, the ideal structure of $0\ltimes_1 M_1$ is the same as the $R$-module structure of $0\ltimes_1 M_1$. Actually, Nagata \cite{MN} used this   to reduce proofs of module-theoretic results to the ideal case. However, for $n\geq 2$, the $R$-module structure of $0\ltimes_n M_1\ltimes \cdots \ltimes M_{n}$ need not be the same as the ideal structure. For instance, consider the $2$-trivial extension $\Z \ltimes_2 \Z \ltimes \Z$ (with the maps induced by the multiplication in $\Z$). Then $\Z (0,1,1)=\{(0,m,m)| m\in \Z \}$ while the ideal of $\Z \ltimes_2 \Z \ltimes \Z$ generated by $(0,1,1)$ is  $0 \ltimes_2 \Z \ltimes \Z$. However, according to Proposition \ref{pro-2-n-exte} (2), $(\Z \ltimes_1 \Z )  (0,1,1) = (\Z \ltimes_2 \Z \ltimes \Z)(0,1,1)$. 
\end{rem}

The notion of extensions of   ideals under  ring homomorphisms is a natural way to construct  examples of ideals. In this context, we use the ring homomorphism $i_m$ (indicated in Proposition \ref{pro-2-n-exte}) to give such examples.

\begin{prop}\label{pro-exten-ideal}
For an  ideal $I$   of $R$, we have the following assertions:
 \begin{enumerate}
   \item The ideal  $I\ltimes_n IM_{1}\ltimes  \cdots \ltimes IM_{n}$ of $R\ltimes_n M$ is the extension of $I$   under the  ring homomorphism $i_n$, and we have the following natural
ring  isomorphism:
  $$(R\ltimes_n M)/(I\ltimes_n IM_{1}\ltimes  \cdots \ltimes IM_{n})\cong (R/I)\ltimes_n(M_1/IM_1)\ltimes \cdots \ltimes (M_n/IM_n)$$
where the multiplications are well-defined as follows:
$$\begin{array}{rccl}
\overline{\varphi}_{i,j}:& M_i/IM_i\times M_j/IM_j &\longrightarrow & M_{i+j}/IM_{i+j}  \\
               & (\overline{m_i},\overline{m_j}) &\longmapsto & \overline{m_i}\, \overline{m_j}: =\overline{\varphi}_{i,j} (\overline{m_i},\overline{m_j}):= \overline{\varphi_{i,j}(m_i,m_j)}= \overline{m_i \,m_j}\, .
\end{array}$$
  \item The ideal  $I\ltimes_n IM_{1}\ltimes  \cdots \ltimes IM_{n}$ is  finitely generated if and only if  $I$ is finitely generated.
    \end{enumerate}
\end{prop}
\pr 1.    The proof is  straightforward.\\
2.   Using $\pi_0$ it is clear that if  $ I\ltimes_n IM_1\ltimes \cdots \ltimes IM_n$ is  generated by elements $(r_j,m_{j,1},...,m_{j,n}) $  with   $j\in E$ for some set $E$,   then $I$ is generated by the $r_j$'s. 
Conversely, if  $I$ is generated by elements $r_j$ with   $j\in E$ for some set $E$, then  $ I\ltimes_n IM_1\ltimes \cdots \ltimes IM_n$ is  generated by  the $(r_j,0,...,0)$'s.\cqfd

Now, we determine the radical, prime and maximal ideals of  $R\ltimes_n M$. As in the classical case, we show that these ideals are particular cases of the homogenous ones, which   are characterized  in the next section. However, we give  these  particular cases here because of their simplicity which is reflected, using the following lemma, on the fact that they contain the  nilpotent ideal $0\ltimes_n M$  (of index $n+1$).

\begin{lem}\label{lem-rad}
Every ideal of $R\ltimes_n M$ which contains  $0\ltimes_n M$ has the form $I\ltimes_n M$ for some ideal $I$ of $R$. In this case, we have the following natural ring isomorphism:
$$R\ltimes_n M/I\ltimes_n M \cong R/I.$$
\end{lem}
\pr Let $J$ be an ideal of  $R\ltimes_n M$ which contains  $0\ltimes_n M$ and consider  the ideal  $I=\pi_{0}(J)$ of $R$ where  $\pi_{0}$ is the surjective ring homomorphism used in   Proposition \ref{pro-2-n-exte}. Then $J\subseteq I \ltimes_n M$ and by the fact that $0\ltimes_n M \subseteq J$, we deduce that $J= I \ltimes_n M$. Finally, using   $\pi_{0}$ and the fact that $\pi_{0}^{-1}(I)=J$, we get the desired isomorphism.  \cqfd

The following result is an extension of \cite[Theorem 3.2]{AW}.

\begin{thm}\label{thm-prime}
Radical ideals of $R\ltimes_n M$ have the form $I\ltimes_n M$ where $I$ is a radical ideal of $R$.\\
In particular, the maximal (resp., the prime) ideals of $R\ltimes_n M$ have the form $\mathcal{M}\ltimes_n M$ (resp, $P\ltimes_n M$) where $\mathcal{M}$ (resp., $P$) is a maximal (resp., a  prime) ideal of $R$.
\end{thm}
\pr
Using Lemma \ref{lem-rad}, it is sufficient to note that every radical ideal contains $0\ltimes_n M$ since  $(0\ltimes_n M)^{n+1}=0$.\cqfd \medskip

 Theorem \ref{thm-prime} allows us to easily determine     both the Jacobson radical and  the nilradical of $R\ltimes_n M$.

\begin{cor}\label{cor-prime}
The Jacobson radical $J(R\ltimes_n M)$ (resp., the nilradical $\nil(R\ltimes_n M)$) of $R\ltimes_n M$ is $J(R)\ltimes_n M$ (resp., $\nil(R)\ltimes_n M$) and the Krull dimension of $R\ltimes_n M$ is equal to that of $R$.
\end{cor}

We end this section with an extension of \cite[Theorems 3.5 and 3.7]{AW}  which determines, respectively, the set of zero divisors $Z(R\ltimes_n M)$, the set of units $U(R\ltimes_n M)$ and the set of idempotents $Id(R\ltimes_n M)$  of $R\ltimes_n M$. It is worth noting  that trivial extensions have been used to construct examples of rings with zero divisors that  satisfies certain properties. As mentioned in the introduction, particular $2$-trivial extensions are used  to settle some questions. Recently, in \cite{BMT},  a $2$-trivial extension is used in the context of zero-divisor graphs   to give an appropriate  example.

\begin{prop}\label{pro-zero-U-Id}
The following assertions are true.
\begin{enumerate}
    \item The set of zero divisors of $R\ltimes_n M$ is
   $$Z(R\ltimes_n M)=\{(r,m_1,...,m_n)| r\in Z(R)\cup Z(M_1)\cup \cdots \cup Z(M_n),m_i\in M_i  \; for \;  i\in\{1,..., n\}\}. $$  Hence $S\ltimes_n M$ where $S=R-(Z(R)\cup Z(M_1)\cup \cdots \cup Z(M_n))$ is the set of regular elements  of $R\ltimes_n M$.
  \item The set  of units  of $R\ltimes_n M$ is
$U(R\ltimes_n M)=U(R)\ltimes_n M$.
  \item The set  of idempotents  of $R\ltimes_n M$ is
$Id(R\ltimes_n M)=Id(R)\ltimes_n 0$.
\end{enumerate}
\end{prop}
\pr All the proofs are similar to the corresponding ones for the classical case. For completeness,   we give   a proof of the first assertion.\\
 \indent Let $(r,m_1,...,m_n)\in R\ltimes_n M$ such that $r\in Z(R)\cup Z(M_1)\cup \cdots \cup Z(M_n)$. If $r=0$, then $(0,m_1,...,m_n)(0,...,m'_n)=(0,...,0)$ for every $m'_n\in M_n$. Hence $(r,m_1,...,m_n)\in Z(R\ltimes_n M)$. Suppose $r\neq 0$. If $r\in Z(R)$, there exists a nonzero element $s\in R$ such that $rs=0$, so $(r,0,...,0)(s,0,...,0)=(0,...,0)$ and hence $(r,0,...,0)\in Z(R\ltimes_n M)$. If $r\in Z(M_i)$, for some $ i\in \{1,..., n\}$, there exists a nonzero element $m''_i$ of $ M_i$ such that $rm''_i=0$,  so $$(r,0,...,0)(0,...,0,m''_i,0,...,0)=(0,...,0).$$ Hence $(r,0,...,0)\in Z(R\ltimes_n M)$. Now, since $Z(R\ltimes_n M)$ is a union of prime ideals and $\nil(R\ltimes_n M)$ is contained in each prime ideal and using the fact that  $(0,m_1,...,m_n)\in \nil(R\ltimes_n M)$, we  conclude that $(r,m_1,...,m_n)=(r,0,...,0)+(0,m_1,...,m_n)\in Z(R\ltimes_n M)$. This gives the first inclusion.\\
\indent Conversely, let   $(r,m_1,...,m_n)\in Z(R\ltimes_n M)$. Then there is $(s,m'_1,...,m'_n) \in R\ltimes_n M-\{(0,...,0)\}$ such that  $(0,...,0)=(r,m_1,...,m_n)(s,m'_1,...,m'_n)=(rs,rm'_1+sm_1,rm'_2+m_1m'_1+sm_2,...,rm'_n+\underset{i+j=n}{ \sum} m_im'_j+sm_n)$. If
  $s\neq 0$, then  $r\in Z(R)$,  and if $s=0$, we get  $r\in Z(M_1)$ if $m'_1\neq 0$, otherwise we pass to  $m'_2$ and  so on we continue until we arrive at     $s=0 $  and $m'_i=0$  for all   $i\in \{1,..., n-1\}$. Then $rm'_n=0$ and   $m'_n  \neq 0$, so $r\in Z(M_n)$. This gives the desired inclusion.  \cqfd


\section{Homogeneous ideals of $n$-trivial extensions}  
  
The study of the classical trivial extension as a graded ring     established some interesting properties (see, for instance,  \cite[Section 3]{AW}). Namely, in \cite{AW}, studying homogeneous ideals of the trivial extension shed more light on the structure of their ideals. Then naturally one would like to extend this study to the context of $n$-trivial extensions. 
In this section we extend this study to the context of $n$-trivial extensions, where here $R\ltimes_n M$ is  a ($\N_0$-)graded ring with, as indicated in Section 3, $(R\ltimes_n M)_0=R $, $(R\ltimes_n M)_i=M_i$, for every $i\in \{1,..., n\}$,    and $(R\ltimes_n M)_i=0$  for every $i\geq n+1$. Note that we could also consider  $R\ltimes_n M$ as a $\Z_{n+1}$-graded ring or $\Gamma_{n+1}$-graded ring as mentioned in Section 3.\medskip

 For that, it is convenient to recall the following definitions: Let $\Gamma$ be a commutative additive monoid and $S=\underset{\alpha \in \Gamma}{\oplus} S_{\alpha}$ be  a    $\Gamma$-graded ring.  Let  $N=\underset{\alpha \in \Gamma}{\oplus} N_{\alpha}$  be a   $\Gamma$-graded $S$-module. For every $ \alpha\in \Gamma $, the elements of $N_{\alpha}$ are said to be \textit{homogeneous of degree $ \alpha$}.  A submodule $N'$ of $N$ is said to be \textit{homogeneous} if one of the following equivalent assertions is true.
\begin{itemize}
  \item [(1)]  $N'$ is generated by homogeneous elements,
  \item [(2)]  If $\underset{\alpha \in G'}{\sum}   n_{\alpha}\in N'$, where $G'$ is a finite subset of $\Gamma$  and each $n_{\alpha}$ is homogeneous of degree $ \alpha$, then   $n_{\alpha}\in N'$ for every $ \alpha\in G'$, or
  \item [(3)] $N'=\underset{\alpha \in \Gamma}{\oplus} (N'\cap N_{\alpha})$.
\end{itemize}

In particular, an ideal $J$ of $R\ltimes_n M$  is   homogeneous if and only if $J=(J\cap R)\oplus(J\cap M_1)\oplus \cdots \oplus(J\cap M_n)$. 
Note that    $I:=J\cap R$ is an ideal of $R$ and, for $i\in \{1,..., n\}$,   $N_i:=J\cap M_i$ is an $R$-submodule of $M_i$  which satisfies   $IM_i\subseteq N_i$ and $N_iM_j\subseteq N_{i+j}$  for evey $i$, $j\in \{1,..., n\}$. \medskip

The next result extends \cite[Theorem 3.3 (1)]{AW}. Namely, it determines the structure  of  the homogeneous ideals  of $n$-trivial extensions.\medskip

In what follows, we use the ring homomorphism  $\Pi_{0}:=\pi_{0}$ (used in Proposition \ref{pro-2-n-exte}) and, for  $i\in \{1,..., n\}$,  the following homomorphism of $R$-modules:
$$\begin{array}{rccl}
\Pi_{i}:& R\ltimes_n M_{1}\ltimes  \cdots \ltimes M_{n} &\longrightarrow & M_i  \\
               & (r,m_1,...,m_n) &\longmapsto &  m_i\, .
\end{array}$$

 \begin{thm}\label{thm1-homoge} The following assertions are true.
 \begin{enumerate}
   \item Let $I$ be an ideal of $R$  and let $C=(C_i)_{i\in \{1,..., n\}}$ be a family of $R$-modules such that $C_i\subseteq M_i$ for every $i\in \{1,..., n\}$.
Then $I\ltimes_n C$ is a (homogeneous) ideal of $R\ltimes_n M$ if and only if $IM_i\subseteq C_i$ and $C_iM_j\subseteq C_{i+j}$ for all $i,j\in \{1,..., n\}$ with $i+j\leq n$.\\
Thus  if  $I\ltimes_n C$ is an ideal of $R\ltimes_n M$, then $M_i/C_i$   is an $R/I$-module for every $i\in \{1,..., n\}$, and    we have a natural ring isomorphism $$(R\ltimes_n M_{1}\ltimes \cdots\ltimes M_{n})/(I\ltimes_n C_1\ltimes \cdots \ltimes C_n)\cong(R /I)\ltimes_n (M_1/C_1)\ltimes \cdots\ltimes (M_n /C_n)$$
where the multiplications are well-defined as follows:
$$\begin{array}{rccl}
\overline{\varphi}_{i,j}:&M_i/C_i\times M_j/C_j &\longrightarrow & M_{i+j}/C_{i+j} \\
               &  (\overline{m_i}, \overline{m_j}) &\longmapsto & \overline{m_i\, m_j}.
\end{array}$$
 In particular, $(R\ltimes_n M_{1}\ltimes \cdots\ltimes M_{n})/(0\ltimes_n C_1\ltimes\cdots\ltimes C_n)\cong R\ltimes_n (M_1/C_1)\ltimes \cdots \ltimes(M_n /C_n)$.
   \item Let $J$ be an ideal of $R\ltimes_n M$ and consider $K:=\Pi_0(J)$ and  $N_i:= \Pi_i(J)$ for every   $i\in \{1,..., n\}$.
    Then,
\begin{enumerate}
    \item $K$ is an ideal of $R$ and $N_i$ is a submodule of $M_i$ for every   $i\in \{1,..., n\}$ such that  $KM_i\subseteq N_i$  and $N_iM_j\subseteq N_{i+j}$ for every  $j\in \{1,..., n\}$  with $i+j\leq n$. Thus $K\ltimes_n N_1\ltimes \cdots \ltimes N_n$ is a homogeneous ideal of $R\ltimes_n M_{1}\ltimes \cdots\ltimes M_{n}$.
 \item $J  \subseteq K\ltimes_n N_1\ltimes\cdots\ltimes N_n  $.
 \item The ideal $J$ is homogeneous if and only if $J  = K\ltimes_n N_1\ltimes\cdots\ltimes N_n  $.
\end{enumerate}
 \end{enumerate}
 \end{thm}
 \pr
   1.  If $I\ltimes_n C_1\ltimes\cdots\ltimes C_n$ is an ideal of $R\ltimes_n M$, then  $(R\ltimes_n M_{1}\ltimes \cdots\ltimes M_{n})(I\ltimes_n C_1\ltimes\cdots\ltimes C_n)=I\ltimes_n (IM_1+ C_1)\ltimes (IM_2+C_2+C_1M_1)\ltimes \cdots \ltimes (IM_n+C_n+\underset{i+j=n}{ \sum} C_iM_j)$. Thus  $IM_i\subseteq C_i$ and $C_iM_j\subseteq C_{i+j}$ for every $i,j\in \{1,..., n\}$. \\
\indent Conversely, suppose that  we have $IM_i\subseteq C_i$ and $C_iM_j\subseteq C_{i+j}$ for all $i,j\in \{1,..., n\}$ with $i+j\leq n$. Then     $M_i/C_i$ is an $R/I$-module for every  $i\in \{1,..., n\}$ and the map
   $$\begin{array}{rccl}
f:& R\ltimes_n M_1\ltimes \cdots \ltimes M_n &\longrightarrow & (R /I)\ltimes_n (M_1/C_1)\ltimes \cdots \ltimes (M_n/C_n) \\
               & (r,m_1,...,m_n) &\longmapsto & (r+I,m_1+C_1,...,m_n+C_n)
\end{array}$$
is a well-defined surjective homomorphism  with $Kerf=I\ltimes_n C_1\ltimes \cdots \ltimes C_n$, so $I\ltimes_n C_1\ltimes \cdots \ltimes C_n$ is an ideal of $R\ltimes_n M_1\ltimes \cdots \ltimes M_n$ and $$(R\ltimes_n M_1\ltimes \cdots \ltimes M_n)/(I\ltimes_n C_1\ltimes \cdots \ltimes C_n)\cong(R /I)\ltimes_n (M_1/C_1)\ltimes \cdots \ltimes (M_n /C_n).$$
  In particular, $(R\ltimes_n M_{1}\ltimes \cdots\ltimes M_{n})/(0\ltimes_n C_1\ltimes\cdots\ltimes C_n)\cong R\ltimes_n (M_1/C_1)\ltimes \cdots \ltimes(M_n /C_n)$.\smallskip

   2.  All of the three statements are easily checked.    \cqfd


The following result  presents some properties of homogeneous ideals of $R\ltimes_n M $. It is an extension of both \cite[Theorem 3.2 (3)]{AW} and \cite[Theorem 3.3 (2) and (3)]{AW}. In particular, we determine, as an extension of \cite[Theorem 3.3  (3)]{AW},  the form of homogeneous principal ideals.  In fact,  the characterization  of homogeneous  principal ideals plays a key role in studying   homogeneous   ideals. This is due to (the easily checked) fact that an ideal $I$ of a    graded ring is homogeneous if   every principal ideal generated by an element of $I$ is homogeneous.

\begin{prop}\label{pro-homogeneous}
The following assertions are true.
\begin{enumerate}
  \item Let $I\ltimes_n N_1\ltimes \cdots \ltimes N_n$ and $I'\ltimes_n N'_1\ltimes \cdots \ltimes N'_n$ be two homogeneous ideals of $R\ltimes_n M$. Then we have the following   homogeneous ideals of $R\ltimes_n M$:
\begin{enumerate}
    \item  $(I\ltimes_n N_1\ltimes \cdots \ltimes N_n)+ (I'\ltimes_n N'_1\ltimes \cdots \ltimes N'_n)=(I+ I')\ltimes_n (N_1+ N'_1)\ltimes \cdots \ltimes(N_n+ N'_n)$,
\item $(I\ltimes_n N_1\ltimes \cdots \ltimes N_n)\cap (I'\ltimes_n N'_1\ltimes \cdots \ltimes N'_n)=(I\cap I')\ltimes_n (N_1\cap N'_1)\ltimes \cdots \ltimes(N_n\cap N'_n)$,  
\item $(I\ltimes_n N_1\ltimes \cdots \ltimes N_n)(I'\ltimes_n N'_1\ltimes \cdots \ltimes N'_n)=II'\ltimes_n (IN'_1+I'N_1)\ltimes(IN'_2+I'N_2+N_1N'_1)\ltimes \cdots \ltimes (IN'_n+I'N_n+    \underset{i+j=n}{ \sum}   N_iN'_j)$, and
\item $(I\ltimes_n N_1\ltimes \cdots \ltimes N_n) : (I'\ltimes_n N'_1\ltimes \cdots \ltimes N'_n)=((I :_{R}I')\cap (N_1:_{R} N'_1) \cap\cdots \cap(N_n:_{R} N'_n)) \ltimes_n  ((N_1:_{M_1}I')\cap (N_2 :_{ M_1} N'_1) \cap\cdots \cap (N_n :_{ M_1} N'_{n-1}))\ltimes\cdots \ltimes (N_n :_{ M_n} I') $ where  $(N_{i+j} :_{ M_i} N'_j) :=\{m_i\in M_i|  m_i N'_j \subseteq N_{i+j}\} $ for every $i,j\in \{0,..., n\}$ with $i+j\leq n$ (here $M_0 = R$, $N_0 = I$ and $N'_0 = I'$).
\end{enumerate}
  \item  A principal ideal $\langle(a,m_1,...,m_n)\rangle$ of $R\ltimes_n M$ is homogeneous if and only if  $\langle(a,m_1,...,m_n)\rangle=aR\ltimes_n(Rm_1+aM_1)\ltimes(Rm_2+aM_2+m_1M_1)\ltimes \cdots \ltimes(Rm_n+aM_n+\underset{i+j=n}{ \sum} m_iM_j)$.
  \item    For   an ideal  $J$ of $R\ltimes_n M$,   $\sqrt{J}=\sqrt{\Pi_0(J)}\ltimes_n M$. In particular, if $I\ltimes_n C_1\ltimes \cdots \ltimes C_n$ is a homogeneous ideal of  $R\ltimes_n M$, then $\sqrt{I\ltimes_n C_1\ltimes \cdots \ltimes C_n}=\sqrt{I}\ltimes_n M$.
\end{enumerate}
\end{prop}
\pr
     1. The proof for each of the first three statements is similar to the corresponding one of \cite[Theorem 25.1 (2)]{HH}.  The   last statement easily follows from the fact that the  residual of two homogeneous ideals is again homogeneous.\smallskip

     2. Apply assertion   $(1)$ and Theorem \ref{thm1-homoge} (1).\smallskip

     3.   The proof is similar to the one of \cite[Theorem 3.2 (3)]{AW}.\cqfd

It is a known fact that,   in case where $n=1$,  even if a homogeneous ideal $I\ltimes   C$ is finitely generated, the $R$-module   $C$  is not necessarily finitely generated (you can consider $\Z \ltimes \Q$ and the principal ideal  $\langle(2,0)\rangle=2\Z \ltimes \Q$ as an example). The following result presents, in this context, some particular cases obtained  using standard arguments.

\begin{prop}\label{pro-homogeneous-fg}
The following assertions are true.
\begin{enumerate}
  \item  The ideal $0\ltimes_n M$ of $R\ltimes_n M$ is finitely generated if and only if each      $R$-module $M_i$ is finitely generated.
  \item If a  homogeneous ideal $I\ltimes_n C_1\ltimes \cdots \ltimes C_n$ of $R\ltimes_n M$ is   finitely generated, then $I$ is a finitely generated ideal of $R$.\\
\indent The converse implication is true when  $C_i$ is a finitely generated $R$-module for every  $i\in \{1,..., n\}$.
\end{enumerate}
\end{prop}

From the previous section, we note that every radical (hence prime) ideal  of $R\ltimes_n M$  is homogeneous. However, it is well-known that the ideals of the classical trivial extensions are not in general homogeneous (see \cite{AW}). Then natural questions arise:\medskip

\noindent\textbf{Question 1:}\label{question1}  When every ideal of a given class $ \mathscr{I} $  of ideals of $R\ltimes_n M$  is homogeneous? \medskip

\noindent\textbf{Question 2:} For a given ring $R$ and a family of $R$-modules  $M=(M_i)_{i=1}^{n}$, what is the   class of all homogeneous  ideals of $R\ltimes_n M$? \medskip

It is clear that these questions depend on  the structure of both $R$ and   each   $M_i$. For instance, for $n=1$, if $R$ is a quasi-local  ring  with maximal $m$, then 
a proper homogeneous  ideal of  $R\ltimes R/m$  has either the form $I\ltimes  R/m$ or  $I\ltimes  0$ where $I$ is a proper  ideal of $R$. And a proper homogeneous principal ideal of  $R\ltimes R/m$  has either the form $0\ltimes  R/m$ or  $I\ltimes  0$ where $I$ is a principal ideal of $R$.   Then, for instance, a principal ideal of $R\ltimes R/m$ generated by an element  $(a,e)$ where $a$ and $e$ are both nonzero with $a\in m$, is not homogeneous.\medskip

Question 1  was   investigated    in    \cite{AW} for the case where $ \mathscr{I} $ is the class of regular ideals of  $R\ltimes_1 M$ \cite[Theorem 3.9]{AW}. Also, under the condition that $R$ is an integral domain,  a  characterization of trivial extension rings over which  every ideal is  homogeneous is given  (see \cite[Theorem 3.3 and Corollary 3.4]{AW}). Our aim in the   remainder of this section is to  extend this study to $n$-trivial extensions. It is worth noting that in the classical case (where $n=1$),  ideals $J$ with $\Pi_0(J) = 0$ are homogeneous. This shows that   the condition that all ideals  $J$ with $\Pi_0(J)\neq 0$ are homogeneous   implies that all ideals of $R\ltimes_1 M$ are homogeneous.  In the context of  $R\ltimes_n M$ for $n\geq 2$ we show that more situations can occur.\medskip

Let us begin with the  class   of ideals $J$ of $R\ltimes_n M$ with   $\Pi_0(J)  \cap  S \neq \varnothing  $   for a given subset $S$ of regular elements of $R$. \medskip

Recall that a ring $S$ is said to be \textit{pr\'{e}simplifiable} if,  for every  $a$ and $b$ in $S$: $ab=a$ implies $ a=0$ or $ b \in U(S)$. Pr\'{e}simplifiable rings were introduced and studied by Bouvier in  a series of papers  (see references)  and  they have also been investigated in \cite{AFa1,AFa2}. In \cite{AW}, the notion of a pr\'{e}simplifiable ring is used when    homogeneous ideals of the classical trivial extensions   were studied. For example,  we have that if $ R$ is pr\'{e}simplifiable but not an integral domain, then every ideal of $R\ltimes_1 M$  is homogeneous if and only if $M_1 = 0$ (see  \cite[Theorems 3.3 (4)]{AW}). This is why we first consider just  subsets  of regular elements.

\begin{thm}\label{thm-homo-regular}   Let  $S$ be a nonempty subset of  $R - Z(R)$ and let $ \mathscr{I} $ be the class of ideals $J$ of  $R\ltimes_n M$ with  $\Pi_0(J)  \cap  S \neq \varnothing  $.  Then  the following assertions are equivalent.
\begin{enumerate}
  \item Every  ideal  in $ \mathscr{I} $   is homogeneous.
\item Every principal ideal in $ \mathscr{I} $  is homogeneous.
  \item For every $s\in S$ and $i\in \{1,..., n\}$, $sM_i=M_i$.
  \item Every principal   ideal $\langle(s,m_1,...,m_n) \rangle$ with $s\in S$ has the form $I\ltimes_n M$ where $I$ is a principal ideal of $R$ with $I\cap S\neq  \varnothing $.
  \item Every ideal in $ \mathscr{I} $  has the form $I\ltimes_n M$ where $I$ is an ideal of $R$ with $I\cap S\neq  \varnothing $.
\end{enumerate}
\end{thm}
\pr $(1)  \Rightarrow (2)$. Obvious.\\
$(2)  \Rightarrow (3)$.
 Let $s\in S$ and $i\in \{1,..., n\}$. We only need to prove that $ M_i  \subseteq s M_i$. Consider an element $m_i$ of $ M_i$. Since  $s\in S$,  $\langle(s,0,...,0,m_i,0,...,0)\rangle$ is homogeneous. Then  $(s,0,...,0)\in \langle(s,0,...,0,m_i,0,...,0)\rangle$, so there is $(x,e_1,...,e_n)\in R\ltimes_n M_1\ltimes \cdots \ltimes M_n$ such that $$(s,0,...,0,m_i,0,...,0)(x,e_1,...,e_n)=(s,0,...,0).$$ Since $s$ is regular, $x=1$. Then $m_i=(-s) e_i$, as desired.\\
$(3)  \Rightarrow (4)$.  Let  $\langle(s,m_1,...,m_n) \rangle$ be a principal   ideal of $R\ltimes_n M$  with $s\in S$. By $(3)$, $$(s,m_1,...,m_n) (0\ltimes_n 0\ltimes \cdots \ltimes 0\ltimes M_n)=0\ltimes_n 0\ltimes \cdots \ltimes 0\ltimes M_n.$$
This implies that $0\ltimes_n 0\ltimes \cdots \ltimes 0\ltimes M_n \subset \langle(s,m_1,...,m_n) \rangle$. Using this inclusion and $(3)$, we get $0 \ltimes_n 0 \ltimes\cdots \ltimes 0\ltimes M_{n-1} \ltimes 0\subset \langle(s,m_1,...,m_n) \rangle$. Then inductively we get $$0 \ltimes_n 0 \ltimes\cdots\ltimes 0 \ltimes  M_{i}  \ltimes 0\ltimes \cdots\ltimes 0\subset \langle(s,m_1,...,m_n) \rangle$$
for every $i\in \{1,..., n\}$. Thus  $0\ltimes_n  M_{1}  \ltimes   \cdots\ltimes M_{n}  \subset \langle(s,m_1,...,m_n) \rangle$. Therefore by Lemma  \ref{lem-rad} and Proposition \ref{pro-homogeneous} (2),  $\langle(s,m_1,...,m_n) \rangle$ has the form $I\ltimes_n M$ where $I=sR$.\\
$(4)  \Rightarrow (5)$. Consider an ideal $J$ in $ \mathscr{I} $. Then there is an element $(s,m_1,...,m_n)\in J $ such that $s\in  \Pi_0(J)  \cap  S   $. Therefore using $(4)$ and   Lemma \ref{lem-rad}, we get the desired result.\\
 $(5)  \Rightarrow (1)$. Obvious.\cqfd

As an example, we can consider the trivial extension $S:=\Z\ltimes_2 \Z_W    \ltimes \Q$ where $\Z_W $ is  the ring of fractions of $\Z$ with respect to the multiplicatively closed subset  $W=\{2^k|k\in \N\}$ of $\Z$. Then the principal ideal $\langle(3,1,0)\rangle$ of $S$ is not homogeneous. Deny, we must have
 $(3,0,0)\in \langle(3,1,0)\rangle  $. Thus there is $(a,e,f)\in S$ such that  $(3,0,0)=  (3,1,0)(a,e,f)$. But this implies that $a=1$ and then $e=\frac{-1}{3}$, which is absurd.\bigskip

The following result is an extension of \cite[Theorem 3.9]{AW}. Recall that an ideal is said to be \textit{regular} if it contains a regular element. Here, from Proposition \ref{pro-zero-U-Id}, an ideal of $R\ltimes_n M$  is regular if and only if it contains an element $(s,m_1,...,m_n)  $ with $s\in R-(Z(R)\cup Z(M_1)\cup \cdots \cup Z(M_n))$.

\begin{cor}\label{cor-reg-id-hom}
Let $S=R-(Z(R)\cup Z(M_1)\cup \cdots \cup Z(M_n))$. Then the following assertions are equivalent.
\begin{enumerate}
   \item Every  regular ideal  of $R\ltimes_n M$  is homogeneous.
\item Every principal regular ideal  of $R\ltimes_n M$ is homogeneous.
  \item For every $s\in S$ and $i\in \{1,..., n\}$, $sM_i=M_i$  (or equivalently, $M_{iS}=M_i$).
  \item Every principal   ideal $\langle(s,m_1,...,m_n) \rangle$ with $s\in S$ has the form $I\ltimes_n M$ where $I$ is a principal ideal of $R$ with $I\cap S\neq  \varnothing $.
  \item Every   regular ideal  of $R\ltimes_n M$    has the form $I\ltimes_n M$ where $I$ is an ideal of $R$ with $I\cap S\neq  \varnothing $.
\end{enumerate}
Consequently, if $R\ltimes_n M$ is root closed (in particular, integrally closed), then every regular ideal of $R\ltimes_n M$ has the form given in $(5)$.
\end{cor}
\pr  The proof  is similar to the one of \cite[Theorem 3.9]{AW}.\cqfd


Compare the following result with \cite[Corollary 3.4]{AW}.

\begin{cor}\label{cor-domain-id-hom}
Assume that $R$ is an integral domain. Then the following assertions are equivalent.
\begin{enumerate}
  \item Every  ideal $J$ of $R\ltimes_n M$  with  $\Pi_0(J)\neq 0$  is homogeneous.
\item Every principal ideal $J$ of $R\ltimes_n M$  with  $\Pi_0(J)\neq 0$  is homogeneous.
  \item For every  $i\in \{1,..., n\}$, $ M_i$ is divisible.
  \item Every principal   ideal $\langle(s,m_1,...,m_n) \rangle$ of $R\ltimes_n M$  with $s\neq 0$ has the form $I\ltimes_n M$ where $I$ is a  nonzero principal ideal of $R$.
  \item Every ideal $J$ of $R\ltimes_n M$  with  $\Pi_0(J)\neq 0$ has the form $I\ltimes_n M$ where $I$ is a  nonzero ideal of $R$.
  \item  Every ideal of $R\ltimes_n M$  is comparable to  $0\ltimes_n M$.
\end{enumerate}
\end{cor}
\pr  The equivalence $(5)\Leftrightarrow  (6)$ is a simple consequence of Lemma \ref{lem-rad}.\cqfd
 
The proof of Theorem \ref{thm-homo-regular} shows that another situation can be considered. This is given in the following result. We use $Ann_R(H)$ to denote the annihilator of an $R$-module $H$.

\begin{thm}\label{thm-homo-regular2} Let $ \mathscr{I} $ be the class of ideals $J$  of  $R\ltimes_n M$ with  $\Pi_0(J)  \cap  S \neq \varnothing  $ where   $S$ is a nonempty subset of  $R-\{0\} $  such that,   for every  $s\in S$,  $Ann_R(s)\subseteq Ann_R(M_i)$.   Then  the following assertions are equivalent.
\begin{enumerate}
  \item Every  ideal  in $ \mathscr{I} $   is homogeneous.
\item Every principal ideal in $ \mathscr{I} $  is homogeneous.
  \item For every $s\in S$ and $i\in \{1,..., n\}$, $sM_i=M_i$.
  \item Every principal   ideal $\langle(s,m_1,...,m_n) \rangle$ with $s\in S$ has the form $I\ltimes_n M$ where $I$ is a principal ideal of $R$ with $I\cap S\neq  \varnothing $.
  \item Every ideal in $ \mathscr{I} $  has the form $I\ltimes_n M$ where $I$ is an ideal of $R$ with $I\cap S\neq  \varnothing $.
\end{enumerate}
\end{thm}
\pr We only need to prove  the implication
$(2)  \Rightarrow (3)$. Let $s\in S$ and $i\in \{1,..., n\}$ and  consider an element $m_i$ of $ M_i-\{0\}$. Since  $s\in S $,  $\langle(s,0,...,0,m_i,0...,0)\rangle$ is homogeneous. Then $(s,0,...,0)\in \langle(s,0,...,0,m_i,0...,0)\rangle$, so there is $(x,e_1,...,e_n)\in R\ltimes_n M_1\ltimes \cdots \ltimes M_n$ such that  $(s,0,...,0,m_i,0...,0)(x,e_1,...,e_n)=(s,0,...,0).$   Then $sx=s$  and, by the hypothesis on $S$,   $(x-1)m_i=0$. Therefore $m_i=xm_i=(-s) e_i$, as desired. \cqfd

For an example of a ring that satisfies the condition of the previous result,    consider a ring $R$ with   an idempotent  $e \in  R - \{1,0\}$ and set $S=\{e\}$ and $M_i=Re$ for every $i\in \{1,..., n\}$. Thus, since   $e M_i=M_i$ for every $i\in \{1,..., n\}$, every ideal $J$ of  $R\ltimes_n M$ with $e \in \Pi_0(J) $ is  homogeneous.\\

Unlike the classical case (where $n=1$), the fact that, for every  $i\in \{1,..., n\}$, $ M_i$ is divisible does not  necessarily imply that every ideal is  homogeneous. For that, we consider the $2$-trivial extension $S:=k \ltimes_2  (k\times k)  \ltimes (k\times k )$ where $k$ is a field. Then the principal ideal $\langle(0,(1,0),(0,1))\rangle$ of $S$ is not   homogeneous. Indeed, if it were homogeneous, we must have   $(0,(1,0),(0,0))\in \langle(0,(1,0),(0,1))\rangle  $. Thus there is $(a,(e,f),(e',f'))\in S$ such that  $(0,(1,0),(0,0))=  (0,(1,0),(0,1))(a,(e,f),(e',f'))$. But this implies that $(a,0)=(1,0)$ and  $(0,0)= (e,a)$, which is absurd.\medskip

This example naturally leads us to investigate when every ideal $J$ of $R\ltimes_n M$ with $\Pi_{0}(J)=0$ is homogeneous. In this context,  the notion of a pr\'{e}simplifiable module is   used. For that, recall that an $R$-module $H$ is called  $R$-\textit{pr\'{e}simplifiable} if, for every  $r\in R$ and $h\in H$, $rh=h $ implies $h=0$ or $ r \in U(R)$. For example, over an integral domain, every torsion-free module is pr\'{e}simplifiable  (see \cite{AFa1}  and also \cite{AAFS}). \medskip

In studying the question when every ideal $J$ of $R\ltimes_n M$ with $\Pi_{0}(J)=0$ is homogeneous,  several different   cases   occur.  For this we use the following lemma.

\begin{lem}\label{lem-comp2}
Let $J$ be an ideal of $R\ltimes_n M$ such that, for $i\in \{1,...,n\}$, $\Pi_0(J)=0$,...,$\Pi_{i-1}(J)=0$, $\Pi_{i}(J)\neq 0$. Then the following assertions are true.
\begin{enumerate}
    \item For $i=n$, the ideal $J$ is homogeneous and it has the form $  0\ltimes_n  0\ltimes \cdots \ltimes 0\ltimes  \Pi_{n}(J) $.
    \item  For $i\not =n$,  if  $0\ltimes_n   0\ltimes \cdots \ltimes  0  \ltimes M_{i+1}\ltimes \cdots \ltimes M_n \subset J$, then $J$ is homogeneous and it has the form $  0\ltimes_n  0\ltimes \cdots \ltimes 0\ltimes  \Pi_{i}(J) \ltimes M_{i+1} \ltimes \cdots \ltimes M_n $. 
\end{enumerate}
\end{lem}
\pr Straightforward.\cqfd

\begin{thm}\label{thm-M=mM}
Assume that $n\geq 2$ and $M_{j}$ is pr\'{e}simplifiable for a given $j\in \{1,...,n-1\} $.  Let $ \mathscr{I} $ be the class of ideals $J$  of  $R\ltimes_n M$    with $\Pi_{i}(J)=0$ for every $i\in \{0,...,j-1\} $ and $\Pi_{j}(J)\not=0$. Then the following assertions are equivalent.
\begin{enumerate}
  \item Every  ideal  in $ \mathscr{I} $   is homogeneous.
\item Every principal ideal in $ \mathscr{I} $  is homogeneous.
 \item  For every  $k\in \{j+1,...,n\} $ and every  $ m_j\in M_j-\{0\}$,  $M_k=m_jM_{k-j}$.
  \item Every principal   ideal $\langle(0,0,...,0,m_{j},...,m_n) \rangle$ with $ m_{j} \neq  0$ has the form $0\ltimes_n 0 \ltimes  \cdots\ltimes 0 \ltimes  N \ltimes M_{j+1} \ltimes \cdots  \ltimes M_n  $ where $N$ is a nonzero cyclic submodule of $M_{j}$.
  \item Every ideal in $ \mathscr{I} $  has the form $0\ltimes_n 0 \ltimes  \cdots\ltimes 0 \ltimes  N \ltimes M_{j+1}\ltimes \cdots  \ltimes M_n  $ where $N$ is a nonzero  submodule of $M_{j}$.
  \item   Every  ideal  in $ \mathscr{I} $  contains $0\ltimes_n 0 \ltimes  \cdots\ltimes 0 \ltimes    M_{j+1} \ltimes\cdots  \ltimes M_n  $.
\end{enumerate}
\end{thm}
\pr The implication $(3)\Longrightarrow(4)$ is proved similarly to the implication $(3)  \Rightarrow (4)$ of Theorem \ref{thm-homo-regular}. The implication $(6)\Longrightarrow(1)$ is a simple consequence of Lemma \ref{lem-comp2}. Then only the implication $(2)\Longrightarrow(3)$ needs a proof. Let  $k\in \{j+1,...,n\} $, $  m_{j} \in M_{j}-\{0\}$ and    $  m_{k} \in M_{k}-\{0\}$. Then the  principal   ideal $p=\langle(0,...,0,m_{j},0,...,0,m_k,0, ...,0) \rangle$ is homogeneous. This implies that $(0,...,0,m_{j},0,...,0) \in p$ and so there exists $(r,e_1,... , e_n) \in  R\ltimes_n M$  such that $$(0,...,0,m_{j},0,...,0) = (r,e_1,... , e_n) (0,...,0,m_{j},0,...,0,m_k,0, ...,0)  .$$  Then $r  m_{j} =m_{j}$ and  $r  m_{k}+ e_{k-j} m_{j}  =0$. Since $M_{j}$ is pr\'{e}simplifiable, $r$ is invertible and then $m_{k}=-r^{-1} e_{k-j} m_{j}$, as desired.    \cqfd

For examples of rings that satisfy the conditions of the previous result, we can consider  the  following two $2$-trivial extensions:  $ \Z\ltimes_2 \Z_W    \ltimes \Q$ and  $ \Z\ltimes_2 \Z_W    \ltimes \Z_W $ where $ \Z_W$  is  the ring of fractions of $\Z$ with respect to the multiplicatively closed subset  $W=\{2^k|k\in \N\}$ of $\Z$.\medskip

The following   particular cases are of interest.

\begin{cor}\label{cor-M=mM-3}
Assume that $n\geq 2$ and $M_{n-1}$ is pr\'{e}simplifiable.  Let $ \mathscr{I} $ be the class of ideals $J$  of  $R\ltimes_n M$    with $\Pi_{i}(J)=0$ for every $i\in \{0,...,n-2\} $. 
Then the following assertions are equivalent.
\begin{enumerate}
  \item Every  ideal  in $ \mathscr{I} $   is homogeneous. 
 \item  For every  $ m_{n-1}\in M_{n-1}-\{0\}$,  $M_n=m_{n-1} M_{1}$. 
  \item   Every  ideal  in $ \mathscr{I} $  is comparable to $0\ltimes_n 0 \ltimes  \cdots\ltimes 0 \ltimes  M_n  $.
\end{enumerate}
\end{cor}
\pr  $(1)\Longrightarrow(2)$. This is a particular case of the corresponding one in Theorem \ref{thm-M=mM}.\\ 
$(2)\Longrightarrow(3)$. Let $I$ be an ideal of $R\ltimes_n M$ in $ \mathscr{I} $. If $\Pi_{n-1}(I)\not =0$, then Theorem \ref{thm-M=mM} shows that $I$ contains $0\ltimes_n 0 \ltimes  \cdots\ltimes 0 \ltimes  M_n  $. Otherwise,  $\Pi_{n-1}(I)  =0$ which  means that  $0\ltimes_n 0 \ltimes  \cdots\ltimes 0 \ltimes  M_n  $   contains $I$.\\
$(3)  \Rightarrow (1)$. Let $I$ be a nonzero ideal of $R\ltimes_n M$ in $ \mathscr{I} $. If  $\Pi_{n-1}(I)\not =0$, then Theorem \ref{thm-M=mM} shows that $I$  is homogeneous. The other case is a consequence of the assertion $(1)$ of Lemma \ref{lem-comp2}.\cqfd\medskip

When $n=2$, we get the following particular case of   Corollary \ref{cor-M=mM-3}.

\begin{cor}\label{cor-M=mM-4}
Assume that  $M_{1}$ is pr\'{e}simplifiable and $n=2$.  Let $ \mathscr{I} $ be the class of ideals $J$ of  $R\ltimes_2 M$    with $\Pi_{0}(J)=0$. Then the following assertions are equivalent.
\begin{enumerate}
  \item Every  ideal  in $ \mathscr{I} $   is homogeneous. 
 \item  For every  $ m_{1}\in M_{1}-\{0\}$,  $M_2=m_1M_{1}$. 
  \item   Every  ideal  in $ \mathscr{I} $  is comparable to $0\ltimes_2 0 \ltimes    M_2 $.
\end{enumerate}
\end{cor}

When  $j=1$ in Theorem \ref{thm-M=mM}, there are additional conditions equivalent to $(1)$-$(6)$.  
 The study of this case leads us to  introduce the following notion in order to avoid   trivial situations. 

\begin{defn}\label{def-module-integral} Assume that $n\geq 2$. For $i\in \{1,..., n-1\}$ and $j\in \{2,..., n\}$ with product $ij\leq n$, $M_i$ is said to be $\varphi$-$j$-integral (where
$ \varphi=\{\varphi_{i,j}\}_{ \underset{1\leq i,j\leq n-1}{i+j\leq n}}$  is the family of multiplications) if, for any $j$ elements $  m_{i_1},..., m_{i_j}$ of $M_i$, if the product $m_{i_1}\cdots m_{i_j}=0$, then at least one of the $m_{i_k}$'s is zero. If no ambiguity arises, $M_i$ is simply called $j$-integral.
\end{defn}

\begin{cor}\label{cor-2-M=mM}
Assume that $n\geq 2$, $M_{1}$ is pr\'{e}simplifiable and  $k$-integral for every  $k\in \{2,..., n-1\}$. Let $ \mathscr{I} $ be the class of ideals $J$ of  $R\ltimes_n M$    with $\Pi_{0}(J)=0$  and $\Pi_{1}(J)\not=0$. Then the following assertions are equivalent.
\begin{enumerate}
  \item Every  ideal  in $ \mathscr{I} $   is homogeneous.
\item Every principal ideal in $ \mathscr{I} $  is homogeneous.
 \item  For every  $k\in \{2,...,n\} $ and every  $ m_1\in M_1-\{0\}$,  $M_k=m_1M_{k-1}$.
  \item For every  $k\in \{2,...,n\} $ and every nonzero elements $  m_{1_1}, ... ,  m_{1_{k-1}}\in M_1-\{0\}$, $M_k= m_{1_1}  \cdots  m_{1_{k-1}}  M_{1}$.
  \item For every  $k\in \{2,...,n\} $ and every nonzero element  $  m \in M_1-\{0\}$, $M_k=     m^{k-1} M_{1}$.
  \item Every principal   ideal $\langle(0,m_{1},...,m_n) \rangle$ with $ m_{1} \neq  0$ has the form $0\ltimes_n    N   \ltimes     M_{2}    \cdots  \ltimes M_n  $ where $N$ is a nonzero cyclic submodule of $M_{1}$.
  \item Every ideal in $ \mathscr{I} $  has the form $0\ltimes_n    N \ltimes        M_{2}   \cdots  \ltimes M_n  $ where $N$ is a nonzero submodule of $M_{1}$.
  \item   Every  ideal  in $ \mathscr{I} $  contains $0\ltimes_n 0 \ltimes     M_{2} \ltimes\cdots  \ltimes M_n  $.
\end{enumerate}
\end{cor}
\pr   The equivalences  $(3)\Leftrightarrow(4 ) \Leftrightarrow(5)$ are easily proved.\cqfd\smallskip

The following result shows that, in fact, the conditions of Corollary \ref{cor-2-M=mM} above are necessary  and sufficient to show that every ideal $J$ of $R\ltimes_n M$ with $\Pi_{0}(J)=0$ is homogeneous. Note that Corollary \ref{cor-M=mM-4} presents the case $n=2$. Thus in the following result we may assume that $n\geq 3$.

\begin{cor}\label{cor-M=mM-all-0}
Assume that $n\geq 3$ and $M_{1}$ is pr\'{e}simplifiable and  $k$-integral for every  $k\in \{2,..., n-1\}$. Then the following assertions are equivalent.
\begin{enumerate}
  \item Every ideal $J$ of $R\ltimes_n M$ with $\Pi_{0}(J)=0$ and $\Pi_{1}(J)\not=0$ is homogeneous.
 \item  For every $j\in \{1,...,n-1\} $, $M_{j}$ is pr\'{e}simplifiable and every  ideal $J$  of  $R\ltimes_n M$    with $\Pi_{0}(J)=0$  is homogeneous.
\end{enumerate}
\end{cor}
\pr  We only need to prove that   $(1)\Rightarrow (2)$. Let  $j\in \{1,...,n-1\} $ and consider $m_j\in M_j-\{0\}$. Let $r\in R$ such that $rm_j=m_j$. By 
  Corollary \ref{cor-2-M=mM} (4), there are  $m_{1_1}, ..., m_{1_j}\in M_1-\{0\}$ such that $m_j=m_{1_1} \cdots m_{1_j}$. Then $rm_{1_1}  \cdots m_{1_j}=m_{1_1}  \cdots m_{1_j}$ which implies that $(rm_{1_1}- m_{1_1})  m_{1_2}  \cdots m_{1_j}=0$. Now, since $M_{1}$ is $k$-integral for every  $k\in \{2,..., n-1\}$, $rm_{1_1}- m_{1_1}=0$. Therefore  $r$ is invertible since $M_{1}$ is pr\'{e}simplifiable. So $M_{j}$ is pr\'{e}simplifiable.\\
\indent Now, to prove that every  ideal   $J$  of  $R\ltimes_n M$    with $\Pi_{0}(J)=0$   is homogeneous, it suffices to prove  that   $M_k = m_jM_{k-j}$ for every $k \in \{2, ..., n\}$, every $j \in \{1, ..., k-1\}$ and every $m_j\in M_j-\{0\}$ (by Theorem \ref{thm-M=mM}).  The case where $k=2$ is trivial. Thus fix $k \in \{3, ..., n\}$ and $j \in \{1, ..., k-1\}$.  Consider $m_j\in M_j-\{0\}$  and $m_k\in M_k-\{0\}$. We prove that $m_k=m_jm_{k-j}$ for some  $m_{k-j}\in M_{k-j}-\{0\}$. By Corollary \ref{cor-2-M=mM} $(5)$, $m_j=     m^{j-1} m_{1}$ for some  $m,m_{1}\in M_{1}-\{0\}$. And, by Corollary \ref{cor-2-M=mM} $(3)$,  $m_k=m_1m_{k-1}$ for some  $m_{k-1}\in M_{k-1}-\{0\}$. Also, by  Corollary \ref{cor-2-M=mM} $(5)$, $m_{k-1}= m^{k-2} m'_{1}$ for some  $m'_{1}\in M_{1}-\{0\}$. Then $m_k=m^{k-2} m_1 m'_{1}=(m^{j-1} m_1 ) (m^{k-j-1}m'_{1})=m_jm_{k-j}$ where $m_{k-j}=m^{k-j-1}m'_{1}\in  M_{k-j}-\{0\}$, as desired.\cqfd\smallskip

Finally, we give a case when we can characterize rings in which every ideal is homogeneous. 
Note that, when $R$ is  a ring with    $aM_i=M_i$ for every  $i\in \{1,..., n-1\}$ and every $a\in R-\{0\}$, and $M_i=     m^{i-1} M_{1}$  for every $i\in \{2,...,n\} $ and every nonzero element  $  m \in M_1-\{0\}$, then $R$ is an integral domain and $M_i$ must be torsion-free for every  $i\in \{1,..., n-1\}$. 

\begin{cor}\label{cor-global-homog}
   Suppose that  $n\geq 2$ and $R$ is an integral domain. Assume that   $M_i $ is   torsion-free, for every $i\in \{1,..., n-1\}$, and  that $M_{1}$  is  $k$-integral for every  $k\in \{2,..., n-1\}$. Then the following assertions are equivalent.
 \begin{enumerate}
    \item Every ideal of $R\ltimes_n M$ is homogeneous.
    \item The following two conditions are satisfied:
\begin{itemize}
    \item[i.] For every  $i\in \{1,..., n\}$, $M_i$ is divisible, and
    \item[ii.] For every $i\in \{2,...,n\} $ and every   $ m_1\in M_1-\{0\}$, $M_{i}=m_1M_{i-1}$.
\end{itemize}
  \end{enumerate}
 \end{cor}
\pr
 Simply use Corollaries \ref{cor-domain-id-hom} and \ref{cor-2-M=mM} and  Theorem  \ref{thm-M=mM}.\cqfd

It is easy to show that the two $n$-trivial extensions $\Z\ltimes_n \Q\ltimes \cdots\ltimes \Q$ and $\Z\ltimes_n \Q\ltimes \cdots\ltimes \Q\ltimes \Q/\Z$ satisfy the conditions of the   above result and so every ideal of  these rings  is   homogeneous.\medskip

We end this section with the following particular case.

\begin{cor}\label{cor-M=mM-field}  Suppose that  $n\geq 2$. Consider the $n$-trivial extension $S:=k \ltimes_n   E_1 \ltimes  \cdots\ltimes E_n $ where $k$ is a field and, for $i\in \{1,...,n\} $, $E_i$  is a $k$-vector space. Suppose that   $E_{1}$  is  $k$-integral for every  $k\in \{2,..., n-1\}$. Then the following assertions are equivalent. 
\begin{enumerate}
  \item Every  ideal of $S $   is homogeneous.
 \item  For every  $k\in \{2,...,n\} $,     every   $j\in \{1,...,k-1\} $ and every  $  e_j\in E_j - \{0\}$,  $E_k=e_jE_{k-j}$.  
\end{enumerate} 
\end{cor}

As a particular case, we can consider a field extension  $K\subseteq F$, then  every ideal of $S:=K \ltimes_n   F \ltimes  \cdots\ltimes F $  is homogeneous. Namely, every proper ideal of $S$ has the form  $0 \ltimes_n  0 \ltimes\cdots \ltimes 0\ltimes  N \ltimes  F \ltimes  \cdots\ltimes F $ where $N$ is a  $K$-subspace of $F$.\medskip



\section{Some ring-theoritic  properties of $R \ltimes_{n} M$}
In this section, we   determine when $R \ltimes_{n} M$   has certain ring properties such
as being Noetherian, Artinian,   Manis valuation, Pr\"{u}fer, chained, arithmetical,  a $\pi$-ring, a  generalized ZPI-ring or a PIR. We end the section with a remark on a question posed in \cite{A86} concerning $m$-Boolean rings.\medskip

We begin by   characterizing  when the $n$-trivial extensions are Noetherian (resp., Artinian).
 The following result extends \cite[Theorem 4.8]{AW}.

\begin{thm}\label{thm-noeth}
 The ring $R\ltimes_n M$ is Noetherian (resp., Artinian) if and only if $R$ is Noetherian (resp., Artinian) and, for every $i\in \{1,...,n\}$, $M_i$   is finitely generated.
\end{thm}
\pr  Similar  to the proof of  \cite[Theorem 4.8]{AW}.\cqfd\medskip


 The following result is an extension of \cite[Theorem 4.2 and Corollary 4.3]{AW}. It  investigates the integral closure of $R\ltimes_n M$ in  the total quotient ring $T(R\ltimes_n M)$ of  $R\ltimes_n M$.

\begin{thm}\label{thm-closure}
 Let $S=R-(Z(R)\cup Z(M_1)\cup \cdots \cup Z(M_n))$. If $R'$ is the integral closure of $R$ in $T(R)$, then $(R'\cap R_{S})\ltimes_n M_{1S}\ltimes \cdots \ltimes M_{nS}$ is the integral closure of $R\ltimes_n M $ in $T(R\ltimes_n M)$.\\
In particular,
\begin{enumerate}
  \item  If $R$ is an integrally closed ring, then $R\ltimes_n M_{1S}\ltimes \cdots \ltimes M_{nS}$ is the integral closure of $R\ltimes_n M_1\ltimes \cdots \ltimes M_n$ in $T(R\ltimes_n M_1\ltimes \cdots \ltimes M_n)$, and
  \item If $Z(M_i)\subseteq Z(R)$ for all $i\in \{1,..., n\}$, then $R\ltimes_n M_{1S}\ltimes \cdots \ltimes M_{nS}$ is integrally closed if and only if $R$ is integrally closed.
\end{enumerate}
\end{thm}
\pr All statements can be proved similarly to the corresponding ones  of \cite[Theorem 4.2 and Corollary 4.3]{AW}.
\cqfd

It is worth noting as in the classical case  that $R\ltimes_n M$ can be integrally closed without $R$ being  integrally closed (see the example given after \cite[Corollary 4.3]{AW}).\medskip

Similar to the classical case \cite[Theorem 4.16 (1) and (2)]{AW}, as a consequence of Theorem \ref{thm-closure} and Corollary \ref{cor-reg-id-hom}, we give the following result which  characterizes when $R\ltimes_n M $ is (Manis) valuation and when it is Pr\"{u}fer. First, recall these two notions.\smallskip

 Let $S$ be a subring of a ring $T$, and let $P$ be a prime ideal of $S$. Then $(S,P)$ is called a \textit{valuation pair} on $T$ (or just $S$ is a \textit{valuation}
ring on $T$) if there is a surjective valuation    $v:T\longrightarrow G\cup \{\infty\}$ ($v(xy)=v(x)+v(y)$, $v(x+y)\geq min\{v(x),v(y)\}$, $v(1)=0$ and $v(0)=\infty$) where $G$ is a totally ordered abelian group, with $S=\{x\in T|v(x)\geq 0\}$ and $P=\{x\in T|v(x)>0\}$. This is equivalent to if $x\in T-S$, then there exists $x'\in P$ with $xx'\in S-P$. A valuation ring $S$ is called a \textit{(Manis) valuation} ring if $T=T(S)$. Also, $S$ is called a \textit{Pr\"{u}fer} ring if every finitely generated regular ideal of $S$ is invertible. This is equivalent to every overring of $S$ being integrally closed (see \cite{HH} for more details).

\begin{cor}\label{cor-Manis-valua-prufer}
Let $S=R-(Z(R)\cup Z(M_1)\cup \cdots \cup Z(M_n))$.
\begin{enumerate}
  \item  $R\ltimes_n M$ is a Manis valuation ring if and only if $R$ is a valuation ring on $R_S$ and $M_i= M_{iS}$ for every $i\in \{1,..., n\}$.
  \item  $R\ltimes_n M$ is a Pr\"{u}fer ring if and only if, for every finitely generated ideal $I$ of $R$ with $I\cap S\neq  \varnothing $, $I$ is invertible  and $M_i= M_{iS}$ for every $i\in \{1,..., n\}$.
\end{enumerate}
      \end{cor} 

Now,  as an extension of  \cite[Theorem 4.16 (3)]{AW},   we  characterize when $R\ltimes_n M $ is  a chained ring.  Recall that a  ring $S$ is said to be \textit{chained}  if the set of ideals of $S$ is totally ordered by inclusion.  \\ 
\indent  As an exception to Convention \ref{conv},  in the following  results (Lemma \ref{lemma-chained1}, Theorem  \ref{thm-chained}, Corollary \ref{cor-arith}, Lemma \ref{lem-ZPI} and Theorem \ref{thm-ZPI}), a  module in the family associated to an $n$-trivial extension  can be zero.\\
\indent The proof of the desired result uses the following lemma which gives another   characterization of a particular  $n$-trivial extension with the property that every ideal is homogeneous.

\begin{lem}\label{lemma-chained1} Assume that $R$ is quasi-local with maximal ideal $m$. Suppose also that at least one of the modules of the family $M$ is nonzero. Then 
 every ideal of $R\ltimes_n M $ is homogeneous if and only if  the    following three conditions are satisfied:   
\begin{enumerate}
    \item $R$ is an integral domain.
    \item For every  $i\in \{1,..., n\}$, $M_i$ is divisible.
    \item  For every $1\leq i\leq j\leq n$ (when  $n\geq 2$), if $M_i\not = 0$ and $M_j \not = 0$, then $M_{j-i}\not = 0$ and $eM_i=M_j$ for every $e\in M_{j-i}$.
\end{enumerate} 
In this case, each ideal has one of the forms $I\ltimes_n M$, for some ideal $I$ of $R$, or $0\ltimes_n 0 \ltimes  \cdots\ltimes 0 \ltimes  N \ltimes M_{j+1}\ltimes \cdots  \ltimes M_n  $ where $N$ is a nonzero  submodule of $M_{j}$ for some $j\in \{1,..., n\}$.
\end{lem}
\pr  $\Longrightarrow$ Clearly the first assertion is a simple consequence of the second one. Then we only need to prove the second and the third assertions.\smallskip

$(2)$.   Let $r\in R - \{0\}$ and $i\in \{1,..., n\}$. Consider  an element $m_i\in M_i$. If $r\not\in m$, the maximal ideal of $R$, then $r$ is invertible and trivially we get the result.  Next assume  $r\in m$. By hypothesis, the ideal  $\langle(r,0,...,0,m_i,0,...,0)\rangle$ is  homogeneous, so there is $(r',m'_1,...,m'_n)$ such that $$(r,0,...,0)=(r,0,...,0,m_i,0,...,0)(r',m'_1,...,m'_n).$$
Then $rr'=r$ and $0=rm'_i+r'm_i$. Thus  $r'$ cannot be in $m$, so    $r'$ is invertible and thus  $ m_i=-(r')^{-1}r m'_i$, as desired.\smallskip

\indent $(3)$.  Let  $1\leq i\leq j\leq n$ such that $M_i\not = 0$ and $M_j \not = 0$.   Consider $m_i  \in  M_i - \{0\}$ and $m_{j} \in M_{j} - \{0\}$. By hypothesis, $\langle(0,...,0,m_i,0,...,0,m_{j},0,...,0)\rangle$ is homogeneous. Then  $$(0,...,0,m_{j},0,...,0)=(0,...,0,m_i,0,...,0,m_j,0,...,0)(r',m'_1,...,m'_n)$$
for some $(r',m'_1,...,m'_n) \in  R\ltimes_n M $.  This implies that $$r'm_i=0\quad \mathrm{and}\quad r'm_{j}+m_im'_{j-i} =m_j.$$
If $M_{j-i}=0$, we get $r'm_i=0$ and $(r'-1)m_{j}=0$. This is impossible since either $r'$ or $r'-1$ is invertible. Then  $M_{j-i}\not =0$. Now, suppose that $r'\neq 0$. By $(2)$,  there exists $m''_{j-i}\in M_{j-i}$ such that $m'_{j-i}=r'm''_{j-i}$. Hence using the fact that $r'm_i=0$, the equality $ r'm_{j}+m_im'_{j-i} =m_j$  becomes $r'm_{j}=m_{j}$.  As in the previous case, this is impossible. Therefore $r'=0$ and this gives the desired result.\smallskip

$\Longleftarrow$ We only need to prove that every principal ideal $\langle(s,m_1,...,m_n)\rangle$ of $R\ltimes_n M$ is homogeneous. For this, distinguish two cases $s \neq 0$ and  $s=0$ and follow an argument similar to that of $(3)\Rightarrow (4)$ of Theorem \ref {thm-homo-regular}.\cqfd

\begin{thm}\label{thm-chained} Assume that $n\geq 2$ and that at least one of the modules of the family $M$ is nonzero. Then the ring $R\ltimes_n M$ is chained if and only if the following conditions are satisfied:
\begin{enumerate}
    \item $R$ is a valuation domain,
    \item  For every  $i\in \{1,..., n\}$, $M_i$ is  divisible,
    \item  For every $1\leq i\leq j\leq n$, if $M_i\not = 0$ and $M_j \not = 0$, then $M_{j-i}\not = 0$ and $eM_i=M_j$ for every $e\in M_{j-i}$, and 
  \item  For every  $i\in \{1,..., n\}$, the set of   all (cyclic) submodules of $M_i$ is totally ordered by inclusion.  
\end{enumerate}
\end{thm}
\pr $\Longrightarrow $ First, we prove that $R$ is  a chained ring. Consider two ideals $I$ and $J$ of $R$. Then  $I\ltimes_n M$ and  $J\ltimes_n M$ are two ideals of $R\ltimes_n M$. Then they are comparable and so are $I$ and $J$ as desired. A similar argument can be used to prove the last assertion.\\
\indent  Now, we prove that  every ideal of $R\ltimes_n M $ is homogeneous. Then by Lemma \ref{lemma-chained1}, we get the other assertions. Consider a nonzero ideal $K$ of $R\ltimes_n M $. If $\Pi_0(K)\not=0$, then necessarily $0\ltimes_n M \subset K$. Then by Lemma  \ref{lem-rad}, $K$ is homogeneous. Now, let  $i\geq 1$ be the smallest integer such that  $\Pi_i(K)\not=0$. If $i=n$, then by the first assertion of Lemma  \ref{lem-comp2}, $K$ is homogeneous. If $i\neq n$, then necessarily  $0\ltimes_n   \cdots \ltimes  0  \ltimes M_{i+1}\ltimes \cdots \ltimes M_n \subset K$. Thus by the second assertion of Lemma  \ref{lem-comp2}, $K$ is homogeneous, as desired.\smallskip

$\Longleftarrow$ Using Lemma  \ref{lem-comp2}, we deduce that   any two ideals $I$ and $J$ of $R\ltimes_n M$ have the forms  $I=0\ltimes_n   \cdots \ltimes  0  \ltimes I_i\ltimes M_{i+1}\ltimes \cdots \ltimes M_n$  and $J=0\ltimes_n   \cdots \ltimes  0  \ltimes J_j\ltimes M_{j+1}\ltimes \cdots \ltimes M_n$ for some $i,j\in \{0,..., n\}$ where $I_i$ and $J_j$ are submodules of  $M_i$ and $M_j$ respectively (here $M_0=R$). If $i\not = j$, then  obviously  $I$  and $J$   are comparable. If $i=j$, then  using   the first   and the last assertion, we can show that $I_i$ and $J_j$ are comparable and so are $I$ and $J$, as desired. \cqfd

Using Theorem \ref{thm-chained} and   Corollary \ref{cor-mult}, we get an extension  of \cite[Theorem 4.16 (4)]{AW}  which   characterizes when $R\ltimes_n M$ is arithmetical. Recall that a ring $S$ is \textit{arithmetical} if and only if $S_{P}$ is chained for each prime (maximal) ideal $P$ of $S$. Also, recall that, for a ring $S$, an $S$-module $H$ is called \textit{arithmetical} if, for each prime (maximal) ideal $P$ of $S$, the set of submodules of $H_P$ is totally ordered by inclusion. Finally, recall  the support of an  $S$-module $H$, $supp(H)$, over a ring $S$ is the set of all prime ideals $P$ of $S$  such that $H_P\not=0$.

\begin{cor}\label{cor-arith}
The ring $R\ltimes_n M$ is arithmetical if and only if the following   conditions are satisfied:
 \begin{enumerate}
  \item $R$ is  arithmetical,
  \item  For every  $i\in \{1,..., n\}$, $M_i$ is an arithmetical $R$-module,
    \item  For every $P\in \underset{i}{\cup} supp(M_i)$, $R_P$ is a valuation domain,
    \item   For every  $i\in \{1,..., n\}$ and every $P\in   supp(M_i)$, $M_{iP}$ is a divisible $R_P$-module, and
 \item  For every $1\leq i\leq j\leq n$, if  $P\in   supp(M_i)\cap supp(M_j)$, then  $P\in   supp(M_{j-i}) $ and $eM_{iP}=M_{jP}$ for every $e\in M_{(j-i)P}$.
\end{enumerate}
\end{cor}
 

Recall that a ring $S$ is called  a \textit{generalized ZPI-ring} (resp., a \textit{$\pi$-ring}) if every proper ideal (resp.,
proper principal ideal) of $S$ is a product of prime ideals. An integral domain which is a
$\pi$-ring is called a $\pi$-domain. Clearly, a generalized ZPI-domain is nothing but a Dedekind domain. It is well known (for example, see [28, Sections 39 and
46]) that $S$ is a $\pi$-ring (resp., a generalized ZPI-ring, a principal ideal ring (PIR)) if and only if
$S$ is a finite direct product of the following types of rings: (1) $\pi$-domains (resp., Dedekind
domains, PIDs) which are not fields, (2) special principal ideal rings (SPIRs) and (3) fields.\medskip

Our next results extend  \cite[Lemma 4.9 and Theorem 4.10]{AW}. They characterize when $R \ltimes_{n} M$ is a $\pi$-ring, a generalized ZPI-ring  or a PIR.

\begin{lem}\label{lem-ZPI}
If $R \ltimes_{n} M$ is a $\pi$-ring (resp., a generalized ZPI-ring, a PIR), then $R$ is a $\pi$-ring (resp., a generalized ZPI-ring,
a PIR).  Hence $R = R_{1} \times\cdots \times R_{s}$ where $R_{i}$ is either (1) a $\pi$-domain (resp., a Dedekind domain, a 
PID) but not a field, (2) an SPIR, or (3) a field. Let $ M_{j,i} = (0 \times\cdots 0 \times R_{j} \times 0\times  \cdots  \times 0)M_{i} $ where $1\leq i \leq n$ and $1\leq j \leq s$. If $R_{i}$ is a domain or SPIR, but not a field, then $ M_{j,i}=0$ while if $R_{i}$ is a field, $ M_{j,i}=0$ or $ M_{j,i} \cong R_{i}$.\\
\indent Conversely, if $R = R_{1} \times\cdots \times R_{s}$ and $ M_{i} = M_{1,i} \times\cdots \times M_{s,i}$ are as above and $R$ is a $\pi$-ring
(resp., a generalized ZPI-ring, a PIR), then $R \ltimes_{n} M$ is a $\pi$-ring (resp., a generalized ZPI-ring, a PIR).
\end{lem}
\pr Using Theorem \ref{thm-prod}, the proof is similar to that of  \cite[Lemma 4.9]{AW}. \cqfd

\begin{thm}\label{thm-ZPI}
 $ R \ltimes_{n} M $ is a $\pi$-ring (resp., a generalized ZPI-ring, a PIR) if and only if $R$ is a $\pi$-ring (resp., a generalized ZPI-ring, a 
PIR) and $M_{i}$ is cyclic with annihilator $P_{i_{1}}\cdots P_{i_{s}}$ where $P_{i_{1}}$, ...,$ P_{i_{s}}$ are some idempotent
maximal ideals of $R$ (if $i_s = 0$, $Ann(M_{i}) = R$, that is, $M_{i} = 0$).
\end{thm}
\pr  Similar to the   proof of  \cite[Theorem 4.10]{AW}. \cqfd
 
 We end the section with a remark on a question posed in \cite{A86}. Recall that a ring $R$ is called \textit{$m$-Boolean} for some $m\in \N$, if $char\, R=2$ and $x_1x_2\cdots x_m(1+x_1)\cdots (1+x_m)=0$ for all $x_1, ...,x_m\in R$. Thus Boolean rings are just $1$-Boolean rings. It is shown in \cite[Theorem 10]{A86} that   $2$-Boolean rings can be represented as trivial extensions.  Namely, it is proved that if $R$ is $2$-Boolean, then $R\cong B \ltimes \nil(R)$ where $B=\{b\in R| b^{2}=b\}$ (\cite[Theorem 10]{A86}). Based on this result the following  natural  question is posed (see \cite[page 74]{A86}):  Whether \cite[Theorem 10]{A86} can be extended to $m$-Boolean rings for $m\geq 2$? \\
\indent One can ask whether  the $n$-trivial extension is the suitable construction  to solve this question. Using  \cite[Theorem 6]{A86}, one can show that the amalgamated algebras along an ideal (introduced in \cite{DFF}) resolve partially this question.  Recall that, given a ring homomorphism  $f:A\longrightarrow B$ and an ideal $J$ of $B$, \textit{the amalgamation of $A$ with $B$ along $J$ with respect to $f$} is the following subring of $A\times B$: 
  $$A\bowtie^{f} B=\{(a,f(a)+j)| a\in A, j\in J\}.$$   Note that $A\bowtie^{f} B\cong A\dot{\oplus}J$ where $A\dot{\oplus}J\subseteq A\times B$ is the ring whose underlying group is $A\oplus J$ with multiplication given by $(a,x)(a',x')=(aa',ax'+a'x+xx')$ for all $a,a'\in A$ and $x,x'\in J$. Here  $J$ is an $A$-module via $f$ and then $ax':=f(a)x'$ and  $a'x:=f(a')x$  (see \cite{DFF} for more details). Now, if $R$ is $m$-Boolean for $m\geq 2$, then from \cite[Theorems 6 and 7]{A86}, $R=B\oplus \nil(R)$ where $B=\{b\in R | b^{2}=b\}$. Then $$R\cong B\dot{\oplus} \nil(R)\cong B\bowtie^{i} \nil(R)$$  where $i:B\hookrightarrow R$ is the canonical injection.\\
Actually any $n$-trivial extension $R\ltimes_n M$ can be seen as the  amalgamation of $R$ with $ R\ltimes_n M$ along $0\ltimes_n M$ with respect to the canonical injection. This leads to pose the following question for every  $m\geq 2$: Is any $m$-Boolean ring an $m$-trivial extension?

\section{Divisibility properties of $R \ltimes_{n} M$}

Factorization  in commutative rings with zero divisors was  first investigated in a series of papers by Bouvier, Fletcher  and  Billis   (see References), where the focus had been on the unicity property. The papers  \cite{AFa3,AFa1,AFa2} marked  the start of a systematic study of factorization in commutative rings with zero divisors.  Since then, this theory has attracted the interest of a number of authors. The study of divisibility properties of the classical trivial extension has lead to some interesting examples and then to answering several questions  (see \cite[Section 5]{AW}). In this section we are interested in extending a part of this study to the context of $n$-trivial extensions.\medskip

First, we recall the following definitions. Let $S$  be a commutative ring and $H$ an $S$-module. Two elements $e, f \in  H$ are said to be \textit{associates} (written $e\sim f$) (resp., \textit{strong associates} (written $e\approx f$), \textit{very strong associates} (written $e\cong f$)) if $Se = Sf$ (resp., $e=uf$ for some $u\in U(S)$, $e\sim f$ and either $e=f=0$ or $e=rf$ implies $r\in U(S)$). Taking $H = S$ gives the notions of ``associates" in $S$. We say that $H$ is \textit{strongly associate} if for every $e, f \in H$, $e\sim f\Rightarrow e\approx f$. When  $S$ is strongly associate  as an $S$-module, we also say that $S$ is \textit{strongly associate}. Finally, recall that $H$ is said to be \textit{$S$-pr\'{e}simplifiable} if for $r\in S$ and $e\in H$, $re=e\Rightarrow r\in U(S)$ or $e=0$. If $S$  is  $S$-pr\'{e}simplifiable we only say that $S$  is \textit{pr\'{e}simplifiable}.\medskip 


We begin with an extension of  \cite[Theorem 5.1]{AW}.  

\begin{prop}\label{propo-presim}
Let $R\subseteq S$ be a ring extension such that $U(S)\cap R=U(R)$.
\begin{enumerate}
  \item If $S$ is pr\'{e}simplifiable, then every $R$-submodule of $S$ is pr\'{e}simplifiable. In particular, $R$ is pr\'{e}simplifiable.
  \item  Suppose that $S=R\oplus N$ as an $R$-module where $N$ is a nilpotent ideal of $S$ which satisfies either $N^{2}=0$ or $N=\underset{i \in \N}{\oplus}   N_i   $ as an $R$-module where  $S=R\oplus    N_1  \oplus N_2  \oplus \cdots$ is a graded ring. Then  $S$ is pr\'{e}simplifiable if $R$  is pr\'{e}simplifiable and $N$   is $R$-pr\'{e}simplifiable.
\end{enumerate}
\end {prop}
\pr 1. Let $H$ be an $R$-submodule of $S$. Consider $e=xe$ with $e\in H-\{0\}$ and $x\in R-\{0\}$. Since $S$ is pr\'{e}simplifiable, $x\in U(S)$ and so $x\in U(S)\cap R=U(R)$.\\
2.  Let $x=r_{x}+n_{x}\neq 0$ and $y=r_{y}+n_{y}$ be two elements of $R\oplus N=S$ where $r_x, r_y\in R$ and $n_x, n_y\in N$, such that $x=yx$. Assume that $r_x\neq 0$. Then $r_x=r_yr_x$ implies that $r_y\in U(R)\subseteq U(S)$, and, since $N$ is nilpotent, $y=r_y+n_y$ is invertible in $S$, as desired. Next, assume now that $r_x=0$. Then $n_x\neq 0$ and so  $n_x=r_yn_x+n_yn_x$. In the case $N^{2}=0$, we have $n_x=r_yn_x$. Hence $r_y\in U(R)$ since $N$ is pr\'{e}simplifiable, and  as above $y\in U(S)$. Finally, in the case where 
$S=R\oplus   N_1  \oplus N_2 \oplus \cdots $ is a graded ring, we may set $n_x=n_{i_1}+ \cdots +n_{i_m}$ with $\{i_1,...,i_m\}\subset \N$ and $m\in \N$ such that $i_1\leq \cdots \leq i_m$ and $n_{i_1}\neq 0$. Then $n_{i_1}=r_yn_{i_1}$  which  implies that $r_y\in U(R)$ and similarly as above $y\in U(S)$.\cqfd

\begin{prop}\label{propo-str-ass}
Let $R=\underset{i \in \N_0}{\oplus}   R_i $ be a graded ring.
\begin{enumerate}
  \item If $R$ is strongly associate, then $R_0$ is a strongly associate ring and $R_i$ is a strongly associate $R_0$-module for every $i\in \N$.
  \item  Suppose   there exists $n\in \N$ such that $R_i=0$ for every $i\geq n+1$, that is, $R=R_0\ltimes_n R_1\ltimes \cdots \ltimes R_n$, and assume that $R_0$ is a  pr\'{e}simplifiable ring and $R_1, ... , R_{n-1}$ are pr\'{e}simplifiable $R_0$-modules. Then $R$ is strongly associate if and only if $R_n$ is strongly associate.
\end{enumerate}
\end{prop}
\pr  1. Let $x_i, y_i\in R_i-\{0\}$ for $i\in \N_0$ such that $R_0x_i=R_0y_i$. Then $Rx_i=Ry_i$. Hence there is $u=u_0+u_1+ \cdots \in U(R)$ such that $x_i=uy_i$. Then $u_0\in U(R_0)$ and $x_i=u_0y_i$, as desired.\\
2. Let $x=x_m+\cdots+x_n$ and $y=y_m+ \cdots +y_n$ be  two associate elements of $R$ where $m\in \{0,...,n\}$ and $x_i, y_i\in R_i$ for $i\in \{m,...,n\}$ such that $x_m$ and $y_m$ are nonzero. Then $x_m\sim y_m$. In particular, there is $\alpha=\alpha_0+ \cdots +\alpha_n$ such that $x=\alpha y$. Then $x_m=\alpha_0y_m$. Hence two cases occur. Case $m\neq n$. Since $R_m$ is pr\'{e}simplifiable, $\alpha_0\in U(R_0)$. Then $\alpha\in U(R)$, as desired. Case $m=n$ (i.e., $x=x_m$ and $y=y_m$). Here, the result follows since $R_n$ is strongly associate.\cqfd

Now we can give the  extension of \cite[Theorem 5.1]{AW} to the context of $n$-trivial extensions.

\begin{cor}\label{cor-pres-str-ass}
The following assertions are true.
\begin{enumerate}
  \item $R\ltimes_n M_{1}\ltimes \cdots \ltimes M_{n}$ is pr\'{e}simplifiable if and only if $R$, $M_1,...,M_n$ are pr\'{e}simplifiable.
  \item If $R\ltimes_n M_{1}\ltimes \cdots \ltimes M_{n}$ is strongly associate, then $R$, $M_1, ... , M_n$ are strongly associate.
  \item Suppose that $R$, $M_1 , ... ,  M_n$ are pr\'{e}simplifiable. Then $R\ltimes_n M_{1}\ltimes \cdots \ltimes M_{n}$ is strongly associate if and only if $M_n$ is strongly associate.
\end{enumerate}
\end {cor}

Now we investigate the extension of  \cite[Theorem 5.4]{AW}. It is  convenient to   recall the following definitions.  Let $S$  be a commutative ring. A nonunit $a\in S$ is  said to be \textit{irreducible} or an \textit{atom} (resp., \textit{strongly irreducible}, \textit{very strongly irreducible}) if $a = bc$ implies $ a \sim b $  or $ a\sim c $ (resp., $ a \approx b $ 
 or $ a \approx c$, $a \cong b $ or $ a\cong c$) and $a$ is said to be \textit{$m$-irreducible} if $Sa$ is a maximal element of the set of proper principal ideals of $S$. Note that, for a nonzero nonunit $a\in S$, $a$ very strongly irreducible $\Rightarrow$ $a$ is $m$-irreducible $\Rightarrow$  $a$ is strongly irreducible $\Rightarrow$ $a$ is irreducible, but none of these implications can be reversed. In the case of an $S$-module $H$, we say that $ e \in H$ is \textit{$S$-primitive} (resp., \textit{strongly $S$-primitive}, \textit{very strongly $S$-primitive}) if for $a\in S$ and $ f\in  H$, $e = af$ $\Rightarrow$ $ e \sim f$ (resp., $e \approx f$, $e \cong f$). And $e$ is \textit{$S$-superprimitive} if $be = af$  for $a,b \in S$ and $f\in H$, implies $a \mid b$. Note that (1) $e$ is $S$-primitive $\Leftrightarrow $ $Se$ is a maximal cyclic $S$-submodule of $H$, (2) $e$ is $S$-superprimitive $\Rightarrow$ $e$ is very strongly $S$-primitive $\Rightarrow$ $e$ is strongly $S$-primitive $\Rightarrow$ $e$ is $S$-primitive, (3) if $Ann(e) = 0$, $e$ is $S$-primitive $\Rightarrow$ $e$ is very strongly $S$-primitive, and (4) $e$ is $S$-superprimitive $\Rightarrow$ $Ann(e) = 0$.\medskip

In the following results the homogeneous element $(0,...,0,m_i,0,...,0)\in R\ltimes_n M$ where $i\in \{1, ... ,n\}$ and $m_i\in M_i-\{0\}$, is denoted  by $ \underline{m_i}$. The following result extends  \cite[Theorem 5.4 (1)]{AW}.

\begin{prop} \label{pro-asso}
  Let $i\in \{1, ... ,n\}$ and $m_i, n_i\in M_i-\{0\}$. Then $m_i\sim n_i$ (resp., $m_i\approx n_i$, $m_i\cong n_i$) in $M_i$ if and only if $ \underline{m_i}\sim \underline{n_i}$ (resp., $ \underline{m_i}\approx \underline{n_i}$, $ \underline{m_i}\cong \underline{n_i}$) in $R\ltimes_n M$.
     \end {prop}
\pr The assertion   is  proved similarly to the corresponding classical one.\cqfd

It is worth noting that the analogue of the assertion $(4)$ of \cite[Theorem 5.4]{AW} does not hold   in the context of $n$-trivial extensions with $n\geq 2$. Indeed, consider the $2$-trivial extension $S=\Z_4\ltimes_2 \Z_4\ltimes \Z_4$. It is easy to show that $\bar{1}$ is superprimitive  in the $\Z_4$-module $\Z_4$. However, $(\bar{0},\bar{0},\bar{1})$ is not very strongly irreducible in $S$ (since $(\bar{2},\bar{1},\bar{2})^{2}=(\bar{0},\bar{0},\bar{1}))$. Moreover, even if we assume that $R$ is an integral domain, we still don't have  the desired analogue. For this, take $S=\Z\ltimes_2 \Z\ltimes \Z$. We have $1$ is superprimitive in the $\Z$-module $\Z$. However, $(0,0,1)$ is not very strongly irreducible in $S$ (since $(0,1,0)^{2}=(0,0,1)$). The last example also shows that  the assertion $(2)$ of \cite[Theorem 5.4]{AW} does not hold   in the context of $n$-trivial extensions with $n\geq 2$. 
Namely, if $0\not =m_i= m_jm_k$ where $ (m_i,m_j,m_k)\in M_i\times M_j\times M_k$,  $i\geq 2$ and  $j,k\in \{1,...,i-1\}$ with $j+k=i$,  then $\underline{m_i}$  cannot be irreducible. Indeed, $\underline{m_i}= \underline{m_j} \underline{m_k}$  but neither $\underline{m_j}$ nor  $ \underline{m_k}$  are in  $\langle\underline{m_i}\rangle\subseteq 0\ltimes_n    \cdots\ltimes 0 \ltimes    M_{i} \ltimes\cdots  \ltimes M_n  $. \\
\indent To extend \cite[Theorem 5.1 (2)]{AW}, we need to introduce the following definitions.

\begin{defn}\label{def-indecom} Assume $n\geq 2$ and each multiplication in the family $ \varphi=\{\varphi_{i,j}\}_{ \underset{1\leq i,j\leq n-1}{i+j\leq n}}$ is not trivial.   Let $i\in \{2, ... ,n\}$. An element $m_i \in M_i-\{0\}$ is said to be  $\varphi$-indecomposable  (or indecomposable relative to the family of multiplications $ \varphi$) if,  for every $ (m_j,m_k)\in M_j\times M_k$ (where $j,k\in \{1,...,i-1\}$ with $j+k=i$),   $m_i\not= m_jm_k$. If  no ambiguity can arise,   $\varphi$-indecomposable elements are simply called indecomposables.
\end{defn}

For example, in $ \Z \ltimes_2    \Z  \ltimes \Q  $,   every element in $\Q -\Z$ is indecomposable. However, every element   $x\in\Z$ ($\Z$ as a submodule of $\Q$) is decomposable (since  $(0,1,0)(0,x,0)=(0,0,x)$).

 \begin{defn}\label{def-phi-integr}
 Let $i\in \{2, ... ,n\}$. The $R$-module $M_i$ is said to be $\varphi$-integral (or integral relative  to the family of multiplications $ \varphi$)  if,  for every $ (m_j,m_k)\in M_j\times M_k$ (where $j,k\in \{1,...,i-1\}$ with $j+k=i$),   $ m_jm_k=0$ implies that $m_j=0$ or $m_k=0$. If  no ambiguity can arise, a  $\varphi$-integral $R$-module is simply called integral.
\end{defn}

For example,  for  $\Z \ltimes_2 \Z  \ltimes   \Z  $, $M_2=\Z $ is integral. And, for  $\Z \ltimes_2 \Z  \ltimes   \Z /2 \Z  $, $  \Z /2 \Z $ is not integral since, for instance, 
$\varphi_{1,1}(1,2)=\overline{1}\; \overline{2}=\overline{0}$.

\begin{prop} \label{pro-irred} Assume $n\geq 2$.  Let $i\in \{1, ... ,n\}$ and $m_i\in M_i-\{0\}$.
 If $ \underline{m_i}$ is irreducible (resp., strongly irreducible, very strongly irreducible) in $R\ltimes_n M$, then $m_i$ is primitive (resp., strongly  primitive, very strongly  primitive) in $M_i$.\\
\indent Conversely, three cases occur:
\begin{description}
    \item[Case $i=1$.]  The reverse implication holds if  $R$ is an integral domain and $M_j$ is torsion-free for every $j\in \{2, ... ,n\}$.
 \item[Case $i=2$] (here $n\geq 2$). The reverse implication holds   if  $R$ is an integral domain, $M_j$ is torsion-free for every $j\in \{1, ... ,n\} - \{2\}$ and  $m_2$ is indecomposable.
 \item[Case $i\geq 3$]  (here $n\geq 3$). The reverse implication holds   if   $R$ is an integral domain, $M_j$ is torsion-free for every $j\in \{1, ... ,n\} - \{i\}$, $M_j$ is integral for every $j\in \{2, ... ,i-1\}$, and $m_i$ is indecomposable.
\end{description} 
\end {prop}
\pr We only prove  the  primitive (irreducible) case. The two other cases are proved similarly.\\
\indent $\Longrightarrow  $ Suppose that $\underline{m_i}$ is irreducible and let $m_i=an_i$ for some $a\in R$ and $n_i\in M_i$. Then $\underline{m_i}=(a,0,...,0)\underline{n_i}$ and so $(R\ltimes_n M) \underline{m_i}=(R\ltimes_n M) \underline{n_i}$. This implies that $Rm_i=Rn_i$, as desired.\\
\indent $\Longleftarrow$ Let $\underline{m_i}=(a_j)(n_j)$ for some $(a_j),(n_j)\in R\ltimes_n M$. Then $a_0n_0=0$. First, we show that the case $a_0=0$ and $n_0=0$ is impossible.  Cases $i=1,\,2$ are easy and are left to the reader. So  assume $i\geq 3$. Suppose that $a_0=0$ and $n_0=0$. Then we have the following equalities:
$$\mathrm{For} \ j\in \{2, ... ,i-1\}, \ a_1n_{j-1}+a_2n_{j-2}+ \cdots +a_{j-1}n_1=0\ \mathrm{and}\  a_1n_{i-1}+a_2n_{i-2}+ \cdots +a_{i-1}n_1=m_i.$$  
A recursive argument on these equalities shows that, for $l\in \{2, ... ,i\}$, there is $k\in \{0, ... ,l-1\}$ such that $(a_0,...,a_k)=(0,...,0)$ and $(n_0,...,n_{l-(k+1)})=(0,...,0)$. Indeed, it is clear  this is true for $l=2$. Then suppose this is true for a given  $l\in \{2, ... ,i-1\}$. So the equality $a_1n_{l-1}+a_2n_{l-2}+ \cdots +a_{l-1}n_1=0$ becomes $a_{k+1}n_{l-(k+1)}=0$. Then  since $M_l$ is integral, we get the desired result for $l+1$. Thus for $l=i$, we get $a_{k+1}n_{i-(k+1)}=m_i$ which is absurd since $m_i$ is indecomposable.  \\
\indent Now, we may assume that $a_0\neq 0$ and $n_0=0$, so $a_0n_1=0$. Since $M_1$ is torsion-free, $n_1=0$. Recursively we get $n_j=0$ for $j\in \{1, ... ,i-1\}$. Then $a_0n_i=m_i$, and since $m_i$ is primitive,  there exists $b_0\in R$ such that $n_i=b_0m_i$. It remains to show that there is $b_j\in M_j$ for $j\in \{1, ... ,n-i\}$ such that $n_{i+j}=b_jm_i$ and this implies that $(n_i)=(0,...,0,n_i,n_{i+1},...)=(b_0,...,b_{n-i},0,...,0)\underline{m_i}$, as desired. We have $a_0n_{i+1}+a_1n_i=0$. Then using both $a_0n_i=m_i$ and $n_i=b_0m_i$, we get $a_0n_{i+1}+a_1b_0a_0n_i=0$. Then $a_0(n_{i+1}+a_1b_0n_i)=0$, so $n_{i+1}+a_1b_0n_i=0$ (since $M_{i+1}$ is torsion-free). Then $n_{i+1}=-a_1b^{2}_0m_i$. Then we set $b_1=-a_1b^{2}_0$ and so $n_{i+1}=b_1 m_i$. Similarly, using the equality $a_0n_{i+2}+a_1n_{i+1}+a_2n_i=0$ with the   equalities  $a_0n_i=m_i$ and $n_i=b_0m_i$, we get   $a_0n_{i+2}+a_1b_1a_0n_{i}+a_2b_0a_0n_i=0$. So $n_{i+2}+a_1b_1n_i+a_2b_0n_i=0$, then $n_{i+2}=b_2m_i$ where $b_2=-a_1b_1b_0-a_2b^{2}_0$. Finally, a recursive argument gives the desired result.  \cqfd\medskip

The following result extends  \cite[Theorem 5.4 (3)]{AW}.

\begin{prop} \label{pro-idem-irre}
    Suppose that $R$ has a nontrivial idempotent. Then for every $i\in \{1, ... ,n\}$ and $m_i\in M_i-\{0\}$, $ \underline{m_i}$ is not irreducible in $R\ltimes_n M$.      \end {prop}
\pr The assertion   is  proved similarly to the corresponding classical one.\cqfd

Now we are interested in some factorization properties. Recall that a   ring $S$ is called \textit{atomic}    if every (nonzero) nonunit of $S$ is a product of irreducible  elements  (atoms) of $S$. Note that, as in the domain case, the ascending chain condition on principal ideals (ACCP) implies atomic.\\
\indent  We begin with an extension of \cite[Theorem 5.5 (2)]{AW} which characterizes  when a trivial extension of a ring satisfies ACCP. For this, we need the following lemma.

\begin {lem} \label{lem-Acc}
  Let $i\in \{0,... , n\}$ and consider two elements $a=(0,...,0,a_i,a_{i+1},...,a_n)$ and  $b=(0,...,0,b_i,b_{i+1},...,b_n)$ of $R\ltimes_n M$ with $a_i\neq 0$. Then the implication  ``$\langle a \rangle\subsetneq \langle b \rangle$ $\Rightarrow$ $b_i\neq 0$ and $\langle a_i \rangle\subsetneq \langle b_i \rangle$" is true if either $(1)$   $0\leq i\leq n-1$ and $M_i$ is pr\'{e}simplifiable (here $M_0=R$) or $(2)$   $i=n$. 
\end{lem}
\pr
 Since $\langle a \rangle\subsetneq \langle b \rangle$, there is $c=(c_0,...,c_n)\in R\ltimes_n M - U(R\ltimes_n M)$ such that $a=cb$. Then $a_i=c_0b_i$ and $c_0\notin U(R)$.
 This shows that  $\langle a_i \rangle\subsetneq \langle b_i \rangle$ in both cases.  \cqfd

 \begin{thm}\label{thm-ACCP}
 Assume $n\geq 2$. Suppose that $M_i$ is   pr\'{e}simplifiable  for every $i\in \{0,...,n-1\}$ (here $M_0=R$). Then  $R\ltimes_n M$ satisfies ACCP if and only if $R$ satisfies ACCP and, for every $i\in \{1,... , n\}$, $M_i$ satisfies ACC on cyclic submodules.
 \end{thm}
 \pr The proof of the direct implication is easy. Let us prove the converse. Suppose that $R\ltimes_n M$ admits a strictly ascending chain of principal ideals
 $$\langle (a_{1,i}) \rangle\subsetneq \langle( a_{2,i}) \rangle \subsetneq \cdots$$
 If there exists $j_0\in \N$ such that $a_{j_0,0}\not= 0$. Then for every $k\geq j_0$, $a_{k,0}\neq 0$. Then by Lemma \ref{lem-Acc}, we get the following  strictly ascending chain of principal ideals of $R$:
 $$\langle a_{ j_0,0} \rangle\subsetneq \langle a_{j_0+1,0} \rangle \subsetneq \cdots$$
 This is absurd since $R$ satisfies ACCP. Now, suppose that $a_{ j,0}=0$ for every $j\in \N$ and that there exists $j_1\in \N$ such that  $a_{j_1,1}\neq 0$. Also, by Lemma \ref{lem-Acc}, we obtain the following strictly ascending chain of cyclic submodules of $M_1$:
 $$\langle a_{ j_1,1} \rangle\subsetneq \langle a_{ j_1+1,1} \rangle \subsetneq \cdots$$
 This is absurd since $M_1$ satisfies ACC on cyclic submodules. We continue in this way until the case where we may suppose that $a_{j,i}=0$ for every $i\in \{1,...,n-1\}$ and every $j\in \N$. Therefore, by  Lemma \ref{lem-Acc}, we get the desired result.\cqfd\medskip

Now, we investigate when $R\ltimes_n M$ is atomic. Namely, we give an extension of \cite[Theorem 5.5 (4)]{AW}. Recall that an $R$-module $N$ is said to satisfy  MCC  if every cyclic submodule of $N$ is contained in a maximal (not necessarily
proper) cyclic submodule of $N$.

\begin{thm}\label{thm-atomic}
 Assume $n\geq 2$. Suppose that $M_i$ is   pr\'{e}simplifiable  for every $i\in \{0,...,n-1\}$ (here $M_0=R$). Then $R\ltimes_n M$ is atomic if $R$ satisfies ACCP, $M_i$ satisfies ACC on cyclic submodules, for every $i\in \{1, ... , n-1\}$, and $M_n$ satisfies MCC.
\end{thm}
 \pr The proof is slightly more technical  than the one of \cite[Theorem 5.5 (4)]{AW}. Here, we need to break the proof into the following $n+1$ steps  such that in the step number $k\in \{1, ... , n+1\}$ we prove that  every nonunit element  $(m_i)\in R\ltimes_n M$ with   $m_0=0$,..., $m_{k-2}= 0$  and  $m_{k-1}\neq 0$  is a product of irreducibles.   \smallskip

We use an inductive argument for the first $n$ steps.\\
              \indent \underline{Step 1}.   Suppose there is a nonunit element $(m_i)$ of $R\ltimes_n M$ with $m_0\neq 0$ and such that $(m_i)$ cannot be factored  into irreducibles.
 Then there exist $(a_{1,i})$, $(b_{1,i})\in R\ltimes_n M-U(R\ltimes_n M)$ such that $(m_i)=(a_{1,i})(b_{1,i})$ and neither $(m_i)$ and $(a_{1,i})$ nor  $(m_i)$ and $(b_{1,i})$ are associate.  Since $0\neq m_0=a_{1,0}b_{1,0}$,  $a_{1,0}\neq 0$ and $b_{1,0}\neq 0$. Clearly $(a_{1,i})$  or  $(b_{1,i})$ must be reducible, say $(a_{1,i})$. Also, for $(a_{1,i})$ there are $(a_{2,i})$, $(b_{2,i})\in R\ltimes_n M-U(R\ltimes_n M)$ such that   $(a_{1,i})=(a_{2,i})(b_{2,i})$  and neither $(a_{1,i})$ and $(a_{2,i})$ nor  $(a_{1,i})$ and $(b_{2,i})$ are associate. As above, $a_{2,0}\neq 0$ and $b_{2,0}\neq 0$ and say $(a_{2,i})$ is reducible. So we continue and then we obtain a strictly ascending chain  
$$\langle (m_i) \rangle\subsetneq \langle (a_{1,i}) \rangle \subsetneq \langle(a_{2,i})\rangle \subsetneq \cdots   $$
Using Lemma \ref {lem-Acc}, we get a strictly ascending chain  of principal ideals of $R$ 
$$\langle m_0 \rangle\subsetneq \langle a_{1,0} \rangle \subsetneq \langle a_{2,0}\rangle \subsetneq \cdots   $$
This is absurd since $R$ satisfies ACCP.\smallskip
 
\indent \underline{Step $j$} ($1\leq j  \leq n$). Suppose there  is a nonunit element  $(m_i)\in R\ltimes_n M$ with  $m_0=0$,..., $m_{j-2}= 0$  and  $m_{j-1}\neq 0$ that is not a product of irreducibles. Then  there   are $(a_{1,i})$, $(b_{1,i})\in R\ltimes_n M-U(R\ltimes_n M)$ such that $(m_i)=(a_{1,i})(b_{1,i})$  and neither $(m_i)$ and $(a_{1,i})$ nor  $(m_i)$ and $(b_{1,i})$ are associate.     Then $a_{1,0}b_{1,j-1}+a_{1,1}b_{1,j-2}+\cdots +a_{1,j-2}b_{1,1} +a_{1,j-1}b_{1,0}=m_{j-1}\neq 0$. If $a_{1,k}=0$ for every $k\in \{0,...,j-2\}$, then   necessarily $b_{1,0}\neq 0$. Hence by the preceding steps, $(b_{1,i})$ is a  product of irreducibles and then by hypothesis on $(m_i)$, $(a_{1,i})$ is reducible. If  $a_{1,k}\not=0$ for some $k\in \{0,...,j-2\}$, then  $(a_{1,i})$ is a  product of irreducibles and $(b_{1,i})$ is reducible. Thus, by symmetry, we may assume that  $(a_{1,i})$ is reducible and it is not a   product of irreducibles. So, necessarily $a_{1,0}=0$,..., $a_{1,j-2}= 0$  and  $a_{1,j-1}\neq 0$.  We repeat the last argument so that we obtain a strictly ascending chain  of principal ideals of $R\ltimes_n M$
$$\langle( m_i) \rangle\subsetneq \langle (a_{1,i}) \rangle \subsetneq \langle (a_{2,i})\rangle \subsetneq \cdots   $$  such that, for every $k\in \N$,  $a_{k,0}=0$,..., $a_{k,j-2}= 0$  and  $a_{k,j-1}\neq 0$. Then, by Lemma \ref{lem-Acc}, we get a strictly ascending chain  of cyclic submodules of $M_{j-1}$
$$\langle  m_{j-1} \rangle\subsetneq \langle  a_{1,{j-1}}  \rangle \subsetneq \langle  a_{2,{j-1}} \rangle \subsetneq \cdots   $$  which is absurd by hypothesis on $M_{j-1}$, as desired.\smallskip

                 \indent \underline{Step $n+1$}.   It remains to prove that every element of $R\ltimes_n M$ of the form $(0,... ,0,m_n)$ with $m_n\neq 0$ is a product of irreducibles. Since $M_n$ satisfies MCC, $Rm_n\subseteq Rm$ where $Rm$ is a maximal cyclic submodule of $M_n$. Then $m_n=am$ for some $a\in R-\{0\}$ and so $(0,... ,0,m_n)=(a,0, ... ,0)(0,... ,0,m)$. Now, $a\neq 0$  shows that $(a,0,... ,0)$ is a product of irreducibles (by Step 1) and $Rm$ is maximal shows that either $(0,... ,0,m)$ is irreducible   or $(0,... ,0,m)=(a_{i})(b_{i})$ where $a_k\neq 0$ and $b_l\neq 0$ for some  $k,l\in \{0,...,n-1\}$. Then   by  the preceding steps, $(a_{i})$ and $(b_{i})$ are products of irreducibles and hence so is  $(0,... ,0,m)$. This concludes the proof.\cqfd\medskip

 A ring $S$ is   said to be a  \textit{bounded factorial ring} (BFR)  if, for each nonzero nonunit $x \in S$, there is a natural number $N(x)$ so that for any factorization $x = x_{1} \cdots x_{s}$ where each $x_{i}$ is a nonunit, we have $s \leq N(x)$. For domains we say   BFD  instead of   BFR. Recall that an $S$-module $H$ is said to be  a \textit{BF-module} if, for each nonzero $h\in H$, there exists a natural number $N(h)$ so that $h=a_{1}\cdots a_{s-1}h_{s}$  (each $a_i$ a nonunit) $\Rightarrow s\leq N(h)$.\\
\indent Our next theorem, which is a generalization  of \cite[Theorem 5.5 (4)]{AW}, investigates when $R\ltimes_n M$ is BFR. It is based on the following lemma.

\begin{lem}\label{lem-BFR} For $j\in \N-\{1\}$, a product of $j$ elements of $R\ltimes_n M$ of the form $(0,x_1,...,x_n)$ is   of the form  $(0,...,0,y_{j},...,y_n)$ (where, if $j\geq n+1$ the product is zero). 
 \end{lem}

\begin{thm}\label{thm-BFR}  
 Assume that $n\geq 2$, $R$ is an integral domain and $M_i$ is torsion-free for every $i\in \{1,..., n-1\}$. Then  $R\ltimes_n M$ is a BFR  if and only if  $R$ is a BFD and $M_i$ is a BF-module for every $i\in \{1,..., n\}$.
 \end{thm}
 \pr
$\Longrightarrow $ Clear.\\
$\Longleftarrow  $ Let $(m_i) $ be a nonzero nonunit element of  $R\ltimes_n M$ and suppose we have a factorization into nonunits $(m_i)=(a_{1,i}) \cdots (a_{s,i})$ for some $s\in \N$.  
If $m_0\neq 0$, $m_0=a_{1,0} \cdots a_{s,0}$ implies that $s\leq N(m_0)$. Otherwise, there is $j\in \{1,..., n\}$ such that $m_0=0$,..., $m_{j-1}= 0$  and  $m_{j}\neq 0$.  We may assume that $s\geq j+1$. Since $R$ is an integral domain and by Lemma \ref{lem-BFR}, we  may assume   there is $k\in\{1,..., j\}$ such that $a_{l,0} = 0$ for every  $l\in \{1,...,k\}$ and  $a_{l,0} \not = 0$ for every  $l\in\{k+1,..., s \}$. Let  $(0,...,0,b_k,...,b_{ n})=\overset{k}{\underset{l=1}{ \prod}} (a_{l,i})$ and $(c_0,...,c_n)=\overset{s}{\underset{l=k+1}{ \prod}} (a_{l,i})$. Since   $M_i$ is torsion-free for every $i\in \{1,..., j-1\}$ and $c_0= \overset{s}{\underset{l=k+1}{ \prod}} a_{l,0} \neq 0$,   $b_k=0$,...,$b_{j-1}=0$  and $b_j \neq 0$. 
Then  $m_{j}=c_0 b_j= \overset{s}{\underset{l=k+1}{ \prod}} a_{l,0} b_j $. Therefore  $s\leq N(m_j) +k-1$ (since $M_j$ is a BF-module).\cqfd \medskip

Now we investigate the  notion of a $U$-factorization. It was    introduced by Fletcher  \cite{F1,F2}      and   developed  by Axtell et al.  in \cite{Ax} and \cite{AxFRS}. Let $S$ be a ring and consider a nonunit $ a \in S$. By a \textit{factorization} of $a$ we mean $ a = a_{1} \cdots a_{s}$ where each
$a_{i}$ is a nonunit.  Recall from  \cite{F1} that, for $ a \in S$, $U(a) = \{ r \in S \mid \exists s \in S$ with $rsa =
a \} = \{ r \in S \mid r(a) = (a) \}$. A \textit{$U$-factorization} of $a$ is a factorization $ a = a_{1} \cdots a_{s}b_{1} \cdots b_{t}$
where,  for every $ 1 \leq i \leq s $, $ a_{i}\in U(b_{1} \cdots b_{t})$ and,  for every     $ 1 \leq i \leq t $, $ b_{i} \notin U(b_{1} \cdots \hat{b_{i}} \cdots b_{t})$. We denote this $U$-factorization by $a = a_{1}\cdots a_{s} \lceil b_{1} \cdots b_{t} \rceil$ and call $ a_{1},$...,$ a_{s}$ (resp.,  $b_{1},$...,$ b_{t}$)  
the \textit{irrelevant} (resp.,  the \textit{relevant}) factors.\smallskip

\indent Our next result investigates when an $n$-trivial extension is a $U$-FFR.  First, recall the following definitions.\smallskip

\indent A ring $S$ is called a \textit{finite factorization} ring (FFR) (resp., a \textit{$U$-finite factorization} ring ($U$-FFR)) if every  nonzero nonunit of $S$ has only a finite number of factorizations (resp., $U$-factorizations)  up to order and associates (resp., associates on the relevant factors).    A ring $S$ is called a \textit{weak finite factorization} ring  (WFFR) (resp., a \textit{$U$-weak finite factorization} ring  ($U$-WFFR)) if  every nonzero nonunit of $R$ has only a finite number of nonassociate divisors (resp., nonassociate relevant factors). We have FFR $\Rightarrow$ WFFR and the converse holds  in the domain case. But $ \mathbb{Z}_{2} \times \mathbb{Z}_{2}$ is a WFFR that is not an FFR. However,  from \cite[Theorem 2.9]{Ax}, $U$-FFR $\Leftrightarrow$ $U$-WFFR.\smallskip

\indent  The study of the  notions above on the classical trivial extensions has lead to consider the following notion (see \cite{Ax}).  Let $N$ be an $S$-module.  For a nonzero element $x\in N$, we say that $Sd_1d_2 \cdots d_s x$ is a \textit{reduced submodule factorization} if, for every $j\in \{1,...,s\}$,  $d_j\not \in U(S)$    and for no cancelling and reordering of the $d_j$'s is it the case that $Sd_1d_2 \cdots d_s x= Sd_1d_2 \cdots d_t x$   where $t < s$. The module $N$ is said   to be a\textit{ $U$-FF} module if for every nonzero element $x\in N$, there exist only finitely many reduced submodule factorizations $Sx= Sd_1d_2 \cdots d_t x_k$, up to order and associates on the $d_i$, as well as up to associates on the $ x_k$.  In our context, we introduce the following definition.

\begin{defn}\label{def-U-FFM} Assume $n\geq 2$ and consider $i\in \{1,...,n\}$.
\begin{enumerate}
    \item  Let  $m_i\in M_i-\{0\}$, $s \in \N $   and $(d_{i_1},  ..., d_{i_s})\in M_{i_1} \times \cdots \times M_{i_s}$ where $\{i_1,i_2 ,...,  i_s\}\subseteq  \{0,...,n\}$ with  $i_1+\cdots + i_s=i$.  We say that $Rd_{i_1}d_{i_2} \cdots d_{i_s}m_i\subseteq M_i$ is a $\varphi$-reduced submodule factorization if, for every $j\in \{1,...,s\}$ such that $i_j=0$,   $d_{i_j}\not\in U(R)$    and for no cancelling and reordering of the $d_j$'s is it the case that $Rd_{i_1}d_{i_2} \cdots d_{i_s}= Rd_{i_1}d_{i_2} \cdots d_{i_t}$   where $t < s$. If  no ambiguity can arise, a $\varphi$-reduced submodule factorization is simply called a reduced submodule factorization.
    \item  The $R$-module $M_i$ is said  to be a $\varphi$-$U$-FF module (or simply $U$-FF module) if, for every nonzero element $x\in  M_i$, there exist only finitely many reduced submodule factorizations $Rx= Rd_{i_1}d_{i_2} \cdots d_{i_s}$, up to order and associates on the $d_{i_j}$. 
\end{enumerate}
\end{defn}

It is clear that, for $i=1$,  the notion of $U$-FF module defined here is the same as the Axtell's one.\smallskip

Based on the proof of \cite[Theorem 4.2]{Ax}, it is asserted in \cite[Theorem 3.6]{AxFRS} that, if $R\ltimes_1 M_1$  is  a $U$-FFR, then for every nonzero  nonunit $d \in R$, there are only finitely many distinct principal ideals $\langle(d,m)\rangle$ in  $R\ltimes_1 M_1$. However, a careful reading of this proof shows that the case of ideals $\langle(d,m)\rangle$ with $dM_1=0$ should be also treated. One can confirm the validity of this assertion for reduced rings. However,  the context of $n$-trivial extensions seems to be more complicated. Nevertheless, under some certain conditions, we next investegate  when $R\ltimes_n M$ is  a $U$-FFR.

\begin{lem} \label{lem-U-FFR}
Assume $n\geq 2$ and $M_n$ is  integral. Then for every nonzero  nonunit $d \in R$, the following   assertions are true.
\begin{enumerate}
    \item For every $i\in  \{1,..., n-1\}$, the following assertions are equivalent:
\begin{itemize}
    \item[1.a.] $dM_i=0$.
\item[1.b.]  $dm_i=0$ for some $m_i \in M_i-\{0\}$.
\item[1.c.]  $dM_{n-i}=0$.
\item[1.d.]  $dm_{n-i}=0$ for some $m_{n-i} \in M_{n-i}-\{0\}$.
\end{itemize}
    \item  The following assertions are equivalent:
\begin{itemize}
    \item[2.a.] $dM_i=0$ for some  $i\in  \{1,..., n-1\}$.
\item[2.b.] $dM_i=0$ for every   $i\in  \{1,..., n-1\}$. 
\end{itemize}
   \item If $dM_n=0$, then $dM_i=0$ for every   $i\in  \{1,..., n-1\}$. 
  \item If $M_n$ is torsion-free, then  $M_i$ is torsion-free for every   $i\in  \{1,..., n-1\}$. 
\end{enumerate}
\end{lem}
\pr $(1)$. For the implications  $(1.a) \Rightarrow (1.b)$ and $(1.c) \Rightarrow (1.d)$ there is nothing to prove.\\
 \indent $(1.b) \Rightarrow (1.c)$. Let $m\in  M_{n-i}$. Then $d m_i m =0\in M_n$. Therefore $dm=0$ (since $M_n$ is integral and $m_i\neq 0$).\\
 \indent $(1.d) \Rightarrow (1.a)$. Similar to the previous proof.\\
\indent $(2)$.  For the implication   $(2.b) \Rightarrow (2.a)$ there is nothing to prove.\\
\indent $(2.a) \Rightarrow (2.b)$. First, we  prove that $dM_1=0$. For every $m_1\in M_1-\{0\}$, $dm_1^i=0\in M_i$ and so  $dm_1^n=0\in M_n$. Therefore $dm_1=0$ (since  $M_n$ is integral and $m_1\neq 0$). Now, consider any $j\in \{1,..., n-1\}$ and any  $m_j\in M_j-\{0\}$. Then for every $m_1\in M_1-\{0\}$, $dm_jm_1^{n-j}=0\in M_n$ which shows that  $dm_j=0$.\\
\indent $(3)$. This is proved as above.\\
\indent $(4)$.  If there is $m_1 \in M_1-\{0\}$     and $r \in R -\{0\}$ such that $r m_1 =0$, then  $r m_1^n =0\in M_n$. Since $M_n$ is torsion-free and $r\not =0$, $ m_1^n =0\in M_n$ so  $ m_1 =0$ (since $M_n$ is integral). This is absurd since  $m_1\not =0$. Finally, by assertions $(1)$ and $(2)$, we conclude that   $M_i$ is torsion-free for every   $i\in  \{1,..., n-1\}$. \cqfd

\begin{thm}\label{thm-U-FFR} Assume  $n\geq 2$.  If $R\ltimes_n M$ is  a $U$-FFR   (equivalently, a $U$-WFFR), then    the following   conditions are satisfied:
\begin{enumerate}
    \item $R$ is an FFR.
    \item $M_i$ is a $U$-FF module for every  $i\in \{1,..., n\}$.
\end{enumerate}
 Moreover, if $R$ is  an integral domain and $M_n$ is  integral and torsion-free, then 
\begin{itemize}
    \item[3.] For every nonzero  nonunit $d \in R$, there are only finitely many distinct principal ideals $\langle(d,m_1,...,m_n)\rangle$ in  $R\ltimes_n M$.
\item[4.]  For every  $i\in \{1,..., n-1\}$ and every   $m \in M_i - \{ 0\}$, there are only finitely many distinct principal ideals $\langle(0,...,0,m,m_{i+1},...,m_n)\rangle$  in  $R\ltimes_n M$.
\end{itemize}
\indent Conversely, if $R$ is  an integral domain and $M_n$ is  integral and torsion-free, then the assertions $(1)$-$(4)$ imply that $R\ltimes_n M$ is  a $U$-FFR.
\end{thm} 
\pr The proof of the ``converse" part is similar to the corresponding one  of \cite[Theorem 4.2]{Ax}.\\
\indent   $\Longrightarrow  $ The proof of each $(1)$ and $(2)$ is similar to that given in \cite[Theorem 4.2]{Ax}.\smallskip

\indent $3.$ Suppose, by contradiction, there exists a nonzero  nonunit $d \in R$    for which there is a family of distinct principal ideals of the form $\langle(d,m_{j,1},...,m_{j,n})\rangle$ where $j$ is in an infinite indexing set $\Gamma$. We prove this is impossible by showing that, for every $j\neq k$ in $\Gamma$, there exists $(1,x_{ 1},...,x_{n}) \in R\ltimes_n M$ such that  $(d,m_{j,1},...,m_{j,n})=(1,x_{ 1},...,x_{n})(d,m_{k,1},...,m_{k,n})$. A  recursive argument shows that the fact that every equation $d X= b_i$, with $b_i \in M_i$ admits a solution $X\in M_i$ implies the existence of the desired   $(1,x_{ 1},...,x_{n})$. Note that, from Lemma \ref{lem-U-FFR},  $M_i$ is torsion-free for every   $i\in  \{1,..., n\}$.  First, consider an element $b_n \in M_n- \{ 0\}$. For every $j\in  \Gamma$, $(d,m_{j,1},...,m_{j,n})(0,...,0, b_n)=(0,...,0,db_n) $. Then   $(0,...,0,db_n)= (d,m_{j,1},...,m_{j,n}) \lceil (0,...,0, b_n)  \rceil$ is the only possible corresponding $U$-factorization of $(0,...,0,db_n)$   (since $R\ltimes_n M$ is a $U$-FFR), so there exists $ r \in R$ such that $ b_n= drb_n$. This shows that the above equation admits a solution for $i=n$.   Now consider $k\in \{1,..., n-1\}$ and any  $b_k \in M_k $. For every $b_{n-k}\in M_{n-k}- \{ 0\}$,  $b_k b_{n-k}\in M_{n}- \{ 0\}$ and so there is $r\in R$ such that $b_k b_{n-k}= drb_k b_{n-k}$. Then  $(b_k  -drb_k )b_{n-k}=0$. Therefore $b_k = drb_k$ (since  $M_n$ is  integral).\smallskip

\indent $4.$ Let $i\in \{1,..., n-1\}$. Suppose, by  contradiction, there exists   $m  \in M_i- \{ 0\}$, 
for which there is a family of distinct principal ideals of the form $\langle(0,...,0,m,m_{j,i+1},...,m_{j,n})\rangle$ where $j$ is in an infinite indexing set $\Gamma$. Let $m_{ n-i} \in M_{ n-i}- \{ 0\}$. Necessarily, $m m_{ n-i} \neq 0$. Then $$(0,...,0,m,m_{j,i+1},...,m_{j,n})   (0,...,0,m_{ n-i},0...,0)=(0,...,0, m m_{ n-i}).$$  Then   $(0,...,0, m m_{ n-i})= (0,...,0,m,m_{j,i+1},...,m_{j,n})  \lceil  (0,...,0,m_{ n-i},0...,0)  \rceil$ is the only possible corresponding $U$-factorization of $(0,...,0, m m_{ n-i})$   (since $R\ltimes_n M$ is a $U$-FFR), so there exists  $(r_0,r_{ 1},...,r_{n})$       such that  $$(0,...,0,m,m_{j,i+1},...,m_{j,n})  (r_0,r_{ 1},...,r_{n})(0,...,0,m_{ n-i},0...,0) =(0,...,0,m_{ n-i},0...,0) ,$$  equivalently $(0,...,0,r_0 m m_{ n-i}) =(0,...,0,m_{ n-i},0...,0) $, which is absurd.\cqfd

A ring $S$ is called a \textit{$U$-bounded factorization} ring ($U$-BFR)  if, for each nonzero nonunit $x \in S$, there is a natural number $N(x)$ so that, for any factorization $x = a   \lceil b_{1} \cdots b_{t} \rceil$,    we have $t \leq N(x)$. An $S$-module $H$ is said to be  a \textit{$U$-BF} module if for every  $h  \in  H-\{0\}$  there exists a natural number $N(h)$ so that if $Sh=Sd_1\cdots d_t h'$ where $d_j \not \in U(S)$, $t> N(h)$ and   $h'\in H$,   then, after cancellation and reordering of some of the $d_j$'s  we have  $Sh=Sd_1\cdots d_s h'$ for some $s\leq N(h)$.\smallskip

The question of when   the classical  trivial extension is a $U$-BFR is   still open. However, there is an answer to this question for an integral domain $D$ \cite[Theorem 4.4]{Ax}:  For a $D$-module $N$, $D\ltimes  N$ is a $U$-BFR if and only if $D$ is a BFD and $N$ is a $U$-BF $R$-module. Two more general results for the direct implication were established in \cite[Theorem 3.7  and Lemma 3.8]{AxFRS}. Here, we extend these results to our context. For this we need to introduce the following definition.

\begin{defn}\label{def-U-BFM} Assume $n\geq 2$ and consider $i\in \{1,...,n\}$. The $R$-module $M_i$ is said to be  a $\varphi$-$U$-BF module    (or simply a $U$-BF module)  if, for
  every nonzero element   $x\in  M_i$, there exists a natural number $N(x)$ so that if $Rx=  Rd_{i_1}d_{i_2} \cdots d_{i_t}$ where  $t \in \N $, $(d_{i_1},  ..., d_{i_t})\in M_{i_1} \times \cdots \times M_{i_t}$ for some $\{i_1,i_2 ,...,  i_t\}\subseteq  \{0,...,n\}$ with  $i_1+\cdots + i_t=i$,       $d_{i_j}\not\in U(R)$ when $i_j=0$, and $t> N(x)$,   then, after cancellation and reordering of some of the $d_{i_j}$'s in $R$, we have  $Rx=R d_{i_1}d_{i_2} \cdots d_{i_s}$ for some $s\leq N(h)$.
\end{defn}

 \begin{thm}\label{thm-U-BFR}  If $R\ltimes_n M$ is   a $U$-BFR, then  $R $ is   a $U$-BFR and $M_i$ is a $U$-BF module for every  $i\in \{1,..., n\}$.   Moreover,  if $R$ is pr\'{e}simplifiable, then  $R $ is   a BFR.\\
\indent Conversely,  assume $R$ to be an integral domain. If $R$ is a BFD and  for every  $i\in \{1,..., n\}$, $M_i$ is a $U$-BF module, then  $R\ltimes_n M$ is   a $U$-BFR.
\end{thm}
\pr Similar to the classical case.\cqfd

A ring  $S$ is called \textit{$U$-atomic} if every nonzero nonunit element of $S$ has a $U$-factorization in which all the  relevant  factors are irreducibles. The question of when the classical trivial extension is $U$-atomic is   still unsolved. In \cite[Theorem 4.6]{Ax}, Axtell gave an answer to this question for   an integral domain $D$ with ACCP: For a $D$-module $N$, $D\ltimes  N$ is atomic if and only if  $D\ltimes  N$ is $U$-atomic.  In \cite[Theorem 3.15]{AxFRS}, it is shown that the condition that the ring is an integral domain  could be replaced by the ring is  pr\'{e}simplifiable. The following result gives an extension of \cite[Theorem 3.15]{AxFRS} to the context of $n$-trivial extensions.

\begin{thm}\label{thm-U-atomic}
 Assume $n\geq 2$. Suppose that $M_i$ is is pr\'{e}simplifiable  for every $i\in \{0,...,n-1\}$ (here $M_0=R$), $R$ satisfies ACCP and $M_i$ satisfies ACC on cyclic submodules  for every $i\in \{1, ... , n-1\}$. Then $R\ltimes_n M$ is atomic if and only if  $R\ltimes_n M$ is $U$-atomic.
\end{thm}
\pr $\Longrightarrow $ Clear.\\
$\Longleftarrow  $ Suppose that $R\ltimes_n M$ is not atomic. Then by the proof of Theorem \ref{thm-atomic},    there exists $\underline{m_n}:=(0,... ,0,m_n)  \in R\ltimes_n M$  with $m_n\neq 0$ which is not a product of irreducibles.  Since $R\ltimes_n M$ is $U$-atomic, $\underline{m_n}$ admits a $U$-factorization of the form  $\underline{m_n}  = a_{1}\cdots a_{s} \lceil b_{1} \cdots b_{t} \rceil$ such that the $ b_{l}$'s are irreducibles. Since $\underline{m_n}$ cannot be a product of irreducibles and  by the proof of Theorem \ref{thm-atomic}, necessarily     $s=1$ and $a_1$ has the form $\underline{x_n}:=(0,... ,0,x_n)$. But $\underline{x_n}  \langle b_{1} \cdots b_{t}\rangle = \langle b_{1} \cdots b_{t}\rangle$  and so $ b_{1} \cdots b_{t}$ has the form $\underline{y_n}:=(0,... ,0,y_n)$. This is impossible since $ \underline{m_n} = \underline{x_n}\underline{y_n}=0$.\cqfd


\end{document}